\documentclass[a4paper,11pt]{article}

\usepackage[a4paper,margin=30mm]{geometry}
\usepackage[subsec]{style_unfolding}
\title{Unfolding of symmetric monoidal $(\infty,n)$-categories}
\author{Rune Haugseng and Thomas Nikolaus}

\date{\today}
\begin{document}

\maketitle

\begin{abstract}
We prove Lurie’s proposed description of left rigid symmetric monoidal $(\infty,n)$-cate\-gories in terms of purely $\infty$-categorical data. More precisely, unfolding gives an equivalence between left rigid symmetric monoidal $(\infty,n)$-categories and complete categorical chain complexes of rigid symmetric monoidal $\infty$-categories. More generally, unfolding identifies skeletal sequences with categorical chain complexes without the completeness condition. We extend these equivalences to $(\infty,\infty)$-categories and, stably, obtain an equivalence between complete unbounded categorical chain complexes and left rigid categorical spectra.

We also analyze their groupoid completion and express the homotopy groups of the associated spaces and spectra in terms of categorical cycles and boundaries.

As examples, we construct the categorical chain complexes associated to geometric and algebraic bordism categories, including their equivariant variants, and to iterated span categories. Their localizations recover classical bordism and $L$-groups, while the span examples are related to algebraic $K$-theory.
\end{abstract}

\tableofcontents

\section{Introduction}

The singular homology of a space $X$ is obtained from singular chains by taking cycles modulo boundaries. Bordism homology has the same format: one takes closed manifolds equipped with maps to $X$, modulo those that are boundaries  of a manifold one dimension higher with a map to  $X$. This analogy underlies Atiyah's interpretation of bordism as a generalized homology theory \cite{AtiyahBordism} and the geometric models of Buoncristiano--Rourke--Sanderson \cite{BRS}. This raises the natural question of whether the resulting bordism spectra can themselves be modeled by chain complexes.

In the most na\"ive sense this is not possible: chain complexes of abelian groups model $H\mathbb Z$-module spectra, whereas bordism spectra do not generally admit $H\mathbb Z$-module structures. Nevertheless, replacing abelian groups by abelian monoids, Kochman constructed a chain model for bordism groups \cite{KochmanChain}. Quinn later developed a machine that produces spaces and spectra from bordism-type data \cite{QuinnAssembly}, and Laures--McClure axiomatized this construction in terms of \emph{ad theories}. In this framework, an ordinary chain complex $C$ is a basic example whose associated bordism groups are $H_{*}(C)$, while classical bordism and $L$-theory arise from other examples \cite{LauresMcClure}. \\

We take a different, higher-categorical perspective on bordism spectra:
Manifolds and bordisms in all nonnegative dimensions can be organized into a symmetric monoidal $(\infty,\infty)$-category whose objects are $0$-manifolds, with $k$-dimensional bordisms as $k$-morphisms for all $k$ \cites{LurieCob,CalaqueScheimbauer}. The geometric realization of this full higher bordism category is the grouplike $E_{\infty}$-space corresponding to the connective bordism spectrum
\cite{SPinv}.

Instead of realizing the full higher category, we can also extract its
adjacent-dimensional $\infty$-categorical layers. For every $k\geq0$,
let $\Bord_{k,k+1}$ denote the symmetric monoidal $\infty$-category
of closed $k$-manifolds and $(k+1)$-dimensional bordisms. These
\icats{} will play the role of cycles. For $k\geq1$, the
corresponding \icat{} of chains is the symmetric monoidal
\icat{} $\Bord^\partial_{k,k+1}$ of $k$-manifolds with
boundary and bordisms with corners. Taking boundaries defines a
symmetric monoidal functor
\[
  \partial\colon
  \Bord^\partial_{k,k+1}\longrightarrow\Bord_{k-1,k}.
\]
Its fibre over the empty manifold is precisely $\Bord_{k,k+1}$,
giving a pullback square
\[
  \begin{tikzcd}
    \Bord_{k,k+1}\ar[r]\ar[d]\drpullback
      & \Bord^\partial_{k,k+1}\ar[d,"\partial"]\\
    0\ar[r] & \Bord_{k-1,k}.
  \end{tikzcd}
\]
As $k$ varies, these squares form the categorical analogue of a chain
complex, with $\Bord^\partial_{k,k+1}$ as the category of $k$-chains
and $\Bord_{k,k+1}$ as the category of $k$-cycles.

The same construction applies to the algebraic bordism categories central to Grothendieck--Witt and $L$-theory \cites{CalmesEtAlI,CalmesEtAlII,CalmesEtAl}. One difference is that, in the purely algebraic setting, nothing prevents the ``dimensions'' of such objects from being negative, even though the resulting negative-dimensional objects cannot arise geometrically. The corresponding algebraic bordism categories therefore naturally form both a categorical spectrum and an associated unbounded categorical chain complex. This leads to the central question:
\begin{quote}
Can higher bordism categories (and, stably, categorical spectra) be completely encoded by categorical chain complexes?
\end{quote}
In particular, does the resulting chain data retain enough information to reconstruct the original higher category? Motivated by this picture, Lurie formulated a precise equivalence in his account of the cobordism hypothesis \cite{LurieCob}*{\S3.3}, which we prove here. More precisely, we show that \emph{left rigid} symmetric monoidal \incats{} are equivalent to \emph{complete} categorical chain complexes of rigid symmetric monoidal \icats{}. Stably, this identifies left rigid categorical spectra with complete categorical chain complexes indexed by all integers. This provides a functorial mechanism for constructing bordism-type higher categories and spectra from geometric and algebraic data.

\subsection*{Main results}

To give precise statements, let us fix some notation. 
A \emph{categorical chain complex of length $\ell$} is a sequence of
pullback squares of rigid symmetric monoidal \icats{}
\begin{equation}\label{eq:intro categorical chain complex}
  \begin{tikzcd}
    \oZ_{k} \drpullback \ar[r] \ar[d]
      & \oC_{k} \ar[d,"\partial_{k}"] \\
    0 \ar[r] & \oZ_{k-1},
  \end{tikzcd}
  \qquad 1\leq k\leq\ell,
\end{equation}
where the underlying functor
of $\partial_{k}$ is a cocartesian fibration.  It is called complete if the space of objects
in $\oC_{k}$ is equivalent to the space of objects in $(\oZ_{k-1})_{\bbone/}$.

Our first main result says that these complexes are exactly the
unfolded forms of left rigid higher categories.
\begin{theoremint}\label{introthm:complete}
  Unfolding induces an equivalence of \icats{} between left rigid
  symmetric monoidal \incats{} and complete categorical chain
  complexes of length $n-1$.
\end{theoremint}
Here $\nC$ is \emph{left $k$-rigid} if it has duals for objects (for
$k\geq1$), every morphism $\bbone\to x$ has a left adjoint (for
$k\geq2$), and $\mOm\nC = \nC(\bbone,\bbone)$ is left $(k-1)$-rigid; it is \emph{left
  rigid} if it is left $n$-rigid.  This is weaker than full
dualizability, since only the morphisms detected recursively from the
tensor unit are required to have adjoints.
The finite complete comparison is proved in
\cref{cor:finite complete unfolding}.

Concretely, the unfolding of $\nC$ is the sequence of squares
\begin{equation}\label{eq:intro unfolding squares}
  \begin{tikzcd}
    \core\mOm^{k}\nC \drpullback \ar[r] \ar[d]
      & \core\mP\mOm^{k-1}\nC \ar[d] \\
    0 \ar[r] & \core\mOm^{k-1}\nC,
  \end{tikzcd}
  \qquad 1\leq k\leq n-1.
\end{equation}
Here $\core\nD$ denotes the underlying \icat{} of $\nD$ and the map $\mP\nD\to\core\nD$ is the
cocartesian fibration classifying $\nD(\bbone,\blank)$.  The theorem
asserts both that
\cref{eq:intro unfolding squares} is a categorical chain complex and
that $\nC$, together with all of its compositions, can be reconstructed
from this purely \icatl{} data.  This proves the comparison stated by
Lurie in his account of the cobordism hypothesis
\cite{LurieCob}*{\S3.3}.

We extend these comparisons beyond the results stated in
\cite{LurieCob} in two directions. First, allowing infinite length
identifies left
rigid symmetric monoidal $(\infty,\infty)$-categories with complete
nonnegatively graded categorical chain complexes. Second, stabilizing
with respect to the functor $\mOm$ of endomorphisms of the tensor unit
gives versions of $(\infty,\infty)$-spectra and extends
unfolding to categorical chain complexes that are unbounded in both
directions. The corresponding equivalences are
\cref{cor:unfolding infinite} and \cref{cor:unfolding unbdd}.

Inverting all morphisms in a left rigid symmetric monoidal
$(\infty,\infty)$-category $\nC$ produces a grouplike
$E_{\infty}$-space $\Lgpd{\nC}$, or equivalently a connective
spectrum.  Similarly, localizing an $(\infty,\infty)$-spectrum
produces a possibly non-connective spectrum.  We analyze these
localizations from the perspective of categorical chain complexes and
obtain homotopy-group formulas reminiscent of the homology of an
ordinary chain complex. In the complete case the boundary contribution
vanishes after groupoid completion, giving the following particularly
simple formula. We state the connective case here, leaving the
non-connective case (which requires introducing more notation) to the
body of the paper.

\begin{theoremint}\label{introthm:homotopy-groups}
Let $\nC$ be a left rigid symmetric monoidal \iicat{}. Then its groupoid
completion $\Lgpd{\nC}$ is a grouplike $E_{\infty}$-space, and
\[
  \pi_{i}\Lgpd{\nC}
    \cong\pi_{0}\Lgpd{\core\mOm^{i}\nC}
  \qquad(i\geq0).
\]
\end{theoremint}

This is \cref{cor:complete loc isos}, part~(ii).
The corresponding formulas for arbitrary categorical chain complexes,
which involve a nontrivial boundary term, are given in
\cref{cor:pi of completion of locn}.

\subsection*{Examples and applications}

The bordism complex described above is our basic geometric example.
It need not be complete, but its completion corresponds to a left rigid symmetric
monoidal \iicat{} $\Bord_{0,\infty}$, and the
homotopy-group formula gives
\[
  \pi_{d}\Lgpd{\Bord_{0,\infty}}
    \cong\pi_{0}\Lgpd{\Bord_{d,d+1}}
  \qquad(d\geq0), 
\]
(\cref{cor:geometric bordism localization}).
Thus its realization recovers the usual bordism groups, compatibly
with the identification of the spectra associated to
geometric bordism \iicats{} due to Schommer-Pries \cite{SPinv}.

The same construction has an algebraic counterpart. A closed oriented manifold determines by its singular chains a finite chain complex equipped with a chain-level Poincaré duality structure. 
Replacing manifolds by such Poincaré complexes and geometric bordisms by algebraic cobordisms gives algebraic bordism categories \cite{CalmesEtAlII}.
There are also algebraic variants of manifolds with boundary (so-called metabolic Poincaré complexes) and their cobordism categories, so that we obtain a categorical chain complex. In contrast to the geometric examples, these algebraic examples can be shifted to all integral dimensions so that the unbounded and stable forms of the unfolding equivalence are essential.

More generally, this construction can be formulated for a Poincaré \icat{} $(\mathcal C,\Qoppa)$, where $\mathcal C$ is a stable \icat{} of chain-complex-like objects and $\Qoppa$ encodes a notion of forms and duality on such objects \cite{CalmesEtAlI}. 
The shifted algebraic cobordism categories and their metabolic analogues form an unbounded categorical chain complex, which corresponds under the unfolding equivalence to the categorical spectrum
$
  \operatorname{Cob}_{-\infty,\infty}
    (\mathcal C,\Qoppa).
$
Its localization is canonically equivalent to the $L$-theory spectrum:
\[
  L_{\Sp}\operatorname{Cob}_{-\infty,\infty}
    (\mathcal C,\Qoppa)
    \simeq \mathrm L(\mathcal C,\Qoppa).
\]
In particular, its $i$th homotopy group is the $i$th $L$-group
$
  L_i(\mathcal C, \Qoppa)
$
(\cref{propn:algebraic bordism realization,cor:algebraic bordism homotopy}).
We develop this construction and recall the necessary background on Poincaré \icats{} in \cref{sec:algebraic-bordisms}.

Both the geometric and algebraic constructions admit equivariant
refinements, contravariantly functorial in $\Glo$, which denotes
the full $(2,1)$-subcategory of groupoids spanned by the one-object
groupoids $BG$ for discrete groups $G$. Geometrically, we use manifolds equipped
with actual smooth $G$-actions; algebraically, we use Poincar\'e
objects equipped with coherent $G$-actions by isometries. 
For $G=\mathbb Z$, the oriented variant of the geometric theory is related
to Kreck's bordism groups of orientation-preserving diffeomorphisms. The algebraic
theory is related by
Ranicki's symmetric automorphism signature to automorphism
$L$-theory \cites{KreckDiffeomorphisms,RanickiKnot}.

Finally, iterated span categories give another important example. If $\oE$ has finite limits, then
$\Span_{n}(\oE)$ is fully dualizable, and its unfolded categorical
chain complex has terms
\[
  \oZ_{i}\simeq\Span_{1}(\oE)\quad(i\geq0),
  \qquad
  \oC_{i}\simeq\Span_{1}(\oE^{[1]})\quad(i\geq1),
\]
with $\oC_0\simeq\oZ_0$.
For a stable \icat{} $\oE$, the finite stages recover shifted copies
of its algebraic $K$-groups. These examples, together with the
equivariant geometric and algebraic constructions, are developed in
\cref{sec:examples}.

\subsection*{Outline of the proof}

We prove \cref{introthm:complete} as the complete case of a more
general comparison.  A \emph{skeletal sequence of length $n$} is a
sequence of symmetric monoidal functors
\[
  \nA_{1}\longrightarrow\nA_{2}\longrightarrow\cdots
    \longrightarrow\nA_{n}
\]
such that $\nA_{k}$ is a left rigid symmetric monoidal
$(\infty,k)$-category and $\nA_{k}\to\nA_{k+1}$ is
$(k-1)$-surjective.  The successive categories need not be the cores
of a single higher category and may retain additional flagging data.
The unrestricted form of the unfolding theorem, stated by Lurie in
\cite{LurieCob}*{3.3.33},\footnote{We also prove the variant
  for Lurie's weak skeletal sequences, in which the top level of
  dualizability is omitted; see \cref{obs:weak chain cx}.} is:
\begin{theoremint}\label{introthm:skseq}
  Unfolding induces an equivalence between skeletal sequences of
  length $n$ and categorical chain complexes of length $n-1$.
\end{theoremint}
This is proved as \cref{cor:skseq cch eq}.
Completeness corresponds to requiring each $\nA_{k}$ to be the
$k$-core of $\nA_{k+1}$; restricting \cref{introthm:skseq} to these
towers of cores gives \cref{introthm:complete}.

The proof begins with one-step reconstruction.  If $\nC$ has duals
for objects, the cocartesian fibration
\[
  \mP\nC\longrightarrow\core\nC
\]
classifies the functor $\nC(\bbone,\blank)$.  Duality gives
\[
  \nC(x,y)\simeq
  \nC(\bbone,x^{\vee}\otimes y),
\]
so maps out of the tensor unit determine all maps in 
$\nC$.  Applying this functor to the canonical closed enrichment of
$\core\nC$ reconstructs $\nC$; this is an instance of the enriched
reconstruction theorem of \cite{closed}, recalled in
\cref{subsec:closed}.

To iterate this reconstruction we must understand duals in
$\mP\nC$.  We show that it has duals exactly when $\nC$ is left
$2$-rigid, with the dual of $f\colon\bbone\to x$ given by
$(f^{L})^{\vee}$.  A direct iteration introduces $\mP^{2}\nC$, which
generally has no duals and therefore obstructs induction.

The key two-step argument identifies the part involving
$\mP^{2}\nC$ as a left Kan extension of the corresponding diagram for
$\mP\mOm\nC$.  Consequently, $\nC$ is reconstructed from
\[
  \begin{tikzcd}
    \mP\mOm\nC \ar[d] \\
    \core\mOm\nC \ar[r] \ar[d]
      & \core\mP\nC \ar[d] \\
    0 \ar[r] & \core\nC.
  \end{tikzcd}
\]
Here $\mP\mOm\nC$ has duals exactly when $\nC$ is left $3$-rigid.
Iterating the construction for $\mOm\nC,\mOm^{2}\nC,\ldots$ yields
\cref{eq:intro unfolding squares}; its relative version proves
\cref{introthm:skseq} for skeletal sequences.

The infinite and unbounded comparisons are obtained by passing to
compatible limits and stabilizing with respect to $\mOm$.
Localization is computed by groupoid-completing the terms of the
categorical chain complex and then completing the resulting complex.
This gives the general quotient formulas of
\cref{cor:pi of completion of locn}; for complete complexes the
boundary contribution vanishes, giving
\cref{introthm:homotopy-groups}.

\subsection*{Related work}
The main results of this paper were obtained a long 
time ago and were presented in several talks and circulated. We
apologise for the delay in making the arguments available in written
form. Related results on unfolding monoidal
$(\infty,n)$-categories have also been obtained by Adam Dauser, which will appear  in forthcoming work. 
While this manuscript
was being finalized, Heine posted
\cite{HeineAdjoints}, which also establishes our main results with
similar proofs, including an equivalence between skeletal sequences
of categorical spectra with adjoints and categorical chain
complexes. We were not aware of this work before its appearance.

\subsection*{Organization}

In \cref{sec:background} we review fibrations of higher categories
and the reconstruction theorem for closed symmetric monoidal
enriched categories.  In \cref{sec:unfolding} we establish the
one-step and two-step reconstruction results and prove the finite
unfolding equivalences.  In \cref{sec:extend} we pass to infinite and
unbounded complexes, construct completion, and compute the homotopy
groups after localization.  Finally, \cref{sec:examples} develops
the geometric and algebraic bordism complexes, their equivariant
versions, and the examples coming from iterated spans.

\subsection*{Notation}

\begin{itemize}
\item We denote generic \icats{} by $\oA, \oB, \oC, \dots$ and
  $(\infty,n)$-categories by $\nA, \nB, \nC, \dots$.
\item The \icat{} of \incats{} is denoted
  $\CatIn$.
\item $\SMCatIl$ is the \icat{} of symmetric
  monoidal \icats{} and lax symmetric monoidal functors.
\item We will often write $0$ for the unique symmetric monoidal
  structure on the point.
\item We denote the underlying \icat{} of an \incat{} $\nC$ as $\core
  \nC$, and more generally the underlying $(\infty,k)$-category (or ``$k$-core'') as
  $\core[k] \nC$.
\end{itemize}

\subsection*{Acknowledgments}
We thank Claudia Scheimbauer and Tashi Walde for helpful discussions,
particularly concerning dualizability and coherent chain complexes. We also thank Thorger Geiss for many related discussions and comments on a draft.

The second author was funded by the Deutsche Forschungsgemeinschaft (DFG, German Research
Foundation) -- Project-ID 427320536 -- SFB 1442, as well as under Germany’s Excellence Strategy
EXC2044/2–390685587, Mathematics M\"unster: Dynamics–Geometry–Structure.

\section{Background}
\label{sec:background}
Our goal in this section is to review the main result of
\cite{closed}, specialized to closed symmetric monoidal \incats{}, and
reformulate this in terms of fibrations --- this will provide the main
input to the unfolding equivalence we prove in the next section. We
start by reviewing the particular cases of fibrations of \incats{} and
their straightening that are relevant for us in \cref{sec:fib}, and
use this to describe symmetric monoidal \incats{} fibrationally in
\cref{sec:symmon}. In \cref{subsec:closed} we then review the main
result of \cite{closed} in the case of closed symmetric monoidal
\incats{} and then specialize this to the \emph{rigid} case, where objects
have duals, in \cref{sec:rigid}.

\subsection{Fibrations and straightening}
\label{sec:fib}
In this subsection we recall the definition of fibrations from an
\incat{} to an \icat{} and their straightening, and rephrase this as
itself being an equivalence of fibrations.

\begin{defn}
  Let $\oB$ be an \icat{} and $\eE$ an $(\infty,n)$-category, and
  consider a functor $p \colon \eE \to \oB$. A morphism
  $f \colon x \to y$ is \emph{$p$-cocartesian} if the commutative
  square of $(\infty,n-1)$-categories
  \[
    \begin{tikzcd}
      \eE(y, z) \ar[r, "f^{*}"] \ar[d] & \eE(x,z) \ar[d] \\
      \oB(py,pz) \ar[r, "p(f)^{*}"] & \oB(px,pz)
    \end{tikzcd}
  \]
  is a pullback for all $z \in \eE$; dually $f$ is \emph{$p$-cartesian}
  if it is $p^{\op}$-cocartesian. We say $p$ is a \emph{cocartesian
    fibration} if for every $x \in \eE$ and $f \colon p(x) \to b$ in
  $\oB$ there exists a $p$-cocartesian morphism
  $\bar{f} \colon x \to y$ such that $p(\bar{f}) \simeq f$; dually,
  $p$ is a \emph{cartesian fibration} if $p^{\op}$ is a cocartesian
  fibration. For an \icat{} $\oB$, we write $\CatInsl{\oB}^{\coc}$ for
  the subcategory of $\CatInsl{\oB}$ whose objects are cocartesian
  fibrations over $\oB$ and whose morphisms are functors that preserve
  cocartesian morphisms.
\end{defn}

\begin{warning}
  This is only the correct definition of a cocartesian fibration when
  the target is an \icat{}, which is the only case that will appear in
  this paper; in general there are further conditions involving the
  functors $\eE(x,y) \to \oB(px,py)$ but these are automatic (for all
  choices of variances) when $\oB(px,py)$ is an \igpd{}.
\end{warning}

\begin{thm}[Nuiten~\cite{NuitenStr}]\label{thm:str fib incat}
  For $\oB$ an \icat{}, there is a natural equivalence of \icats{}
  \[\CatInsl{\oB}^{\coc} \simeq \Fun(\oB, \CatIn)\]
  where the (contravariant) functoriality in $\oB$ is with respect to pullback on the left and composition on the right.
\end{thm}

\begin{remark}
  In fact, Nuiten shows that this is an equivalence of
  $(\infty,n+1)$-categories, with $\oB$ itself allowed to be an
  $(\infty,n)$-category, but we only require the given special case. A
  simpler proof of this particular statement can also be obtained by
  extending the argument of
  Blans and Blom~\cite[Proposition A.2.9]{BlansBlom} in the case $n=2$.
\end{remark}

Altogether these equivalences combine to a natural equivalence of
functors of \icats{} $\CatI \to \LCatI$, which we can in turn unstraighten to get
an equivalence of fibrations using the following notation:
\begin{notation}
  We write $\Fib_{n/1}^{\coc}$ for the subcategory of $\Ar(\CatIn)$
  whose
  \begin{itemize}
  \item objects are functors $p \colon \eE \to \oB$ with $\oB$ an
    \icat{} such that $p$ is a cocartesian fibration,
  \item morphisms are commuting squares where the top horizontal
    morphism preserves cocartesian morphisms.
  \end{itemize}
  Since $\ev_{1} \colon \Ar(\CatIn) \to \CatIn$ is the cocartesian
  fibration for the overcategory functor $\CatInsl{\blank} \colon
  \CatIn \to \LCatI$ with functoriality given by composition, it is
  also the cartesian fibration for the right adjoint functors given by
  pullback. Its restriction to the subcategory $\Fib_{n/1}^{\coc}$
  hence gives the cartesian fibration for \[\CatInsl{\blank}^{\coc}
    \colon \CatI^{\op} \to \LCatI \]
  via pullback.
\end{notation}

\begin{notation}
  For $\tB$ an \itcat{} and $b \in \tB$, we write $\tB_{\upslash b}
  \to \tB$ for the cartesian fibration for the mapping \icat{} functor
  $\tB(\blank,b)$; this can also be defined as the fibre at $b$ of the functor
  \[\ev_{1} \colon \FUN([1], \tB)^{\oplax} \to \tB\]
  from the oplax arrow \itcat{} of $\tB$. 
  If $\oB$ is the underlying \icat{} of $\tB$, we
  will abusively write $\oB_{\upslash b}$ for the underlying \icat{}
  of $\tB_{\upslash b}$ when the $(\infty,2)$-upgrade of $\oB$
  intended is clear; then the forgetful functor \[\oB_{\upslash b} \to
    \oB\]
  is the cartesian fibration for the functor
  \[ \tB(\blank,b) \colon \oB^{\op} \to \CatI. \]
  Informally, the objects of
  $\oB_{\upslash b}$ are morphisms with target $b$ in $\oB$, and
  morphisms are oplax triangles
  \[
    \begin{tikzcd}
     x \ar[rr] \ar[dr, ""{name=A}] & & y \ar[dl] \\
     & b.
     \arrow[from=A,to=1-3,Rightarrow,shorten <=8pt, shorten >=20pt]
    \end{tikzcd}
  \]
\end{notation}

We can then reformulate \cref{thm:str fib incat} in terms of
fibrations as:
\begin{cornoproof}\label{cor:str as fib eqce}
  There is an equivalence of cartesian fibrations over $\CatIn$,
  \[
    \begin{tikzcd}
     \Fib_{n/1}^{\coc} \ar[rr, "\sim"] \ar[dr, "\ev_{1}"'] & & \CatIlsl{\CatIn} \ar[dl] \\
      & \CatIn,
    \end{tikzcd}
  \]
  where $\CatIlsl{\CatIn}$ denotes the full subcategory of
  $\LCatIlsl{\CatIn}$ containing functors from small \icats{} to
  $\CatIn$.
\end{cornoproof}

\begin{remark}\label{rmk:str as fib 0}
  The equivalence of \cref{cor:str as fib eqce} makes sense also in
  the case $ n =0$, if we interpret $\Fib_{0/1}^{\coc}$ as the full
  subcategory of $\Ar(\CatI)$ consisting of left fibrations.
\end{remark}

\subsection{Symmetric monoidal $(\infty,n)$-categories and
  straightening}
\label{sec:symmon}

\begin{defn}
  If $\oO$ is an \iopd{}, an \emph{$\uO$-monoidal \incat{}} is a cocartesian
  fibration $\eC^{\otimes} \to \uO$ such that the corresponding
  functor $\oO \to \CatIn$ is an $\oO$-monoid. A \emph{lax
    $\oO$-monoidal functor} is a commutative triangle
  \[
    \begin{tikzcd}
     \eC^{\otimes} \ar[rr, "F"] \ar[dr, "p"'] & & \eD^{\otimes} \ar[dl, "q"] \\
      & \oO
    \end{tikzcd}
  \]
  where $p$ and $q$ are $\oO$-monoidal \incats{} and $F$ preserves
  cocartesian morphisms over inert maps in $\oO$; if $F$ preserves
  \emph{all} cocartesian morphisms it is a (\emph{strong}) \emph{$\oO$-monoidal functor}.
  For $\xF_{*}$, we instead use the term \emph{symmetric monoidal}. We
  write $\SMCatInl$ for the subcategory of $\CatInsl{\xF_{*}}$ with
  objects the symmetric monoidal \incats{} and morphisms the lax
  symmetric monoidal functors, and $\SMCatIn$ for its wide subcategory
  containing only the strong symmetric monoidal functors.
\end{defn}

\begin{notation}
  We write $\SMFibl_{n/1}$ for the full subcategory of the pullback
  \[ \Fib_{n/1}^{\coc} \times_{\CatI} \SMCatIl,\]
  taken via the target and source functors, that is spanned by the
  objects
  \[ \nE \xto{p} \oB^{\otimes} \to \xF_{*}\]
  such that the cocartesian fibration $p$ exhibits $\nE$ as a
  $\oB^{\otimes}$-monoidal \incat{}.
\end{notation}

\begin{observation}
  A morphism in $\SMFibl_{n/1}$ is a
  commutative diagram
  \begin{equation}
    \label{eq:smfib mor}
    \begin{tikzcd}
      \nE^{\otimes} \ar[rr, "F"] \ar[d, "p"] & & \nF^{\otimes} \ar[d, "q"] \\
      \oB^{\otimes} \ar[rr, "G"] \ar[dr] & & \oC^{\otimes} \ar[dl] \\
      & \xF_{*},
    \end{tikzcd}
  \end{equation}
  where $G$ preserves inert cocartesian morphisms over $\xF_{*}$ and
  $F$ takes \emph{all} $p$-cocartesian morphisms to $q$-cocartesian ones.
\end{observation}

With this notation, we get the following consequence of \cref{cor:str as fib eqce}:
\begin{cor}\label{cor:fibeqsmn}
  There is an equivalence of cartesian fibrations over $\SMCatl$
  \[
    \begin{tikzcd}
     \SMFibl_{n/1} \ar[rr, "\sim"] \ar[dr] & & \SMCatIllsl{\CatIn^{\times}} \ar[dl] \\
      & \SMCatIl,
    \end{tikzcd}
  \]
  where $\SMCatIllsl{\CatIn^{\times}}$ denotes the full subcategory of
  $\LSMCatIllsl{\CatIn^{\times}}$ spanned by lax symmetric monoidal
  functors to $\CatIn^{\times}$ from small symmetric monoidal \icats{}.
\end{cor}
\begin{proof}
  Pulling back the equivalence of \cref{cor:str as fib eqce}, we get
  an equivalence of cartesian fibrations over $\SMCatl$ between
  $\SMFibl_{n/1}$ and the full subcategory of
  \[ \CatIlsl{\CatIn} \times_{\CatI} \SMCatl\]
  spanned by the functors $\oB^{\otimes} \to \CatIn$ that are
  $\oB^{\otimes}$-monoids, for all symmetric monoidal \icats{}
  $\oB^{\otimes}$. This is the fibration for the functor
  $(\SMCatl)^{\op} \to \LCatI$ taking $\oB^{\otimes}$ to the full subcategory
  \[\Mon_{\oB^{\otimes}}(\CatIn) \subset \Fun(\oB^{\otimes}, \CatIn)\]
  of $\oB^{\otimes}$-monoids. By \cite[Proposition 2.4.2.5]{HA} we
  have a natural equivalence
  \[ \Mon_{(\blank)}(\CatIn) \simeq \Alg_{(\blank)}(\CatIn^{\times})\]
  for the cartesian symmetric monoidal structure on $\CatIn$; the
  right-hand side gives the mapping \icat{} functor for the \itcat{}
  of \iopds{}, or its full subcategory of symmetric monoidal \icats{}
  and lax symmetric monoidal functors, so the corresponding cartesian fibration is
  \[ \SMCatIllsl{\CatIn^{\times}} \to \SMCatIl,\]
  as required.
\end{proof}

\begin{observation}\label{obs:Bmon via symmon}
  Suppose $\oB^{\otimes} \to \xF_{*}$ is a symmetric monoidal
  \icat{} and we have a commutative triangle
  \[
    \begin{tikzcd}
     \nE \ar[rr, "F"] \ar[dr, "p"'] & & \oB^{\otimes} \ar[dl] \\
      & \xF_{*}
    \end{tikzcd}
  \]
  where $\nE$ is an \incat{}.  Then $F$ exhibits $\nE$ as a
  $\oB^{\otimes}$-monoidal \incat{} \IFF{} $p$ exhibits $\nE$ as a
  symmetric monoidal \incat{},
  $F_{\angled{1}} \colon \nE_{\angled{1}} \to \oB$ is a cocartesian
  fibration, and the tensor product functor preserves
  $F_{\angled{1}}$-cocartesian morphisms in each variable.
  We can therefore identify $\SMFibl_{n/1}$ as a subcategory of the
  arrow \icat{}  $\Ar(\SMCatInl)$:
  \begin{itemize}
  \item Its objects are the lax symmetric monoidal functors
    \[
      \begin{tikzcd}
       \nE^{\otimes} \ar[rr, "\phi"] \ar[dr] & & \oB^{\otimes} \ar[dl] \\
        & \xF_{*}
      \end{tikzcd}
    \]
    such that $\oB^{\otimes}$ is an \icat{}, $\phi_{\angled{1}} \colon
    \nE \to \oB$ is a cocartesian fibration, and the tensor product in $\nE$
    preserves $\phi_{\angled{1}}$-cocartesian morphisms.
  \item Its morphisms are the commutative diagrams of the form
    \cref{eq:smfib mor} where the top horizontal map $F$ preserves all
    cocartesian morphisms over the middle horizontal map $G$.
  \end{itemize}
\end{observation}

\subsection{Closed symmetric monoidal $(\infty,n)$-categories}\label{subsec:closed}

The main result of \cite{closed} gave a description of closed
symmetric monoidal $\oV$-enriched \icats{} in terms of lax symmetric
monoidal functors to $\oV$. Here we recall the specialization of this
result to the case where $\oV$ is $\CatIn$, and then reformulate this
in terms of fibrations using \cref{cor:fibeqsmn}. We first introduce
some notation:

\begin{defn}
  We say that a symmetric monoidal \incat{} $\eC$ is \emph{closed} if
  the functor of \incats{} $x \otimes \blank \colon \eC \to \eC$ has a
  right adjoint for all $x \in \eC$, or equivalently the presheaf of
  \inmcats{} $\eC(\blank \otimes x, y)$ is
  representable by an object $[x,y]$ for all $x,y \in \eC$. We write
  $\SMCatInlcl$ for the full subcategory of $\SMCatInl$ spanned by the
  closed symmetric monoidal \incats{}.
\end{defn}

\begin{observation}
  Suppose $\nC$ is a symmetric monoidal \incat{}. If the underlying
  \icat{} $\core \nC$ is closed, then for all objects $x,y$ there
  exists an object $[x,y]$ and a counit map $x \otimes [x,y] \to
  y$. It is then a property for composition with this map to induce an
  equivalence
  \[ \nC(c, [x,y]) \to \nC(x \otimes c,x \otimes [x,y]) \to \nC(x
    \otimes c, y)\] of \inmcats{} (rather than just on
  underlying \igpds{}). Thus $\nC$ is closed \IFF{} $\core \nC$ is closed
  and the counit maps have this additional property.
\end{observation}

\begin{defn}
  A \emph{pinning} of a
  closed symmetric monoidal \incat{} $\nC$ is a closed symmetric
  monoidal \icat{} $\oB$ together with a strong symmetric
  monoidal functor $\oB \to \core \nC$ that is essentially surjective
  on objects and
  preserves internal homs. We write
  \[ \PSMCatInlcl \subseteq \Ar(\SMCatInlcl) \] for the full
  subcategory whose objects are the \emph{pinned} closed symmetric
  monoidal \incats{}, \ie{} closed symmetric monoidal \incats{}
  equipped with a pinning.
\end{defn}

\begin{thm}[\cite{closed}*{Theorem A}]\label{thm:closedsmenr}
  There is an
  equivalence 
  \[ \SMCatIlcllsl{\CatInm^{\times}} \simeq \PSMCatInlcl \] of
  \icats{} over
  $\SMCatIlcl$. Furthermore, this restricts to identify $\SMCatInlcl$ 
  with the full subcategory of $\SMCatIlcllsl{\CatInm^{\times}}$ whose
  objects are lax symmetric monoidal functors
  $\phi \colon \oC \to \CatInm$ such that the canonical map
  $\oC(\bbone_{\oC}, \blank) \to \core[0] \phi(\blank)$ is
  an equivalence.
\end{thm}

\begin{remark}
  Here the objects of the \icat{} $\SMCatIlcllsl{\CatInm^{\times}}$ are closed
  symmetric monoidal \icats{} equipped with a lax symmetric monoidal
  functor to $\CatInm$, and a morphism $(\oC, \phi) \to (\oD, \psi)$
  consists of a lax symmetric monoidal functor $f \colon \oC \to \oD$
  together with a natural transformation $\phi \to \psi \circ f$. The
  lax symmetric monoidal functor $\phi \colon \oC \to \CatInm$
  corresponds to the \incat{} $\phi_{*}\overline{\oC}$, where
  $\overline{\oC}$ is the self-enrichment of $\oC$, which we push
  forward along the lax symmetric monoidal functor $\phi$; this is
  pinned by a canonical strong symmetric monoidal functor
  $\oC \to U(\phi_{*}\overline{\oC})$; in the other direction,
  $\oC \to \nD$ corresponds to the composite
  \[\oC \to \core \nD \xto{\nD(\bbone_{\nD}, \blank)} \CatInm.\]
\end{remark}

\begin{notation}
  Let \[\SMFiblcl_{n/1} := \SMFibl_{n/1} \times_{\SMCatIl} \SMCatIlcl\]
  denote the full subcategory of $\SMFibl_{n/1}$ whose projection to
  $\SMCatIl$ is a closed symmetric monoidal \icat{}.
\end{notation}

Combining the equivalences of \cref{thm:closedsmenr} with
\cref{cor:fibeqsmn}, we then have:
\begin{cornoproof}\label{cor:fibcatneq1}
  There is an equivalence
  \[
    \SMFiblcl_{n-1/1} \simeq \PSMCatInlcl
  \]
  over $\SMCatIlcl$.
\end{cornoproof}
We next want to identify the full subcategory of $\SMFiblcl_{n-1/1}$
that corresponds to $\SMCatInlcl$ under this equivalence. For this we need some
preliminary observations:

\begin{construction}
  Suppose $\oB^{\otimes}$ is a symmetric monoidal \icat{}. Then the
  forgetful functor $(\oB^{\otimes})_{\angled{0}/} \to \oB^{\otimes}$
  is a $\oB^{\otimes}$-monoidal \icat{}: It is a left fibration, and
  so in particular a cocartesian
  fibration, and the Segal condition follows from the commutative
  square
  \[
    \begin{tikzcd}
     \Map_{\oB^{\otimes}}(\angled{0}, (b_{1},\ldots,b_{n})) \ar[r] \ar[d] & \prod_{i=1}^{n}\Map_{\oB^{\otimes}}(\angled{0}, b_{i})  \ar[d] \\
     \Hom_{\xF_{*}}(\angled{0}, \angled{n}) \ar[r, "\sim"'] & \prod_{i=1}^{n}\Hom_{\xF_{*}}(\angled{0},\angled{1})
    \end{tikzcd}
  \]
  being a pullback, since $\angled{0}$ is initial in $\xF_{*}$. The
  composite $(\oB^{\otimes})_{\angled{0}/} \to \oB^{\otimes} \to
  \xF_{*}$ is a symmetric monoidal structure on $\oB_{\bbone/}$, given
  on objects by \[(\bbone \xto{f} b) \otimes (\bbone \xto{g} b') \simeq
  \bbone \xto{f \otimes g} b \otimes b',\] and we therefore denote it
  \[ \oB_{\bbone/}^{\otimes} := (\oB^{\otimes})_{\angled{0}/}.\]
\end{construction}

\begin{observation}\label{obs:Bunder1initial}
  By \cite[Theorem 1.2]{freepres}, the forgetful functor
  $\oB_{\bbone/}^{\otimes} \to \oB^{\otimes}$ is the free cocartesian
  (and indeed left) fibration on $\{\angled{0}\} \to \oB^{\otimes}$,
  so if $\eE^{\otimes} \to \oB^{\otimes}$ is a
  $\oB^{\otimes}$-monoidal \incat{}, then the unique object of
  $\nE_{\angled{0}} \simeq *$ corresponds to a unique
  $\oB^{\otimes}$-monoidal functor
  $\oB_{\bbone/}^{\otimes} \to \nE^{\otimes}$. Straightening over
  $\oB^{\otimes}$, the left fibration
  $\oB_{\bbone/}^{\otimes} \to \oB^{\otimes}$ corresponds to a lax
  symmetric monoidal structure on the functor
  $\oB(\bbone, \blank) \colon \oB \to \Spc$, so we can conclude that
  \begin{itemize}
  \item the functor $\oB(\bbone,\blank)$ has a unique lax symmetric
    monoidal structure,
  \item and any lax symmetric monoidal functor $F$ from $\oB$ to $\Spc$ (or
    more generally to $\CatIn$) admits a unique symmetric monoidal
    natural transformation $\oB(\bbone, \blank) \to
    F$.\footnote{The underlying natural transformation corresponds
      via the Yoneda lemma to the map $* \to F(\bbone)$ given by the
      lax symmetric monoidal structure on $F$.}
  \end{itemize}
  In other words, $\oB(\bbone,\blank)$ is the initial lax symmetric
  monoidal functor $\oB \to \Spc$.
\end{observation}

\begin{notation}
  Let $\CSMFiblcl_{n/1}$ denote the full subcategory of
  $\SMFiblcl_{n/1}$ that is spanned by the objects
  $\nE^{\otimes} \to \oB^{\otimes}$ such that the canonical map
  $\oB_{\bbone/}^{\otimes} \to \nE^{\otimes}$ restricts to an
  equivalence on the underlying left fibration of $\nE$, or
  equivalently gives an equivalence on \igpd{} cores; we refer to
  these as the \emph{complete} objects of $\SMFiblcl_{n/1}$.
\end{notation}

\begin{cor}\label{cor:fibcatneq}
  The equivalence of \cref{cor:fibcatneq1} restricts to 
  an equivalence
  \[ \CSMFiblcl_{n-1/1} \simeq \SMCatInlcl\]
  on full subcategories.
\end{cor}
\begin{proof}
  We know the equivalence in question restricts to one between
  $\SMCatInlcl$ and the
  full sub-\itcat{} of $\SMFiblcl_{n-1/1}$ spanned by objects
  \[ \nE^{\otimes} \xto{p} \oB^{\otimes} \to \xF_{*}\]
  where $p$ corresponds to a lax symmetric monoidal functor $\phi
  \colon \oB \to
  \CatIn$ such that the canonical natural transformation
  $\oB(\bbone,\blank) \to \phi(\blank)^{\simeq}$ is an
  equivalence. By the uniqueness in \cref{obs:Bunder1initial}, this corresponds precisely
  to the unique map $\oB_{\bbone/}^{\otimes} \to \nE^{\otimes}$ giving an
  equivalence on underlying left fibrations.
\end{proof}

\begin{observation}
  For $n = 1$, we can identify $\CSMFiblcl_{0/1}$ (interpreted
  appropriately as in \cref{rmk:str as fib 0}) as the \icat{} of
  left fibrations $\oE^{\otimes} \to
  \oB^{\otimes}$ where $\oB_{\bbone/}^{\otimes} \isoto \oE^{\otimes}$;
  this is evidently equivalent to $\SMCatIlcl$.
  For $n = 2$ we have identified closed
  symmetric monoidal $(\infty,2)$-categories with the data of a closed
  symmetric monoidal \icat{} $\oB^{\otimes}$ and a
  $\oB^{\otimes}$-monoidal \icat{} $\eE^{\otimes}$ such that
  $\oB_{\bbone/} \to \oB$ is the underlying left fibration of $\eE$.
\end{observation}

We can restrict the equivalences above to strong symmetric monoidal
functors, for which we introduce the following notation:
\begin{notation}
  We write $\SMFib_{n/1}$ for the wide subcategory of $\SMFibl_{n/1}$
  where the morphisms are the commutative diagrams
  \[
    \begin{tikzcd}
      \nE^{\otimes} \ar[rr] \ar[d] & & \nF^{\otimes} \ar[d] \\
      \oB^{\otimes} \ar[rr] \ar[dr] & & \oC^{\otimes} \ar[dl] \\
       & \xF_{*},
    \end{tikzcd}
  \]
  where \emph{both} horizontal functors preserve \emph{all} cocartesian
  morphisms. Similarly, $\SMFibcl_{n/1}$ and $\CSMFibcl_{n/1}$ denote
  the intersection of $\SMFib_{n/1}$ with $\SMFiblcl_{n/1}$ and
  $\CSMFiblcl_{n/1}$, respectively. We also write $\PSMCatIncl$ for
  the wide subcategory of $\PSMCatInlcl$ where the morphisms are
  strong symmetric monoidal for both the symmetric monoidal \incat{}
  and its pinning and $\SMCatIncl$ for the full subcategory of
  $\SMCatIn$ spanned by the closed symmetric monoidal \incats{}.
\end{notation}

\begin{cor}
  The equivalence of \cref{cor:fibcatneq1} restricts on subcategories
  to equivalences
  \[ \SMFibcl_{n-1/1} \simeq \PSMCatIncl, \qquad \CSMFibcl_{n-1/1}
    \simeq \SMCatIncl.\]
\end{cor}
\begin{proof}
  This follows from \cite{closed}*{Corollary 6.11}, interpreted in the fibrational context.
\end{proof}

\subsection{The rigid case}
\label{sec:rigid}

In this subsection we briefly recall the notion of duals for objects
in a symmetric monoidal \incat{}, and specialize the equivalences of
the previous section to this case, which will be our concern for
the rest of this paper.

\begin{defn}
  Recall that if $\oC$ is a symmetric monoidal \icat{}, a \emph{dual}
  of an object $x \in \oC$ is another object $x^{\vee}$ together with
  maps
  \[ \ev : x \otimes x^{\vee} \to \bbone, \quad \coev \colon \bbone
    \to x^{\vee} \otimes x,\] called the \emph{evaluation} and
  \emph{coevaluation}, so that the composites
  \[ x \xto{\id \otimes \coev} x \otimes x^{\vee} \otimes x \xto{\ev
      \otimes \id} x,\]
  \[ x^{\vee} \xto{\coev \otimes \id} x^{\vee} \otimes x \otimes
    x^{\vee} \xto{\id \otimes \ev} x^{\vee},\] are both equivalent to
  the corresponding identity map; this data is unique if it exists by
  \cite[Lemma 4.6.1.10]{HA}. We say that $\oC$ is a \emph{rigid}
  symmetric monoidal \icat{} if every object of $\oC$ has a dual. If
  $\nC$ is a symmetric monoidal \incat{}, we say that $\nC$ is
  \emph{1-rigid} if the underlying symmetric monoidal \icat{}
  $\core \nC$ is rigid.
\end{defn}

\begin{observation}
  If $\oC$ is a symmetric monoidal \icat{} and $x \in \oC$ has a dual
  $x^{\vee}$, then the evaluation and coevaluation maps exhibit $x
  \otimes \blank$ as left adjoint to $x^{\vee} \otimes
  \blank$. If $\oC$ is rigid, it is therefore in particular closed,
  with internal Homs given by
  \[ [x,y] = x^{\vee} \otimes y.\]
  Moreover, if $\nC$ is a 1-rigid symmetric monoidal \incat{}, then
  the evaluation and coevaluation for $x$ similarly give equivalences
  \[ \nC(x \otimes y, z) \simeq \nC(y, x^{\vee} \otimes z),\]
  so that $\nC$ is in fact a closed symmetric monoidal \incat{}.
\end{observation}

\begin{observation}
  Suppose $F \colon \oC \to \oD$ is a strong symmetric monoidal
  functor of \icats{}, and $x \in \oC$ has a dual $x^{\vee}$. Then
  $F(\ev)$ and $F(\coev)$ exhibit $F(x^{\vee})$ as the dual of
  $F(x)$. In particular, strong symmetric monoidal functors preserve
  dualizable objects.
\end{observation}

\begin{observation}\label{obs:rigid gpd}
  Suppose $X$ is a symmetric monoidal \igpd{}. Then an object of $X$
  is dualizable \IFF{} it is invertible. Thus $X$ is rigid \IFF{} it
  is \emph{grouplike} (or a commutative group in $\Spc$).
\end{observation}

\begin{notation}
  We write
  \begin{itemize}
  \item   $\SMFibpltrig_{n/1}$ for the full sub-\itcat{} of
    $\SMFibpl_{n/1}$ spanned by objects
  \[ \nE^{\otimes} \to \oB^{\otimes} \to \xF_{*}\]
  where $\oB$ is rigid,
\item $\CSMFibpltrig_{n/1}$ for the intersection of
  $\SMFibpltrig_{n/1}$
  with $\CSMFibplcl_{n/1}$,
\item $\SMCatInplorig$ for the
  full sub-\itcat{} of $\SMCatInplcl$ on the 1-rigid symmetric monoidal \incats{},
\item $\PSMCatInplorig$ for the
  full sub-\itcat{} of $\PSMCatInplcl$ on the pinned
  symmetric monoidal \incats{} $\oB \to \nC$ such that both $\oB$ and
  $\nC$ have duals for objects.
  \end{itemize}
\end{notation}

\begin{observation}
  Since any strong symmetric monoidal functor preserves duals, we can
  describe $\PSMCatInplorig$ as the full subcategory of
  $\Ar(\SMCatInpl)$
  spanned by the strong symmetric monoidal functors $F \colon \oB \to \nC$
  where $\oB$ is an \icat{}, $F$ is essentially surjective on
  objects, and $\oB$ is rigid (since then $\nC$ is also 1-rigid).
\end{observation}

\begin{cornoproof}\label{cor:fibeqduals}
  The equivalence of \cref{cor:fibcatneq1} restricts to equivalences
  \[ \SMFibpltrig_{n-1/1} \simeq \PSMCatInplorig,\quad \CSMFibpltrig_{n-1/1} \simeq
    \SMCatInplorig\]
  on subcategories.
\end{cornoproof}

\section{Unfolding for $(\infty,n)$-categories}
\label{sec:unfolding}

Our goal in this section is to prove \cref{introthm:skseq} (and its
special case \cref{introthm:complete}). We start by identifying duals
in \icats{} of the form $\mP\nC$ in \cref{sec:dual lax slice}. This
allows us to unpack left 2-rigid \incats{} to a square, which we carry
out in \cref{sec:unpack square}. We next show in \cref{sec:two-step}
that these squares are obtained by left Kan extension, using a
technical result that we prove in \cref{sec:kan ext result}. This is
the key input required to formulate the inductive step for skeletal
sequences in \cref{sec:skeletal}. Finally, we define categorical chain
complexes and carry out the induction to complete the proof in
\cref{sec:chain complex unfolding}.

\subsection{Duals in lax slices}\label{sec:dual lax slice}
Let $\nC$ be a closed symmetric monoidal \incat{}. We write
\[ (\mP\nC)^{\otimes} \to (\core\nC)^{\otimes} \]
for the $(\core \nC)^{\otimes}$-monoidal $(\infty,n-1)$-category it
corresponds to under the equivalence of \cref{cor:fibcatneq}.
Here
$(\mP\nC)^{\otimes}$ is a symmetric monoidal $(\infty,n-1)$-category,
so if it also happens to be closed, then we can regard the fibration
as a morphism in $\SMCat_{(\infty,n-1)}^{\orig}$. We can then
apply \cref{cor:fibeqduals} again, and recover this functor from the
corresponding morphism in $\CSMFib_{n-2/1}^{\trig}$, which is the
commutative square
\[
  \begin{tikzcd}
    (\mP^{2} \nC)^{\otimes} \ar[r] \ar[d] & (\core\nC)_{\bbone/}^{\otimes} \ar[d] \\
    (\core\mP\nC)^{\otimes} \ar[r] & (\core \nC)^{\otimes}
  \end{tikzcd}
\]
of symmetric monoidal $(\infty,n-2)$-categories.
Our goal in this section is to characterize precisely when this happens in the rigid case:
when $\nC$ is 1-rigid, we will establish a necessary and sufficient
criterion for $\mP\nC$ to also be 1-rigid. We will also 
identify the squares of symmetric monoidal $(\infty,n-2)$-categories that correspond to (pinned)
symmetric monoidal \incats{} in this way.

\begin{notation}
  For a pinned symmetric monoidal \incat{} $\oA \to \nC$, we write
  $\mP\nC|_{\oA}$ for the pullback $\oA \times_{\core \nC}
  \mP\nC$. Then the corresponding fibration is $\mP\nC|_{\oA} \to \oA$.
\end{notation}

\begin{observation}
  The underlying functor $\mP\nC \to \core\nC$ is the cocartesian
  fibration for 
  \[ \nC(\bbone, \blank) \colon \core\nC \to \Cat_{(\infty,n-1)}.\] We
  can describe $\mP\nC$ in low dimensions as follows:
  \begin{itemize}
  \item The objects of $\mP\nC$ are morphisms $\bbone \to x$ in $\nC$,
    projecting to $x$ in $\core\nC$.
  \item The morphisms of $\mP\nC$ are triangles
    \[
      \begin{tikzcd}
        {} & \bbone \ar[dl, ""{name=A}] \ar[dr, ""{swap, name=B}] \\
        x \ar[rr] & & y \ar[from=A, to=B, Rightarrow]
      \end{tikzcd}
    \]
    containing a 2-morphism in $\nC$.
  \item The 2-morphisms
    \[
      \begin{tikzcd}
        {} & \bbone \ar[dl, "f"', ""{name=A}] \ar[dr, "g", ""{swap, name=B}] \\
        x \ar[rr, "h"] & & y \ar[from=A, to=B, Rightarrow, "\phi"]
      \end{tikzcd} \to \begin{tikzcd}
        {} & \bbone \ar[dl, "f"', ""{name=A}] \ar[dr, "g", ""{swap, name=B}] \\
        x \ar[rr, "h'"] & & y \ar[from=A, to=B, Rightarrow, "\psi"]
      \end{tikzcd}
    \]
    consist of an equivalence $h\simeq h'$ and (after identifying the
    sources of $\phi$ and $\psi$ via this equivalence) a 3-morphism
    $\phi \to \psi$ in $\nC$.
  \end{itemize}
\end{observation}

\begin{defn}
  We say a symmetric monoidal \incat{} $\nC$ is \emph{left 2-rigid} if
  \begin{itemize}
  \item $\nC$ is 1-rigid,
  \item and every morphism in $\nC$ of the form $f \colon \bbone \to X$ has a left
    adjoint $f^{L} \colon X \to \bbone$.
  \end{itemize}
\end{defn}

Every rigid symmetric monoidal \igpd{} (or grouplike commutative
monoid in spaces, \cref{obs:rigid gpd}) is also left 2-rigid, as the
additional condition is vacuous. However, these are the only
symmetric monoidal \icats{} that are left 2-rigid:
\begin{propn}\label{propn:left 2-rigid cat is gpd}
  Suppose $\oC$ is a symmetric monoidal \icat{} that is left 2-rigid
  when regarded as an \itcat{}. Then $\oC$ is an \igpd{}.
\end{propn}
\begin{proof}
  Since $\oC$ is left 2-rigid, we know that every morphism $\bbone \to X$ in
  $\oC$ has a left adjoint; since $\oC$ is an \icat{}, all such
  morphisms are therefore invertible. In
  particular, the coevaluation $\bbone \to X^{\vee} \otimes X$ is
  invertible for every $X$. Given a morphism $f \colon A
  \to B$ such that $\id_{X} \otimes f$ is invertible for some object $X$, it follows that so is
  $f$, by applying the 2-of-3 property to the square
  \[
    \begin{tikzcd}
     A \ar[r, "\mathrm{coev}_{X} \otimes \id_{A}", "\sim"'] \ar[d,
     "f"'] & X^{\vee}\otimes X \otimes A  \ar[d, "\id_{X^{\vee}}
     \otimes \id_{X} \otimes f", "\sim"'] \\
     B \ar[r, "\mathrm{coev}_{X} \otimes \id_{B}"', "\sim"] & X^{\vee}\otimes X \otimes B.
    \end{tikzcd}
  \]
  For an 
  arbitrary morphism $f \colon A \to B$, we can now consider the
  composite 
  \[ \bbone \xto{\coev_{A}} A^{\vee} \otimes A \xto{\id_{A^{\vee}} \otimes f}
    A^{\vee} \otimes B, \]
  where we know that the first map and the composite are
  invertible. It then follows by 2-out-of-3 that the second map is invertible, and
  hence by what we just observed that $f$ itself is also
  invertible. Since $f$ was arbitrary, this completes the proof.
\end{proof}

\begin{propn}[\cite{LurieCob}*{3.3.29}]\label{propn:PCdualadj}
  Suppose $\nC$ is a 1-rigid symmetric monoidal \incat{}. Then the symmetric monoidal
  \inmcat{}
  $\mP \nC$ has duals for objects \IFF{} $\nC$ is left 2-rigid; the
  dual of $f \colon \bbone \to X$ in $\mP\nC$ is
  $(f^{L})^{\vee} \colon \bbone \to X^{\vee}$.
\end{propn}

\begin{notation}
  Suppose $x$ is a dualizable object in $\nC$. Then the functor
  \[ T_{x} \colon \nC(x,y) \to \nC(\bbone, x^{\vee} \otimes y)\]
  defined by
  \[T_{x}(f) = \bbone \xto{\coev_{x}} x^{\vee} \otimes x \xto{\id
      \otimes f} x^{\vee} \otimes y\]
  is an equivalence, with inverse given by
  \[ T_{x}^{-1}(g) = x \xto{\id_{x} \otimes g} x \otimes x^{\vee} \otimes
    y \xto{\ev_{x} \otimes \id_{y}} y.\]
  For $f \colon x \to y$ we abbreviate $f^{t} := T_{x}(f)$ (and call this the
  \emph{transpose} of $f$) and (when $y$ is also dualizable) $f^{\vee}
  := T_{y^{\vee}}^{-1}T_{x}(f)\colon y^{\vee} \to x^{\vee}$ (and call
  this the \emph{dual} of $f$); note that if $y = \bbone$ then $f^{t}
  = f^{\vee}$. 
\end{notation}

\begin{lemma}\label{lem:dual comp eq}
  Suppose $x$ and $y$ are dualizable objects of an \incat{}
  $\nC$, and let \[\gamma \colon \nC(\bbone,x^{\vee} \otimes y) \times \nC(\bbone, y^{\vee}
    \otimes z) \to \nC(\bbone, x^{\vee} \otimes z)\]
  be the functor given by
  \[
    \gamma(f,g) = (\id_{x^{\vee}} \otimes \ev_{y} \otimes \id_{z})(f
    \otimes g). \]Then there is a commutative diagram
  \[
    \begin{tikzcd}
      \nC(x,y) \times \nC(y,z) \ar[r, "T_{x} \times T_{y}", "\sim"'] \ar[d, "\circ"'] &
     \nC(\bbone,x^{\vee} \otimes y) \times \nC(\bbone, y^{\vee}
     \otimes z)  \ar[d, "\gamma"] \\
     \nC(x,z) \ar[r, "T_{x}"', "\sim"] & \nC(\bbone,x^{\vee}\otimes z).
    \end{tikzcd}
  \]
  In other words, $T_{x}(gf) \simeq \gamma(T_{x}f, T_{y}g)$, and
  similarly for higher morphisms.
\end{lemma}
\begin{proof}
  This follows from the triangle identities.
\end{proof}

\begin{observation}\label{obs:dual maps comp}
  For maps $f \colon \bbone \to x$ and $g
  \colon x \to \bbone$ we have in particular $\ev_{x} \circ (f \otimes g^{\vee})
  \simeq gf$ and $(\ev_{x} \otimes \id_{x}) \circ (\id_{x} \otimes g^{\vee}
  \otimes f) \simeq g \otimes f \simeq fg$. 
\end{observation}

\begin{proof}[Proof of \cref{propn:PCdualadj}]
  First suppose $\mP\nC$ is 1-rigid. Since the forgetful functor to
  $\core \nC$ is symmetric monoidal, and so preserves duality data, if
  $f \colon \bbone \to x$ has a dual $g$, then $g$ is necessarily a
  map $\bbone \to x^{\vee}$, and the coevaluation and evaluation maps
  are triangles of the form
  \[
    \begin{tikzcd}
      {} & \bbone \ar[dl, "="', ""{name=A}] \ar[dr, "g \otimes f", ""{swap, name=B}] \\
      \bbone \ar[rr, "\coev_{x}"] & & x^{\vee} \otimes x, \ar[from=A, to=B, Rightarrow, "\phi"]
    \end{tikzcd} \quad \begin{tikzcd}
      {} & \bbone \ar[dl, "f \otimes g"', ""{name=A}] \ar[dr, "=", ""{swap, name=B}] \\
      x \otimes x^{\vee} \ar[rr, "\ev"] & & \bbone. \ar[from=A, to=B, Rightarrow, "\epsilon"]
    \end{tikzcd}
  \]
  Here \cref{obs:dual maps comp} identifies $\epsilon$ as a 2-morphism
  $g^{\vee}f \to \id$, while from $\phi \colon \coev \to g \otimes f$
  we get a 2-morphism
  $\eta := T_{x}^{-1}(\phi) \colon \id \to fg^{\vee}$. We must show
  that the triangle identities for the duality imply that these are
  the unit and counit of an adjunction between $f$ and $g^{\vee}$.
  This follows from \cref{lem:dual comp eq} since we have
  \[(\ev_{x} \otimes \id_{x})\circ (f \otimes \phi) \simeq f\eta \colon
    f \to fg^{\vee}f,\]
  \[ \epsilon \otimes f \simeq \epsilon f \colon fg^{\vee}f \to f,\]
\[ (\id_{x^{\vee}} \otimes \ev_{x}) \circ (\phi \otimes g) \simeq
  \eta^{\vee} g \colon g \to gf^{\vee}g,\]
\[ g \otimes \epsilon \simeq g \epsilon^{\vee} \colon gf^{\vee}g \to g.\]
The converse is proved in the same way, applying these identifications
in reverse.
\end{proof}

\begin{notation}
  Let $\PSMCatIn^{2\dlrig}$ denote the full subcategory of
  $\PSMCatIn^{\rig}$ whose objects are the pinned symmetric
  monoidal $(\infty,n)$-categories that are left 2-rigid, and define
  $\SMCatIn^{2\dlrig}$ similarly. We also write
  $\SMFibstrig_{n/1}$ and $\CSMFibstrig_{n/1}$ for the full subcategories of
  $\SMFibtrig_{n/1}$ and $\CSMFibtrig_{n/1}$, respectively spanned by the
  objects whose source also has duals.
\end{notation}
\begin{cor}\label{cor:pin2lrig}
  The equivalences of \cref{cor:fibeqduals} restrict on full
  subcategories to equivalences
  \[ \PSMCatIntworig \simeq \SMFibstrig_{n-1/1}, \quad
    \SMCatIntworig \simeq \CSMFibstrig_{n-1/1}. \]
\end{cor}
\begin{proof}
  By construction, an object $\oB \xto{\phi} \nC$ of $\PSMCatIntworig$
  corresponds to the fibration $\mP\nC|_{\oB} := \oB \times_{\core \nC} \mP\nC \to
  \oB$. We know from \cref{propn:PCdualadj} that $\nC$ is left 2-rigid
  \IFF{} $\mP\nC$ is 1-rigid, so it only remains to show that the
  latter is equivalent to the pullback $\mP\nC|_{\oB}$
  being 1-rigid. It follows from \cite[Proposition 4.6.1.11]{HA} that
  an object in the pullback is dualizable \IFF{} its image in $\oB$
  and $\mP\nC$ is dualizable (as then the image in $\core \nC$ is
  automatically dualizable, since the functors in the pullback are
  strong symmetric monoidal); the conditions are therefore equivalent
  since $\phi$ is by assumption essentially surjective on objects.
\end{proof}

\begin{remark}
  The version of \cref{cor:pin2lrig} with lax symmetric monoidal
  functors is also true, by the same proof, but we omit this case from
  now on for simplicity, as we will in any case have to drop it when
  we get to the inductive unfolding process below.
\end{remark}

\subsection{Unpacking to a square}
\label{sec:unpack square}

In \cref{cor:pin2lrig}, an object of $\SMFibstrig_{n-1/1}$ can also be viewed as an arrow in $\SMCat_{(\infty,n-1)}^{\orig}$ (whose
target happens to be an \icat{}), and we can identify the former with a
full subcategory of $\Ar(\SMCat_{(\infty,n-1)}^{\orig})$. We
can thus apply the equivalence of \cref{cor:fibeqduals} again and hence
identify $\PSMCatIntworig$ with some full subcategory of
$\Ar(\SMFibtrig_{n-2/1})$. Our goal in this subsection is to 
identify this subcategory. More generally, we want to consider also a slight
variant of this situation, where we include some pinning data
on the right-hand side of the equivalence. On the left-hand side,
this requires the following definition:
\begin{defn}
  Let $\nC$ be a symmetric monoidal $(\infty,n)$-category that is left
  2-rigid. A
  \emph{2-pinning} of $\nC$ is a sequence of symmetric monoidal
  functors
  \[ \oA \xto{F} \tB \xto{G} \nC\]
  such that
  \begin{itemize}
  \item $\oA$ is a 1-rigid symmetric monoidal \icat{},
  \item $\tB$ is a left 2-rigid symmetric monoidal \itcat{}
  \item $F$ is essentially surjective on objects,
  \item $G$ is essentially surjective on objects and
    1-morphisms.
  \end{itemize}
  The data of $\nC$ together with its 2-pinning we refer to as a
  \emph{2-pinned} symmetric monoidal $(\infty,n)$-category.  Let
  $\PSMCatIn[2]^{2\dlrig}$ denote the full subcategory of
  $\Fun([2], \SMCatIn)$ spanned by the 2-pinned left 2-rigid
  symmetric monoidal $(\infty,n)$-categories.
\end{defn}

\begin{observation}
  Here we can regard $\PSMCatIn[2]^{2\dlrig}$ as a full
  subcategory of $\Ar(\PSMCatIn^{2\dlrig})$ with objects the squares
  \[
    \begin{tikzcd}
     \oA \ar[r, "\alpha"] \ar[d, "F"'] & \oA'  \ar[d, "F'"] \\
     \tB \ar[r, "G"'] & \nC
    \end{tikzcd}
  \]
  such that $\alpha$ is an equivalence and
  $\oA \to \tB \to \nC$ is a 2-pinning.
\end{observation}

\begin{notation}
  Let $\PSMFib_{n/1}$ denote the full subcategory of
  $\Ar(\SMFib_{n/1})$ with objects the squares
    \[
    \begin{tikzcd}
     \oB^{\otimes} \ar[r, "G"] \ar[d] & \nC^{\otimes}  \ar[d] \\
     \oA^{\otimes} \ar[r, "\alpha"'] & \oA'^{\otimes}
    \end{tikzcd}
  \]
  such that $\alpha$ is an equivalence, $\oB$ is an \icat{}, and $G$
  is a pinning of $\nC$. We introduce similar notation for
  subcategories, in particular $\PSMFibstrig_{n/1}$. 
\end{notation}

Combining this observation with \cref{cor:pin2lrig}, we get:
\begin{cor}\label{cor:unf to square step 1}
  The equivalence of \cref{cor:pin2lrig} induces an equivalence
  \[ \PSMCatIn[2]^{2\dlrig} \simeq 
    \PSMFibstrig_{n-1/1}.\]
\end{cor}
\begin{proof}
  It is clear that an object 
      \[
    \begin{tikzcd}
     \oB^{\otimes} \ar[r, "G"] \ar[d] & \nC^{\otimes}  \ar[d] \\
     \oA^{\otimes} \ar[r, "\alpha"'] & \oA'^{\otimes}
    \end{tikzcd}
  \]
  of $\PSMCatIn[2]^{2\dlrig}$
  corresponds to a morphism
  \[
    \begin{tikzcd}
     \oA \ar[r, "\alpha"] \ar[d, "F"'] & \oA'  \ar[d, "F'"] \\
     \nD \ar[r, "\gamma"'] & \nE
    \end{tikzcd}
  \]
  where $\alpha$ is an equivalence and $F$ and $F'$ are pinnings. It
  remains to check that under this equivalence, $\gamma$ is a
  2-pinning \IFF{} $G$ is essentially surjective on objects. Since
  both $F$ and $F'$ are essentially surjective on objects, so is
  $\gamma$. Observe that we can think of an object of $\oB$ as a pair
  $(a, \bbone \xto{\phi} F(a))$ and \[G(a, \phi) \simeq (\alpha(a),
    \gamma(\phi) \colon \bbone \to F'\alpha(a)).\]
  It follows that if $\gamma$ is essentially surjective on morphisms
  then $G$ is so on objects. Conversely, if $G$ is surjective on
  objects then every morphism in $\nE$ of the form $\bbone \to e$ is
  in the image of $\gamma$. Hence so is every morphism since we can
  transpose it to one of this form and $\gamma$ preserves transposes.
\end{proof}

\begin{observation}
  Here the \icat{} $\PSMFibstrig_{n-1/1}$ can also be viewed as a subcategory of
  $\Ar(\PSMCat_{(\infty,n-1)}^{\orig})$, whose objects are the squares
  \[
    \begin{tikzcd}
      \oB \ar[r, "F"] \ar[d, "G"'] & \oA  \ar[d, "p"] \\
      \nC \ar[r, "q"'] & \nD
    \end{tikzcd}
  \]
  where $p$ is an equivalence and $F$ and $q$ are symmetric monoidal
  cocartesian fibrations over $\oA$ with $G$ preserving cocartesian
  morphisms; the morphisms are similarly required to preserve
  cocartesian morphisms.
\end{observation}

The equivalence of
\cref{cor:fibeqduals} therefore identifies $\PSMFibstrig_{n-1/1}$ with some subcategory of
$\Ar(\SMFibtrig_{n-2/1})$. Our next goal is to describe this
subcategory; we note the following in order to state this
characterization:

\begin{observation}\label{obs:duals implies cartesian}
  Suppose $\eE^{\otimes} \to \oB^{\otimes}$ is a
  $\oB^{\otimes}$-monoidal \incat{}, where both $\eE$ and $\oB$ have
  duals for objects. Then the underlying functor $p \colon \eE \to
  \oB$, which is by assumption a cocartesian fibration, is also a
  \emph{cartesian} fibration --- taking duals gives a commutative
  square
  \[
    \begin{tikzcd}
      \eE \ar[r, "\vee"{above}, "\sim"{below}] \ar[d, "p"] & \eE^{\op}
      \ar[d, "p^{\op}"] \\
      \oB \ar[r, "\vee"{below}, "\sim"{above}] & \oB^{\op}
    \end{tikzcd}
  \]
  where $p^{\op}$ is a cartesian fibration. In particular, the
  $p$-cartesian morphisms in $\eE$ are precisely the duals of the
  $p$-cocartesian morphisms.
\end{observation}

\begin{defn}\label{defn:smfibsq}
  Let $\tSMFibSqtrig_{n/1}$ be the
   sub-\itcat{} of
   $\Ar(\SMFibtrig_{n/1})$ whose objects are commutative squares
    \[
      \begin{tikzcd}
        \nC^{\otimes} \ar[r] \ar[d] & \oE^{\otimes} \ar[d]
        \\
        \oB^{\otimes} \ar[r] & \oA^{\otimes}
      \end{tikzcd}
    \]
    where
    \begin{itemize}
    \item $\oB^{\otimes} \to \oA^{\otimes}$ is a  cocartesian
      fibration of symmetric monoidal \icats{},
    \item the canonical symmetric monoidal functor $\oA_{\bbone/} \to
      \oE$ is an equivalence,
    \item $\oB$ and $\oA$ have duals for objects,
    \item $\nC$ is a $\oB^{\otimes}$-monoidal $(\infty,n)$-category,
    \item the functor $\nC^{\otimes} \to \oE^{\otimes}$ preserves cocartesian
    morphisms,
  \end{itemize}
  and whose morphisms are commutative cubes
  \[
    \begin{tikzcd}[row sep=small,column sep=small]
     \nC_{1}^{\otimes} \ar[rr] \ar[dr] \ar[dd] & & \oE_{1}^{\otimes} \ar[dr] \ar[dd] \\
      & \nC_{2}^{\otimes} \ar[rr,crossing over] & & \oE_{2}^{\otimes} \ar[dd] \\
      \oB_{1}^{\otimes} \ar[rr] \ar[dr] & & \oA_{1}^{\otimes}  \ar[dr] \\
      & \oB_{2}^{\otimes}  \ar[rr]  \ar[uu,leftarrow,crossing over] & & \oA_{2}^{\otimes} 
    \end{tikzcd}
  \]
  where the map $\oB_{1} \to \oB_{2}$ preserves cocartesian morphisms
  over $\oA_{i}$ (and by default all functors are strong symmetric
  monoidal and the morphisms in the top square all preserve all
  cocartesian morphisms).
  Let $\SMFibSqtrig_{n/1}$ be the full subcategory of
  $\tSMFibSqtrig_{n/1}$ whose objects satisfy the following
  further condition:
  \begin{itemize}
  \item For every \emph{cartesian}\footnote{Cf.~\cref{obs:duals
      implies cartesian}.} morphism $\bar{f} \colon x \to y$ in $\oB$
  over $f \colon a \to b$ in $\oA$, the commutative square
    \[
      \begin{tikzcd}
        \nC_{x} \ar[r, "\bar{f}_{!}"]
        \ar[d] & \nC_{y} \ar[d] \\
        \oE_{a} \ar[r, "f_{!}"] & \oE_{b}
      \end{tikzcd}
    \]
    is a pullback.
  \end{itemize}
\end{defn}

\begin{propn}\label{propn:P2SMCat2lrig to square}
  Applying \cref{cor:fibeqduals} twice gives equivalences
  \[ \PSMCatIn[2]^{(\lax,)2\dlrig} \simeq \PSMFibstrig_{n-1/1}
    \simeq \SMFibSqtrig_{n-2/1}.\]
\end{propn}
\begin{proof}
  The first equivalence is \cref{cor:unf to square step 1}, so
  it remains to show that the conditions for
  objects and morphisms of \[\PSMFibstrig_{n-1/1} \subseteq
  \Ar(\PSMCat_{(\infty,n-1)}^{\orig})\] correspond to those defining
  \[\SMFibSqtrig_{n-2/1} \subseteq \Ar(\SMFib_{n-2/1}^{\trig}).\]
  Consider an object
  \begin{equation}
    \label{eq:ob of ar psmcat}
    \begin{tikzcd}
      \oB \ar[r, "F"] \ar[d, "G"'] & \oA  \ar[d, "p"] \\
      \nC \ar[r, "q"'] & \nD
    \end{tikzcd}
  \end{equation}
  of $\Ar(\PSMCat_{(\infty,n-1)}^{\orig})$; this corresponds to the object
  \begin{equation}
    \label{eq:ob of ar smfib}
    \begin{tikzcd}
      \mP\nC|_{\oB} \ar[r, "\pi"] \ar[d] & \mP\nD|_{\oA}  \ar[d] \\
      \oB \ar[r, "F"'] & \oA
    \end{tikzcd}
  \end{equation}
  in $\Ar(\SMFib_{n-2/1}^{\trig})$. Here we observe that $p$ is an
  equivalence \IFF{} the canonical symmetric monoidal functor
  $\oA_{\bbone/} \to \mP\nD|_{\oA}$ is an equivalence. We therefore
  need to show that an $F$-cocartesian morphism in $\oB$ is mapped by
  $G$ to a $q$-cocartesian morphism in $\nC$ \IFF{} the pullback
  condition in \cref{defn:smfibsq} holds for $\mP\nC|_{\oB}$. (This
  implies that the tensor product of cocartesian morphisms is
  cocartesian since $G$ is surjective on objects and strong symmetric
  monoidal; the remaining conditions also obviously match up.)

  Given a cocartesian morphism $\bar{f} \colon x \to y$ in $\oB$ over
  $f \colon a \to b$ in $\oA$, we want to show that $G(\bar{f})$ is
  also a cocartesian morphism in $\nC$, meaning that composition with
  $G(\bar{f})$ gives a pullback square
  \[
    \begin{tikzcd}
      \nC(G(y),z) \ar[r, "\bar{f}^{*}"] \ar[d] & \nC(G(x),z)  \ar[d] \\
     \oA(b,c) \ar[r, "f^{*}"'] & \oA(a,c)
    \end{tikzcd}
  \]
  for every object $z$ in $\nC$ over $c \in \oA$. Using duality for
  objects in $\nC$, and using that any $z$ is $G(z')$ for some $z' \in
  \oB$, we can identify the square above as
  \[
    \begin{tikzcd}
      \nC(\bbone,Gy^{\vee} \otimes z) \ar[r, "(G\bar{f}^{\vee} \otimes
      \id)_{*}"] \ar[d] & \nC(\bbone,Gx^{\vee} \otimes z)  \ar[d] \\
     \oA(\bbone,b^{\vee}\otimes c) \ar[r, "(f^{\vee} \otimes
     \id)_{*}"'] & \oA(\bbone,a^{\vee} \otimes c)
   \end{tikzcd}
   \quad = \quad
    \begin{tikzcd}
      (\mP\nC|_{\oB})_{y^{\vee} \otimes z'} \ar[r, "(\bar{f}^{\vee} \otimes
      \id)_{!}"] \ar[d] & (\mP\nC|_{\oB})_{x^{\vee} \otimes z'}  \ar[d] \\
     \oA(\bbone,b^{\vee}\otimes c) \ar[r, "(f^{\vee} \otimes
     \id)_{*}"'] & \oA(\bbone,a^{\vee} \otimes c),
   \end{tikzcd}
  \]
  where $\bar{f}^{\vee} \otimes \id_{z'}$ is cartesian, since it is the
  dual of the cocartesian morphism $\bar{f}
  \otimes\id_{z'^{\vee}}$. The given condition for cartesian morphisms
  therefore implies that this square is a pullback. Conversely, if
  $\nC \to \oA$ is a cocartesian fibration, $G$ preserves cocartesian
  morphisms, and $\bar{f} \colon x \to
  y$ in $\oB$ is cartesian over $f \colon a \to b$, then we similarly
  have an identification of squares
  \[
    \begin{tikzcd}
      (\mP\nC|_{\oB})_{x} \ar[r, "\bar{f}_{!}"] \ar[d] & (\mP\nC|_{\oB})_{y}  \ar[d] \\
     \oA(\bbone,a) \ar[r, "f_{*}"'] & \oA(\bbone, b),
   \end{tikzcd}
   \quad = \quad
    \begin{tikzcd}
      \nC(Gx^{\vee},\bbone) \ar[r, "\bar{f}^{\vee,*}"] \ar[d] & \nC(Gy^{\vee},\bbone)  \ar[d] \\
     \oA(a^{\vee},\bbone) \ar[r, "f^{\vee,*}"'] & \oA(b^{\vee},\bbone)
    \end{tikzcd}   
  \]
  where the right-hand square is a pullback since $\bar{f}^{\vee}$ is a
  cocartesian morphism.  This shows that objects match up under our
  desired equivalences; to see that the conditions on morphisms also agree, we
  just have to observe that the preservation of cocartesian morphisms
  by a functor of $(\infty,n-1)$-categories can be checked on the
  underlying functor of \icats{}.
\end{proof}

\begin{warning}
  Even if $\nC$ has more adjoints, we can't repeat the same process
  again: $\mP \nC$ does not have many adjoints for 1-morphisms, since
  these must map to adjoints in $\core\nC$, \ie{} equivalences. In
  particular, $\mP\nC$ does not usually have adjoints for all maps
  from the monoidal unit, which by \cref{propn:PCdualadj} means that
  $\mP^{2}\nC$ will not have duals.
\end{warning}

\subsection{Two-step unfolding}\label{sec:two-step}
We have seen that if $\nC$ is a left 2-rigid symmetric monoidal \incat{},
then we can recover $\nC$ from the commutative square
\[
  \begin{tikzcd}
   \mP^{2}\nC \ar[r] \ar[d] & (\core\nC)_{\bbone/}  \ar[d] \\
   \core\mP\nC \ar[r] & \core\nC
  \end{tikzcd}
\]
of symmetric monoidal $(\infty,n-2)$-categories. Our goal in this
section is to show that we can reconstruct this square, and thus $\nC$
itself, from the diagram
\[
  \begin{tikzcd}
    \mP\mOm\nC  \ar[d]  \\
    \core\mOm \nC \ar[r] \ar[d] \drpullback & \core\mP\nC \ar[d] \\
     0 \ar[r] & \core\nC
  \end{tikzcd}
\]
where $\mOm \nC := \nC(\bbone, \bbone)$ is the $(\infty,n-1)$-category
of endomorphisms of the unit in $\nC$. This equivalence provides
the inductive step for the ``unfolding'' of $\nC$ to a chain complex
of symmetric monoidal \icats{}, which we consider in the next section.
More precisely, we will obtain a new description of the \icat{}
$\SMFibSqtrig_{n/1}$, starting with the following
observations:
\begin{observation}
  We have a commutative triangle
  \[
    \begin{tikzcd}
      \tSMFibSqtrig_{n/1}
     \ar[rr] \ar[dr] & & \SMFibstrig_{1/1} \ar[dl] \\
      & \SMCatI^{\rig}
    \end{tikzcd}
  \]
  where the top horizontal functor takes a square
  \[
    \begin{tikzcd}
      \nC^{\otimes} \ar[r] \ar[d] & \oA_{\bbone/}^{\otimes} \ar[d]
      \\
      \oB^{\otimes} \ar[r] & \oA^{\otimes}
    \end{tikzcd}
  \]
  to $\oB^{\otimes} \to \oA^{\otimes}$, and the diagonal functors pick
  out $\oA^{\otimes}$. Observe that all 3 functors in this diagram are cartesian fibrations:
  \begin{itemize}
  \item in $\tSMFibSqtrig_{n/1}$ the cartesian morphism over
    $\oA' \to \oA$ is given by the cube
    \[
      \begin{tikzcd}[row sep=small,column sep=small]
        \nC'^{\otimes} \ar[rr] \ar[dr] \ar[dd] & & \oA'^{\otimes}_{\bbone/} \ar[dr] \ar[dd] \\
        & \nC^{\otimes} \ar[rr,crossing over] & & \oA^{\otimes}_{\bbone/} \ar[dd] \\
        \oB'^{\otimes} \ar[rr] \ar[dr] & & \oA'^{\otimes}  \ar[dr] \\
        & \oB^{\otimes}  \ar[rr]  \ar[uu,leftarrow,crossing over] & & \oA^{\otimes} 
      \end{tikzcd}
    \]
    where the top and bottom faces are pullbacks,
  \item in $\SMFibstrig_{1/1}$ the cartesian morphism over $\oA' \to
    \oA$ with target $\oB \to \oA$ is given by pulling back $\oB$ to $\oA'$,
  \item in $\tSMFibSqtrig_{n/1}$ the cartesian morphism over
    \[
      \begin{tikzcd}
        \oB' \ar[r] \ar[d] & \oA'  \ar[d] \\
        \oB \ar[r] & \oA
      \end{tikzcd}
    \]
    in $\SMFibstrig_{1/1}$ with target $\nC$ is given by 
    pulling back $\nC$ along the induced functor \[\oB'
      \times_{\oA'} \oA'_{\bbone/} \to \oB \times_{\oA} \oA_{\bbone/}.\]
  \end{itemize}
  Note that the horizontal functor moreover preserves cartesian morphisms.
  
  Here $\SMCatI^{\rig}$ has $\xF_{*} = *^{\otimes}$ as its initial object, and we
  can identify the map on fibres over this as
  $\ev_{1} \colon \SMFib_{n/1}^{\trig} \to \SMCatI^{\rig}$. We therefore
  have a commutative square
  \[
    \begin{tikzcd}
      \tSMFibSqtrig_{n/1} \ar[r] \ar[d] & \SMFib_{n/1}^{\trig}  \ar[d] \\
     \SMFib_{1/1}^{\strig} \ar[r] & \SMCatI^{\rig},
    \end{tikzcd}
  \]
  with the horizontal maps given by cartesian transport to the fibre
  over $*^{\otimes}$. The top horizontal functor is thus given by
  \[    \left( \begin{tikzcd}
      \nC^{\otimes} \ar[r] \ar[d] & \oA_{\bbone/}^{\otimes} \ar[d]
      \\
      \oB^{\otimes} \ar[r] & \oA^{\otimes}
    \end{tikzcd}
    \right) \mapsto
    \left(\begin{tikzcd}
      (\nC \times_{\oA_{\bbone/}} \oA(\bbone,\bbone))^{\otimes} \ar[d]
      \\\oB_{\bbone}^{\otimes}
    \end{tikzcd}
    \right).
 \]
  Let $\frD_{n}$ denote the pullback, and write
  \[ R \colon \tSMFibSqtrig_{n/1} \to \frD_{n}\]
  for the induced functor over $\SMFib_{1/1}^{\strig}$. Then  $R$
  is a morphism of cartesian fibrations.
\end{observation}

The main result we will prove in this section is then the following description of $\SMFibSqtrig_{n/1}$:
\begin{thm}\label{thm:2-step}
  The functor $R$ has a fully faithful relative left adjoint over
  $\SMFib_{1/1}^{\strig}$, whose image is $\SMFibSqtrig_{n/1}$. In particular, $R$
  restricts to an equivalence $\SMFibSqtrig_{n/1} \isoto \frD_{n}$.
\end{thm}
Combining this equivalence with those of \cref{propn:P2SMCat2lrig to
  square} and \cref{cor:fibeqduals}, we get:
\begin{cornoproof}\label{cor:P2 to CC}
  There are equivalences
  \[
    \begin{split}
      \PSMCatIn[2]^{2\dlrig} & \simeq \SMFib_{n-2/1}^{\trig}
      \times_{\SMCatI^{\rig}} \SMFib_{1/1}^{\strig} \\ &  \simeq
      \PSMCat_{(\infty,n-1)}^{\orig} \times_{\SMCatI^{\rig}} \SMFib_{1/1}^{\strig}
    \end{split}
  \]
  of \icats{}.
\end{cornoproof}

\begin{observation}
  Here an object in the right-hand pullback is a cocartesian fibration
  of symmetric monoidal \icats{} with duals $\oB^{\otimes} \to
  \oA^{\otimes}$ together with a pinned symmetric monoidal
  $(\infty,n-1)$-category $\oB_{\bbone}^{\otimes} \to
  \nC^{\otimes}$. An object $\oA \to \nB \to \nC$ in
  $\PSMCatIn[2]^{2\dlrig}$ corresponds to $\mOm\nB \to \mOm\nC$ in
  $\PSMCat_{(\infty,n-1)}^{\rig}$ together with $\oA \times_{\core
    \nB} \mP\nB \to \oA$ in $\SMFib_{1/1}^{\strig}$.
\end{observation}

\begin{observation}\label{obs:2-rig to complete CC}
  We can identify $\SMCatIn^{2\dlrig}$ with a full subcategory of
  $\PSMCatIn[2]^{2\dlrig}$ by sending $\nC$ to the canonical 2-pinning
  $\core \nC \to \core[2] \nC \to \nC$. Then the equivalence of
  \cref{cor:P2 to CC} restricts to an equivalence
  \[ \SMCatIn^{2\dlrig} \simeq 
    \SMCat_{(\infty,n-1)}^{\rig} \times_{\SMCatI^{\rig}}
    \CSMFib_{1/1}^{\strig},
    \]
  with $\nC \in \SMCatIn^{2\dlrig}$ corresponding to $\mOm \nC$
  together with $\mP\nC \to \core \nC$ and the equivalence
  $(\mP\nC)_{\bbone} \simeq \core \mOm \nC$.
\end{observation}

Since $R$ is a morphism of cartesian fibrations, it follows from
\cite{HA}*{7.3.2.6} that in order to show that $R$ admits a relative left
adjoint, it suffices to show that the induced functor on each fibre
has a left adjoint.

\begin{observation}
  On the fibres over $p \colon \oB^{\otimes} \to \oA^{\otimes}$, the functor
  $R$ is equivalent to the functor
  \[
    R_{p} \colon \Mon_{\oB^{\otimes}}(\Cat_{(\infty,n)})_{/\oA(\bbone,p(\blank))}
    \to \Mon_{\oB_{\bbone}^{\otimes}}(\Cat_{(\infty,n)}) \] given by
  restricting to the fibre $\oB_{\bbone}^{\otimes}$ and then taking
  the pullback along $\{\id_{\bbone}\} \to \oA(\bbone, \bbone)$ (where
  the restriction of $\oA(\bbone,p(\blank))$ is constant at
  $\oA(\bbone, \bbone)$).
\end{observation}

The proof of \cref{thm:2-step} therefore reduces to the following
fibrewise statement, which we will prove in the next section: 
\begin{propn}\label{propn:2-step fibrewise}
  Suppose $p \colon \oB^{\otimes} \to \oA^{\otimes}$ is a cocartesian fibration
  of symmetric monoidal \icats{}, such that the underlying functor is
  also a cartesian fibration. Then (operadic) left Kan extension along
  the inclusion $i \colon \oB_{\bbone}^{\otimes} \to \oB^{\otimes}$ of the
  fibre over the unit $\bbone \in \oA$ induces a fully faithful
  functor
  \[ \Mon_{\oB_{\bbone}^{\otimes}}(\CatIn) \hookrightarrow
    \Mon_{\oB^{\otimes}}(\CatIn)_{/\oA(\bbone,p(\blank))}, \]
  left adjoint to $R_{p}$,
  whose image consists of those monoidal natural transformations that
  satisfy the condition from \cref{cor:fibrelkeff}.
\end{propn}

\subsection{A technical result on Kan extensions}
\label{sec:kan ext result}
In this subsection we will prove \cref{propn:2-step fibrewise}. We
start by looking at ordinary left Kan extensions along functors of the form $\oB_{\bbone} \to
\oB$, which we will later see also compute operadic Kan extensions in
this case. This is a special case of a more general result on left Kan
extensions compatible with factorization systems:

\begin{propn}\label{propn:lkefact}
  Suppose $\oB$ is an \icat{} with a factorization system
  $(\oB_{L},\oB_{R})$. If $p \colon \oC \to \oB$ has cocartesian lifts
  of morphisms in $\oB_{L}$, then for a functor $F \colon \oC \to \oX$
  to a cocomplete \icat{} $\oX$, there is a natural equivalence
  between the left Kan extension $p_{!}F$ and the functor 
  \[ x \mapsto \colim_{\oC_{/x}^{R}} F \] where $\oC_{/x}^{R}$ denotes
  the fibre at $x$ of the functor
  $\oC \times_{\oB} \Ar_{R}(\oB) \xto{\ev_{1}} \oB$, where $\Ar_{R}(\oB)$
  is the full subcategory of $\Ar(\oB)$ spanned by the morphisms in
  $\oB_{R}$.
\end{propn}

For the proof we start with the following:
\begin{lemma}\label{lem:factLcocfreeadj}
  Suppose $\oB$ is an \icat{} with a factorization system
  $(\oB_{L},\oB_{R})$, and that $p \colon \oC \to \oB$ has cocartesian lifts
  of morphisms in $\oB_{L}$.
  \begin{enumerate}[(i)]
  \item The functor $\oC \times_{\oB} \Ar_{R}(\oB)
  \xto{\ev_{1}} \oB$ is a cocartesian fibration, with cocartesian
  morphisms given by
  \begin{equation}
    \label{eq:mor in C x_B Ar(B)}
    \begin{tikzcd}
      c \ar[d, "\gamma"] & p(c) \ar[r, "f"] \ar[d, "p(\gamma)"] & b \ar[d] \\
      c' & p(c') \ar[r, "g"] & b'
    \end{tikzcd}
  \end{equation}
  where $p(\gamma)$ lies in $\oB_{L}$ and $\gamma$ is
  $p$-cocartesian.
\item The inclusion
  $i \colon \oC \times_{\oB} \Ar_{R}(\oB) \to \oC \times_{\oB}
  \Ar(\oB)$ has a relative left adjoint $\lambda$ over $\oB$, and this
  preserves cocartesian morphisms.
  \end{enumerate}
\end{lemma}
\begin{proof}
  The first point is \cite[Proposition 2.2.2]{envelopes}. To see that $i$ has a
  left adjoint, it suffices to show that for every object
  $(c, p(c) \xto{f} b)$ in $\oC \times_{\oB} \Ar(\oB)$, the slice
  $(\oC \times_{\oB} \Ar_{R}(\oB))_{(c,f)/}$ has an initial object. If
  $p(c) \xto{f_{L}} b' \xto{f_{R}} b$ is the factorization of $f$
  according to the given factorization system, it is easy to check
  that $(f_{L,!}c, f_{R})$ gives the desired initial object.
  The cocartesian morphisms in $\oC \times_{\oB} \Ar(\oB)$ are those
  of the form \cref{eq:mor in C x_B Ar(B)} such that $\gamma$ is an
  equivalence. It follows from cancellation properties of $\oB_{L}$
  and cocartesian morphisms that $\lambda$ takes this to a cocartesian
  morphism  in $\oC \times_{\oB} \Ar_{R}(\oB)$.  Finally,
  both the unit and counit clearly lie over the identity of $\oB$ via
  $\ev_{1}$, so this is a relative adjunction.
\end{proof}

\begin{proof}[Proof of \cref{propn:lkefact}]
  From \cref{lem:factLcocfreeadj} we get natural transformations
  \[ p_{!}F(x) \simeq \colim_{\oC_{/x}} F \to \colim_{\oC_{/x}} F
    \circ i\lambda \to \colim_{\oC_{/x}^{R}} F\circ i,\]
  with the first induced by the unit transformation $\id \to i
  \lambda$ composed with $F$ and the second by the functor
  $\lambda$ in the index of the colimit. For a fixed $x$, these fit in
  a commutative diagram
  \[
    \begin{tikzcd}
       & \colim_{\oC_{/x}^{R}} F\circ i \ar[dl, "\sim"'] \ar[dr,
       "\sim"] \ar[drr, bend left, "="] \\
      \colim_{\oC_{/x}} F \ar[rr] & &  \colim_{\oC_{/x}} F  \circ
      i\lambda \ar[r] &  \colim_{\oC_{/x}^{R}} F\circ i.
    \end{tikzcd}
  \]
  Here the two left-most diagonal functors are induced by the functor
  $i \colon \oC^{R}_{/x} \to \oC_{/x}$; this is cofinal since $i$ is a
  right adjoint. Moreover, the composite $\lambda i$ is equivalent to
  the identity, hence so is the right composite as indicated. It
  follows that the maps in the bottom row are equivalences, which
  completes the proof.
\end{proof}

\begin{cor}\label{cor:cartfibinclLKE}
  Let $q \colon \oB \to \oA$ be a cartesian fibration of \icats{}, and
  let $j \colon \oB_{a} \to \oB$ denote the inclusion of the fibre at
  $a$. Then the left Kan extension along $j$ of a functor $F \colon \oB_{a} \to
  \oX$, where $\oX$ is a cocomplete \icat{}, is given by
  \[ x \mapsto \colim_{(\oB_{a})^{\cart}_{/x}} F.\]
  If $\oX = \CatIn$, then this left Kan extension is given in terms of
  cocartesian fibrations by
  \[ (\nE \to \oB_{a}) \mapsto (\nE \times_{\oB} \Ar_{\cart}(\oB)
    \xto{\ev_{1}} \oB).\]
  In particular, the left Kan extension of $\id_{\oB_{a}}$ corresponds
  to
  \[ \oB_{a} \times_{\oB} \Ar_{\cart}(\oB) \simeq \oA_{a/}
      \times_{\oA} \oB \to \oB,\] and that of $\oB_{a}(x,\blank)$ to the
    slice $\oB_{x/} \to \oB$.
\end{cor}
\begin{proof}
  The cartesian morphisms are the right class in a factorization
  system on $\oB$, with the fibrewise morphisms over $\oA$ as the left
  class. Moreover, $\oB_{a} \to \oB$ is trivially cocartesian over the
  fibrewise morphisms. The first statement therefore follows from
  \cref{propn:lkefact}.

  Applying this description to the identity of $\oB_{a}$, we get
  \[ i_{!}(*) \simeq \oA(a, p(\blank)). \]
  To get the description of $i_{!}\oB_{a}(x,\blank)$, we observe that
    \[ (\oB_{a})_{x/} \times_{\oB} \Ar_{\cart}(\oB) \simeq \oB_{x/}\]
    since any morphism $x \to y$ factors uniquely as $x \to x' \to y$
    where $x \to x'$ is in the fibre $\oB_{a}$ and $x' \to y$ is
    cartesian.
\end{proof}

\begin{cor}\label{cor:fibrelkeff}
  Suppose $p \colon \oB \to \oA$ is a cartesian fibration, and let $i
  \colon \oB_{a} \to \oB$ be the inclusion of its fibre at $a \in
  \oA$. Then left Kan extension along $i$ gives a fully faithful
  functor
  \[ \Fun(\oB_{a}, \CatIn) \hookrightarrow \Fun(\oB,
    \CatIn)_{/\oA(a,p(\blank))} \]
  whose image consists of those natural transformations $F \to
  \oA(a,p(\blank))$ such that when $\bar{f} \colon x \to y$ is a
  cartesian morphism over $f \colon b \to c$, the commutative square
  \[
    \begin{tikzcd}
      F(x) \ar[r, "F(\bar{f})"] \ar[d] & F(y) \ar[d]  \\
      \oA(a,b) \ar[r, "f_{*}"] & \oA(a,c)
    \end{tikzcd}
  \]
  is a pullback.
\end{cor}
\begin{proof}
  We have an adjunction
  \[ i_{!} : \Fun(\oB_{a}, \CatIn) \hookrightarrow \Fun(\oB,
    \CatIn) : i^{*} \]
  where $i_{!}$ is left Kan extension and $i^{*}$ is restriction along
  $i$. Using \cref{cor:cartfibinclLKE}, we know that $i_{!}$ takes
  the terminal object to $\oA(a,p(\blank))$, so we get an induced
  adjunction
  \[ i_{!} : \Fun(\oB_{a}, \CatIn) \hookrightarrow \Fun(\oB,
    \CatIn)_{/\oA(a, p(\blank))} : R \]
  where $R$ takes a natural transformation $\gamma \colon F \to \oA(a, p(\blank))$
  to the pullback
  \[
    \begin{tikzcd}
      R(\gamma) \ar[r] \ar[d] \drpullback & i^{*}F \ar[d] \\
      * \ar[r] & \oA(a, a)
    \end{tikzcd}
  \]
  (since $\oA(a, pi(\blank))$ is constant at $\oA(a,a)$, and the
  bottom horizontal map picks out $\id_{a}$).
  Since the colimit in $\CatIn$ of a diagram indexed by an \igpd{} can be computed as the
  domain of the corresponding cocartesian fibration, our formula for
  $i_{!}F$ implies that the natural map $F \to R(i_{!}F)$ is an
  equivalence, since we have pullback squares
  \[
    \begin{tikzcd}
      F \ar[r] \ar[d] \drpullback & i^{*}i_{!}F \ar[d] \\
      * \ar[r] & \oA(a, a).
    \end{tikzcd}
  \]
  Therefore, $i_{!}$ is fully faithful as a functor to
  $\Fun(\oB, \CatIn)_{/\oA(a, p(\blank))}$. Its image consists of
  those $\gamma \colon F \to \oA(a,p(\blank))$ such that for every $x
  \in \oB$ the horizontal map in the commutative triangle
  \[
    \begin{tikzcd}
      \colim_{f \in \oA(a,p(x))} F(f^{*}x)_{\id_{a}} \ar[rr] \ar[dr] &  & F(x)
      \ar[dl] \\
       & \oA(a,p(x))
    \end{tikzcd}
  \]
  is an equivalence. We can detect this equivalence on fibres, so this
  condition is equivalent to
  \[ F(f^{*}x)_{\id_{a}} \to F(x)_{f}\]
  being invertible for every $f \colon a \to p(x)$. We claim this is
  equivalent to our given condition: On the one hand, we can recognize
  the map on fibres as one of the maps on fibres in the commutative
  square
    \[
    \begin{tikzcd}
      F(f^{*}x) \ar[r, "F(\bar{f})"] \ar[d] & F(x) \ar[d]  \\
      \oA(a,a) \ar[r, "f_{*}"] & \oA(a,p(x)),
    \end{tikzcd}
  \]
  where $\bar{f}\colon f^{*}x \to x$ is cartesian over $x$. On the
  other hand, if we consider the fibre over $g \colon a \to b$
  in the square
    \[
    \begin{tikzcd}
      F(x) \ar[r, "F(\bar{f})"] \ar[d] & F(y) \ar[d]  \\
      \oA(a,b) \ar[r, "f_{*}"] & \oA(a,c),
    \end{tikzcd}
  \]
  we can find a cartesian morphism $\bar{g} \colon x' \to x$ over $g$,
  and then we get a commutative triangle
  \[
    \begin{tikzcd}
      & F(x')_{\id_{a}} \ar[dl] \ar[dr] \\
      F(x)_{g} \ar[rr] & & F(y)_{fg},
    \end{tikzcd}
  \]
  which completes the proof.
\end{proof}

\begin{remark}
  The condition characterizing the image in \cref{cor:fibrelkeff} is a
  case of the ``equifibred'' condition from \cite[Proposition 2.3.3]{envelopes}.
\end{remark}

The following observation lets us extend these results to left Kan extensions along symmetric
monoidal functors of the form $\oB_{\bbone}^{\otimes} \to \oB^{\otimes}$:
\begin{lemma}\label{lem:opd lke ordinary}
  Suppose $\phi \colon \oA \to \oB$ is a symmetric monoidal functor
  and $F \colon \oA \to \uX$ is lax symmetric monoidal, with
  $\uX^{\otimes}$ compatible with small colimits. Then the operadic
  left Kan extension $\phi_{!}^{\mathrm{opd}}F$ along $\phi$ has the
  ordinary left Kan extension $\phi_{!}F$ as its underlying functor.
\end{lemma}
\begin{proof}
  We must show that for $b$ in $\oB$, the functor
  \[ \oA_{/b} \to \oA^{\otimes,\act}_{/b} \] is cofinal. Given an
  object in the target of the form
  \[\alpha = ((a_{1},\ldots,a_{n}), \phi(a_{1}, \ldots, a_{n}) \to b),\]
  we must thus show that $(\oA_{/b})_{\alpha/}$ is weakly
  contractible. But since $\phi$ is symmetric monoidal 
  this has an initial object, given by the cocartesian morphism
  $(a_{1},\ldots,a_{n}) \to a_{1} \otimes \cdots \otimes a_{n}$.
\end{proof}

\begin{proof}[Proof of \cref{propn:2-step fibrewise}]
  We can identify $\oB^{\otimes}$-monoids in $\CatIn$ with
  $\oB^{\otimes}$-algebras in the cartesian symmetric monoidal \icat{}
  $\CatIn^{\times}$. Operadic left Kan extension \cite[\S 3.1.2]{HA}
  therefore provides an adjunction
  \[ i_{!} : \Mon_{\oB_{\bbone}^{\otimes}}(\CatIn) \to
    \Mon_{\oB^{\otimes}}(\CatIn) : i^{*},\]
  where $i_{!}$ is given by ordinary left Kan extension on underlying
  functors by \cref{lem:opd lke ordinary}. In particular, $i_{!}$
  takes the terminal object to $\oA(\bbone,p(\blank))$, so there is an
  induced adjunction
  \[ i_{!} : \Mon_{\oB_{\bbone}^{\otimes}}(\CatIn) \to
    \Mon_{\oB^{\otimes}}(\CatIn)_{/\oA(\bbone,p(\blank))} : R_{p}.\]
  To see that $i_{!}$ is fully faithful and identify its image it
  suffices to note that the forgetful functors to
  $\Fun(\oB_{\bbone}, \CatIn)$ and
  $\Fun(\oB, \CatIn)_{/\oA(\bbone,p(\blank))}$ are conservative, and
  the adjunction restricts there to that of \cref{cor:fibrelkeff}.
\end{proof}

\subsection{Skeletal sequences and the inductive step}
\label{sec:skeletal}

Armed with \cref{thm:2-step}, we are now ready to set up the inductive
step in our proof. For this we must first introduce an extended
version of the 2-pinnings of left 2-rigid symmetric monoidal \incats{}
we have considered so far.

\begin{defn}\label{defn:left m-rig}
  We say a symmetric monoidal \incat{} $\nC$ is \emph{left $m$-rigid}
  for $0 \leq m \leq n$ if
  \begin{enumerate}[(1)]
  \item $\nC$ has duals for objects, if $m > 0$,
  \item every $1$-morphism in $\nC$ of the form $\bbone \to X$ has a
    left adjoint, if $m > 1$,
  \item the symmetric monoidal \incat{} $\mOm \nC$ is left
    $(m-1)$-rigid, if $m > 0$.
  \end{enumerate}
  In the maximal case where $\nC$ is left $n$-rigid, we will simply
  say that $\nC$ is \emph{left rigid}.
  Note that in the degenerate case $m = 0$ \emph{any}
  symmetric monoidal \incat{} is left $0$-rigid.
\end{defn}

\begin{observation}\label{obs:m-rig defn}
  The inductive definition is slightly redundant: it follows from
  \cref{propn:PCdualadj} that $\mOm \nC$ automatically has duals if
  (2) holds, so instead of (3) it suffices to ask for condition (2) to
  hold for $\mOm^{k}\nC$ for $k = 0,1,\dots,m-2$. We can furthermore
  use \cref{propn:PCdualadj} to see that $\nC$ is 
  left $m$-rigid for $m>0$ \IFF{}
  $\nC$ and $\mP\mOm^{k}\nC$ have duals for objects for all $k = 0,\dots,m-2$.
\end{observation}

\begin{question}
  Suppose $\nC$ is a symmetric monoidal \incat{} which, when regarded
  as an $(\infty,m)$-category, is left $m$-rigid for some $m>n$. If
  $n = 1$, then the \icat{} $\nC$ is an \igpd{} by
  \cref{propn:left 2-rigid cat is gpd}, and it is tempting to guess
  that this is also true for larger $n$. If $\nC$ is a left 3-rigid
  \itcat{} then $\mOm\nC$ is an \igpd{}, but it is unclear whether or not
  this forces all 2-morphisms to be invertible in general.
\end{question}

\begin{defn}
  We say a functor of \incats{} $F \colon \nC \to \nD$ is
  \emph{$k$-surjective} for $k \geq 0$ if
  \begin{itemize}
  \item $F$ is essentially surjective, \ie{} $\pi_{0}\nC^{\simeq} \to
    \pi_{0}\nD^{\simeq}$ is surjective,
  \item if $k > 0$, then for all objects $X,Y \in \nC$, the functor of
    $(\infty,n-1)$-categories $\nC(X,Y) \to \nD(FX,FY)$ induced by $F$
    is $(k-1)$-surjective.
  \end{itemize}
\end{defn}

\begin{remark}\label{rmk:k-surj gpd}
  In \cite{LurieCob}, $k$-surjective functors are called
  \emph{$k$-connective}. Note that a morphism of \igpds{}
  $F \colon X \to Y$ is $k$-surjective \IFF{} it is surjective on
  connected components and for all base points $x$,
  the map $\pi_{i}(X,x) \to \pi_{i}(Y,F(x))$ is an isomorphism for
  $i < k$ and surjective for $i = k$. This is equivalent to $F$ being
  surjective on connected components and the homotopy groups of the
  fibres of $F$ vanishing in dimensions $< k$, \ie{} to $F$ having
  $k$-connective fibres. 
\end{remark}

\begin{lemma}
  Suppose $\nC$ is an \incat{}, $\nD$ is an $(\infty,n+1)$-category,
  and $F \colon \nC \to \nD$ is a $k$-surjective functor for $k >
  n$. Then:
  \begin{enumerate}[(i)]
  \item $\nD$ is an $(\infty,n)$-category.
  \item The functor $F$ is conservative on $i$-morphisms for all $i$.
  \item For $0 \leq i \leq n$, the induced functor $\core[i] F \colon \core[i]\nC \to \core[i]\nD$ is also
    $k$-surjective, and so
    \begin{itemize}
    \item $\nC^{\simeq} \to \nD^{\simeq}$ is a $k$-surjective morphism
      of \igpds{},
    \item for all $x,y \in \nC$, the map $\nC(x,y)^{\simeq} \to
      \nD(Fx,Fy)^{\simeq}$ is a $(k-1)$-surjective morphisms of
      \igpds{}, etc.
    \end{itemize}
  \end{enumerate}
\end{lemma}
\begin{proof}
  Part (i) is immediate, since any $i$-morphism in the image of $F$
  for $i > n$ is necessarily invertible. We prove (ii) by observing
  that $F$ is trivially conservative on $(n+1)$-morphisms and then
  showing that if $F$ is conservative on $(i+1)$-morphisms then it is
  so on $i$-morphisms: if $f \colon x \to y$ is an $i$-morphism so
  that $F(f)$ is invertible, then in $\nD$ we have an inverse $g$ and
  invertible $(i+1)$-morphisms $gF(f) \isoto \id_{F(x)}$,
  $F(f)g \isoto \id_{F(y)}$. Here $g$ is the image of some morphism
  $g' \colon y \to x$ in $\nC$, and we can lift the invertible
  $(i+1)$-morphisms to $(i+1)$-morphisms $g'f \to \id_{x}$,
  $fg' \to \id_{y}$. Since $F$ is conservative on $(i+1)$-morphisms
  these lifts are necessarily invertible, and so exhibit $g'$ as the
  inverse of $f$ in $\nC$. It follows that an invertible $i$-morphism
  in $\nD$ can always be lifted to an \emph{invertible} $i$-morphism
  in $\nC$, which gives (iii).
\end{proof}

\begin{defn}
  Let $\nC$ be a left $m$-rigid symmetric monoidal \incat{}, where $m
  \leq n$. A $0$-pinning of $\nC$ is, by definition, no additional
  data. For $1 \leq k \leq m$ and $k<n$, a
  \emph{$k$-pinning} of $\nC$ is a sequence
  \[ \nA_{1} \to \nA_{2} \to \dots \to \nA_{k} \to \nC \] of strong
  symmetric monoidal functors, where
  \begin{itemize}
  \item   $\nA_{i}$ is a left $i$-rigid
  symmetric monoidal $(\infty,i)$-category,
\item  $\nA_{i} \to \nA_{i+1}$ is $(i-1)$-surjective,
\item and $\nA_{k} \to \nC$ is $(k-1)$-surjective.
  \end{itemize}
 Let $\PSMCatIn[k]^{m\dlrig}$ denote the \icat{} of
  $k$-pinned left $m$-rigid symmetric monoidal \incats{}, which we can
  regard as a full subcategory of $\Fun([k], \SMCatIn)$. For $m = n$
  we abbreviate $\PSMCatIn[k]^{\lrig} := \PSMCatIn[k]^{n\dlrig}$.
  In particular, $\PSMCatIn[0]^{m\dlrig}=\SMCatIn^{m\dlrig}$.
\end{defn}

\begin{defn}
  Following \cite{LurieCob}, we will refer to an $(n-1)$-pinned left $n$-rigid symmetric
  monoidal \incat{} as a \emph{skeletal sequence}
  \[ \nA_{1} \to \nA_{2} \to \cdots \to \nA_{n}\]
  of length $n$. We write
  \[ \SkSeq_{n} := \PSMCatIn[n-1]^{\lrig}  \]
  for the \icat{} thereof.
\end{defn}

\begin{remark}
  An $(n-1)$-pinned left $(n-1)$-rigid symmetric monoidal \incat{} is
  called a \emph{weak skeletal sequence} of length $n$ in
  \cite{LurieCob}.
\end{remark}

\begin{observation}\label{obs:k-pinning as arrow}
  By viewing
  \[ \nA_{1} \to \nA_{2} \to \dots \to \nA_{k} \to \nC \]
  as a natural transformation
  \[
    \begin{tikzcd}
     \nA_{1} \ar[r] \ar[d, equals] &  \nA_{2} \ar[r] \ar[d, equals] &
     \cdots  \ar[r]  &  \nA_{k} \ar[d] \\
     \nA_{1} \ar[r] &  \nA_{2} \ar[r] &
     \cdots  \ar[r]  &  \nC,
    \end{tikzcd}
  \]
  we can regard $\PSMCatIn[k]^{m\dlrig}$ as the full subcategory of
  $\Ar(\PSMCatIn[k-1]^{k\dlrig})$ containing those objects
    \[
    \begin{tikzcd}
     \nA_{1}' \ar[r] \ar[d] &  \nA_{2}' \ar[r] \ar[d] &
     \cdots  \ar[r]  &  \nB \ar[d] \\
     \nA_{1} \ar[r] &  \nA_{2} \ar[r] &
     \cdots  \ar[r]  &  \nC
    \end{tikzcd}
  \]  
  where $\nA_{i}' \to \nA_{i}$ is an equivalence for $i=1,\ldots,k-1$,
  $\nB$ is an $(\infty,k)$-category, and $\nC$ is left $m$-rigid.
\end{observation}

\begin{lemma}\label{lem:P2 to CC m-lrig}
  The equivalence of \cref{cor:P2 to CC} restricts on full
  subcategories to equivalences
  \[ \PSMCat[2]_{(\infty,n)}^{m\dlrig} \isoto
    \PSMCat_{(\infty,n-1)}^{(m-1)\dlrig} \times_{\SMCatI^{\rig}}
    \SMFib_{1/1}^{\strig}\]
  and
  \[ \SMCat_{(\infty,n)}^{m\dlrig} \isoto
    \SMCat_{(\infty,n-1)}^{(m-1)\dlrig} \times_{\SMCatI^{\rig}}
    \CSMFib_{1/1}^{\strig}.\]
\end{lemma}
\begin{proof}
  This is immediate from the inductive definition of left $m$-rigidity
  (\cref{defn:left m-rig}) and the definition of these functors.
\end{proof}

We are finally ready to state and prove the inductive step in our
``unfolding'' construction:

\begin{propn}\label{propn:1-step m-rig pinned}
  For $k \geq 1$, the equivalence of \cref{cor:P2 to CC} induces an equivalence
    \[ \PSMCat[k+1]_{(\infty,n)}^{m\dlrig} \simeq
    \PSMCat[k]_{(\infty,n-1)}^{(m-1)\dlrig} \times_{\SMCatI^{\rig}}
    \SMFib_{1/1}^{\strig},\]
  which takes
  \[ \nA_{1} \to \nA_{2} \to \dots \to \nA_{k+1} \to \nC \]
  to $\mOm \nA_{2} \to \dots \to \mOm \nA_{k+1} \to \mOm \nC$ together
  with $\nA_{1} \times_{\core \nA_{2}}\mP \nA_{2} \to \nA_{1}$. In
  particular, we have
  \[ \SkSeq_{n} \isoto \SkSeq_{n-1} \times_{\SMCatI^{\rig}}
    \SMFib_{1/1}^{\strig}.\]
\end{propn}
\begin{proof}
  We induct on $k$, with the case $k = 1$ being \cref{lem:P2 to CC
    m-lrig}. By \cref{obs:k-pinning as arrow} we can identify
  $\PSMCat[k+1]_{(\infty,n)}^{m\dlrig}$ as a full subcategory of
  \[ \Ar(\PSMCat[k]_{(\infty,n)}^{k\dlrig}) \simeq
    \Ar(\PSMCat[k-1]_{(\infty,n-1)}^{(k-1)\dlrig}
    \times_{\SMCatI^{\rig}} \SMFib_{1/1}^{\strig}).\] Moreover, the
  description of the \icat{} $\PSMCat[k]_{(\infty,n-1)}^{(m-1)\dlrig}$ as a full
  subcategory of $\Ar(\PSMCat[k-1]_{(\infty,n-1)}^{(k-1)\dlrig})$
  extends to an identification 
  \[\PSMCat[k]_{(\infty,n-1)}^{(m-1)\dlrig} \times_{\SMCatI^{\rig}}
  \SMFib_{1/1}^{\strig} \subseteq 
  \Ar(\PSMCat[k-1]_{(\infty,n-1)}^{(k-1)\dlrig}
  \times_{\SMCatI^{\rig}} \SMFib_{1/1}^{\strig}) \] as the full subcategory given by objects that
  are constant in the other factors of the fibre product. Here this
  gives precisely the image of $\PSMCat[k+1]_{(\infty,n)}^{m\dlrig}$, as required.
\end{proof}

\subsection{Chain complexes of symmetric monoidal
  \icats{}}\label{sec:chain complex unfolding}

We can now apply \cref{propn:1-step m-rig pinned} iteratively to
unfold skeletal sequences and prove
\cref{introthm:skseq}. To spell this out, we first introduce some terminology:

\begin{defn}\label{defn:chain complex cat}
  We inductively define the \icat{} of (non-negatively graded)
  \emph{chain complexes} of symmetric monoidal \icats{} with duals of
  length $n$ by
  \[ \CChp_{0} := \SMCatI^{\rig},\]
  \[ \CChp_{1} := \SMFib_{1/1}^{\strig},\]
  \[ \CChp_{n} := \CChp_{1} \times^{t,s}_{\CChp_{0}} \CChp_{n-1},\]  
  where we inductively define the functors $s$ and $t$ by defining $s \colon \SMFib_{1/1}^{\strig}
  \to \SMCatI^{\rig}$ to be the functor that takes the fibre at
  $\bbone$ and $t$ to be projection to the target, and then setting
  \[ s \colon \CChp_{n} \to \CChp_{1} \xto{s} \CChp_{0},  \quad t \colon
    \CChp_{n} \to \CChp_{n-1} \xto{t} \CChp_{0},\]
  where the first functor in each case is the projection from the
  pullback defining $\CChp_{n}$.
  Note that we then have equivalences
  \[ \CChp_{n+m} \simeq \CChp_{n} \times^{t,s}_{\CChp_{0}} \CChp_{m}\]
  for all $n,m$.
\end{defn}

\begin{observation}
  We can think of an object of $\CChp_{n}$ as a sequence of pullback
  squares of symmetric monoidal \icats{} with duals
  \begin{equation}
    \label{eq:ch cx pb}
    \begin{tikzcd}
      \oZ_{k} \drpullback \ar[r] \ar[d] & \oC_{k} \ar[d, "\partial_{k}"] \\
      0  \ar[r] & \oZ_{k-1},
    \end{tikzcd}
  \end{equation}
  for $1 \leq k \leq n$, where the underlying functor of
  $\partial_{k}$ is a cocartesian fibration, together with the
  degenerate case
  \[
    \begin{tikzcd}
      \oZ_{0} \drpullback \ar[r, "\sim"] \ar[d] & \oC_{0} \ar[d] \\
      0  \ar[r] & 0.
    \end{tikzcd}
  \]
\end{observation}

\begin{warning}
  Our indexing convention differs from that of \cite{LurieCob}:
  Lurie takes $\oZ_{0} = 0$, so $\oZ_{1} \simeq \oC_{1}$, and a chain
  complex of length $n$ in Lurie's sense corresponds to a symmetric
  monoidal $(\infty,n)$-category. We find our convention more
  convenient for the extension to negative chain complexes we consider
  below.
\end{warning}

\begin{remark}\label{rmk:chain complex motivation}
  The term ``chain complex'' comes from the analogy with the following
  description of a non-negatively graded chain complex of abelian groups
  $(C_{\bullet}, \partial)$: if $Z_{i} \subseteq C_{i}$ is the
  subgroup of cycles, then we can describe the chain complex as a
  sequence of pullback squares of abelian groups
    \[
    \begin{tikzcd}
      Z_{n+1} \drpullback \ar[r] \ar[d] & C_{n+1} \ar[d, "\partial"] \\
      0  \ar[r] & Z_{n},
    \end{tikzcd}
  \]
  where the relation $\partial^{2} = 0$ is encoded by viewing the
  differential as taking values in the subgroup of cycles.
\end{remark}

\begin{remark}
  As a special case of \cref{defn:chain complex cat}, we can consider
  those chain complexes of symmetric monoidal \icats{} with duals
  that consist of \igpds{}. Note that a symmetric monoidal \igpd{} has
  duals \IFF{} it is grouplike, so these objects are given by a
  sequence of pullback squares
  \begin{equation}
    \label{eq:cchpb}
    \begin{tikzcd}
      Z_{k} \drpullback \ar[r] \ar[d] & C_{k} \ar[d, "\partial_{k}"] \\
      0  \ar[r] & Z_{k-1},
    \end{tikzcd}
  \end{equation}
  where $C_{i}$ and $Z_{i}$ are grouplike $E_{\infty}$-spaces (or
  connective spectra). 
\end{remark}

\begin{cornoproof}\label{cor:unfold pinned}
  Iterating \cref{propn:1-step m-rig pinned}, we have equivalences
  \[ \PSMCat[k]_{(\infty,n)}^{m\dlrig} \simeq
    \PSMCat[k-i]_{(\infty,n-i)}^{(m-i)\dlrig} \times_{\CChp_{0}} \CChp_{i}
  \]
  for $0 \leq i \leq k$. In particular, for $i = k$ we have
  \[ \PSMCat[k]_{(\infty,n)}^{m\dlrig} \simeq
    \SMCat_{(\infty,n-k)}^{(m-k)\dlrig} \times_{\CChp_{0}}
    \CChp_{k}\]
  for $k < m$.
\end{cornoproof}

Taking $k = n-1$, we get the following special case:
\begin{cornoproof}\label{cor:skseq cch eq}
  There is an equivalence
  \[ \SkSeq_{n} \simeq 
    \CChp_{n-1}\]
  that takes a skeletal sequence $\nA_{1} \to \dots \to \nA_{n}$ to
  the sequence of pullback squares
\[
  \begin{tikzcd}
   \mOm^{i}\nA_{i+1} \ar[r] \ar[d] & \mP\mOm^{i-1}\nA_{i+1}  \ar[d] \\
   0 \ar[r] & \mOm^{i-1}\nA_{i}
  \end{tikzcd}
\]
of rigid symmetric monoidal \icats{}.
\end{cornoproof}

\begin{observation}\label{obs:weak chain cx}
  Dropping the top level of dualizability, \cref{cor:unfold pinned}
  gives an equivalence
  \[ \PSMCat[n-1]_{(\infty,n)}^{(n-1)\dlrig} \simeq
    \PSMCat_{(\infty,2)}^{1\dlrig} \times_{\CChp_{0}}
    \CChp_{n-2}.\]
  Now we can apply \cref{cor:fibeqduals} to get an equivalence
  \[ \PSMCat[n-1]_{(\infty,n)}^{(n-1)\dlrig} \simeq
    \SMFibtrig_{1/1} \times_{\CChp_{0}}
    \CChp_{n-2};\]
  here we can interpret the right-hand side as the \icat{} of chain
  complexes of length $n-1$ where we do not ask for rigidity at the
  top level. These objects are called \emph{weak chain complexes} in
  \cite{LurieCob}, so this provides precisely the equivalence between these and
  weak skeletal sequences stated in \cite[Proposition 3.3.33]{LurieCob}.  
\end{observation}

We now deduce \cref{introthm:complete} by introducing a completeness
condition on chain complexes:
\begin{defn}
  We define a full subcategory $\CChpc_{n} \subseteq \CChp_{n}$ by
  inductively setting
  \[ \CChpc_{0} := \CChp_{0}, \quad \CChpc_{1} := \CSMFib_{1/1}^{\strig}
    \subseteq \CChp_{1},\]
  \[ \CChpc_{n} := \CChpc_{n-1} \times_{\CChpc_{0}}^{s,t} \CChpc_{1}
    \subseteq \CChp_{n}.\]
  We say that a chain complex of length $n$
  is \emph{complete} if it lies in $\CChpc_{n}$.
\end{defn}

\begin{observation}
  Unwinding the definitions, a chain complex given by pullback squares
  \cref{eq:ch cx pb} is complete \IFF{} the canonical map
  \[
    \begin{tikzcd}
     (\oZ_{k-1})_{\bbone/} \ar[rr] \ar[dr] & & \oC_{k} \ar[dl, "\partial_{k}"] \\
      & \oZ_{k-1}
    \end{tikzcd}
  \]
  gives an equivalence on groupoid cores for every $k$.
\end{observation}

\begin{cor}\label{cor:finite complete unfolding}
  The equivalence of \cref{cor:unfold pinned} restricts on full
  subcategories to an equivalence
  \[ \SMCat_{(\infty,n)}^{m\dlrig} \simeq
    \SMCat_{(\infty,n-k)}^{(m-k)\dlrig} \times_{\CChpc_{0}}
    \CChpc_{k}\]
  for $k < m$, and so
  \[ \SMCat_{(\infty,n+1)}^{\lrig} \simeq 
    \CChpc_{n}.\]
\end{cor}
\begin{proof}
  Unpacking the constructions, this follows by applying
  \cref{cor:fibcatneq} inductively.
\end{proof}

\section{Extensions and localizations}
\label{sec:extend}

In the first part of this section we extend the unfolding equivalence,
first to non-negatively graded categorical chain complexes of infinite
length in \cref{sec:iicats}, and then further to the negatively
unbounded case in \cref{sec:unbounded}. We then study the localization
to \emph{complete} categorical chain complexes in
\cref{sec:completion}.
In the last part of the section we study chain complexes of \igpds{}
and the localization to these in \cref{sec:gpdchcx}, and then
specialize the description of completion to \igpds{} in \cref{sec:comp
  igpd}, which gives the homotopy-group formulas in
\cref{cor:pi of completion of locn}.

\subsection{Extending up to $(\infty,\infty)$-categories}
\label{sec:iicats}

In this section we extend the unfolding equivalence to non-negatively
graded chain complexes of infinite length, relating these to
infinitely long skeletal sequences, and in the complete case to left
rigid symmetric monoidal $(\infty,\infty)$-categories.

\begin{defn}
  We have functors $p \colon \SkSeq_{n+1}  \to \SkSeq_{n}$ that forget the last term in an
  $n$-pinned sequence, \ie{}
  \[ (\nA_{1} \to \cdots \to \nA_{n+1}) \mapsto (\nA_{1} \to \cdots
    \to \nA_{n}).\]
  (Formally, these can be defined by restricting the functors
  \[ d_{n}^{*} \colon \Fun([n],\SMCat_{(\infty,n+1)}) \to \Fun([n-1],
    \SMCat_{(\infty,n+1)})\]
  to full subcategories.)
  We define
  \[ \SkSeq_{\infty} = \PSMCat[\infty]_{(\infty,\infty)}^{\lrig} := \limP_{n \rightarrow
      \infty} \PSMCat[n]_{(\infty,n+1)}^{\lrig}\]
  using these functors. An object of
  $\SkSeq_{\infty}$ is then an
  infinite sequence
  \[ \nA_{1} \xto{i_{1}} \nA_{2} \to \cdots \to \nA_{n} \xto{i_{n}}
    \nA_{n+1}  \to \cdots \]
  where
  \begin{itemize}
  \item $\nA_{n}$ is a left $n$-rigid symmetric monoidal
    $(\infty,n)$-category,
  \item $i_{n}$ is a strong symmetric monoidal functor that is $(n-1)$-surjective.
  \end{itemize}
  The functors defining the limit are compatible with the inclusion of
  the full subcategories $\SMCat_{(\infty,n)}^{\lrig}
  \hookrightarrow \PSMCat[n-1]_{(\infty,n)}^{\lrig}$, so in the
  limit we have a full subcategory
  \[ \SMCat_{(\infty,\infty)}^{\lrig} := \limP_{n \rightarrow
      \infty} \SMCat_{(\infty,n)}^{\lrig}\]
  which consists of sequences as above where $i_{n}$ exhibits
  $\nA_{n}$ as the $n$-core $\core[n] \nA_{n+1}$.
\end{defn}

\begin{observation}
  Recall that the \icat{} $\CatII$ of \emph{\iicats{}} is the limit
  $\limP_{n \rightarrow \infty} \CatIn$ along the core functors $\core[n-1]
  \colon \CatIn \to \CatInm$. These functors preserve products (being
  right adjoints), so this also exhibits the cartesian symmetric
  monoidal structure on $\CatII$ as a limit, so that we have
  \[ \SMCatII \simeq \limP_{n \rightarrow \infty} \SMCatIn.\] This limit can
  be identified as the cartesian sections of the corresponding
  cartesian fibration $q \colon \oE \to \mathbb{N}$. Let
  $\catname{Sec}(q) = \Fun_{/\mathbb{N}}(\mathbb{N}, \oE)$ be the
  \icat{} of all sections of $q$; this has as objects sequences of
  objects $\nC_{n} \in \SMCatIn$ together with symmetric monoidal
  functors $\nC_{n} \to \core[n] \nC_{n+1}$, and it is easy to see
  that $\SkSeq_{\infty}$ is a full subcategory of it. Under this
  identification $\SMCatII^{\lrig}$ is precisely the intersection of
  $\SkSeq_{\infty}$ and $\SMCatII$ inside $\catname{Sec}(q)$. In
  particular, we may regard $\SMCatII^{\lrig}$ as a full subcategory of
  $\SMCatII$, as suggested by the notation.
\end{observation}

\begin{defn}
  There are functors
  \[ b \colon \CChp_{n} \to \CChp_{n-1}\]
  that forget the top square in the chain complex, and we set
  \[ \CChp_{\infty} := \limP_{n \rightarrow \infty} \CChp_{n}\] for
  the limit along these functors; this is the \icat{} of
  non-negatively graded categorical chain complexes of infinite
  length. These are compatible with the full subcategories of complete
  objects, so $\CChp_{\infty}$ has a full subcategory
  $\CChpc_{\infty} := \limP_{n \rightarrow \infty} \CChpc_{n}.$
\end{defn}

\begin{defn}
  We write \[\mOm \colon \SkSeq_{n+1} \to \SkSeq_{n}\] for the functor
  that takes a skeletal sequence $\nA_{1} \to \cdots \to \nA_{n+1}$ to
  $\mOm \nA_{2} \to \cdots \to \mOm \nA_{n+1}$; formally, this can be defined by
  restricting the composite functor
  \[ \Fun([n], \SMCat_{(\infty,n+1)}) \xto{d_{0}^{*}} \Fun([n-1],
    \SMCat_{(\infty,n+1)}) \xto{\mOm_{*}} \Fun([n-1], \SMCatIn)\]
  to full subcategories.
\end{defn}

\begin{observation}\label{obs:Omega bottom commute} 
  There is a commutative square
  \[
    \begin{tikzcd}
     \SkSeq_{n+1} \ar[r, "\mOm"] \ar[d, "p"'] & \SkSeq_{n}  \ar[d, "p"] \\
     \SkSeq_{n} \ar[r, "\mOm"'] & \SkSeq_{n-1},
    \end{tikzcd}
  \]
  with both functors taking the skeletal sequence $\nA_{1} \to \cdots
  \to \nA_{n+1}$ to $\mOm \nA_{2} \to \cdots \to \mOm
  \nA_{n}$. Indeed, this is a restriction to full subcategories of
  the outer square in the commutative diagram
  \begin{equation}
    \label{eq:Omega bottom commute diagram}
    \begin{tikzcd}
      \Fun([n], \SMCat_{(\infty,n+1)}) \ar[d, "d_{0}^{*}"] \ar[r, "d_{n}^{*}"] & 
      \Fun([n-1],\SMCat_{(\infty,n+1)}) \ar[d, "d_{0}^{*}"]
      \\
\Fun([n-1], \SMCat_{(\infty,n+1)})  \ar[r, "d_{n-1}^{*}"] \ar[d,
      "\mOm_{*}"] &
            \Fun([n-2], \SMCat_{(\infty,n+1)}) \ar[d,
      "\mOm_{*}"]
      \\
      \Fun([n-1], \SMCatIn) \ar[r, "d_{n-1}^{*}"] 
 & \Fun([n-2], \SMCatIn).
    \end{tikzcd}
  \end{equation}\end{observation}

\begin{propn}\label{propn:bottom commute eqce}
  The equivalences of \cref{propn:1-step m-rig pinned} fit in a
  commutative square
  \[
    \begin{tikzcd}
     \SkSeq_{n} \ar[r, "\sim"] \ar[d, "p"'] &
     \SkSeq_{n-1}\times_{\CChp_{0}} \CChp_{1}  \ar[d, "p \times_{\id} \id"] \\
     \SkSeq_{n-1} \ar[r, "\sim"'] & \SkSeq_{n-2} \times_{\CChp_{0}}
     \CChp_{1}.
    \end{tikzcd}
  \]
\end{propn}
\begin{proof}
  By definition the equivalence
  \[\SkSeq_{n} \isoto \SkSeq_{n-1}\times_{\CChp_{0}} \CChp_{1}\]
  is given by applying $\mOm$ in the first factor and the composite
  \[ \SkSeq_{n} \xto{p^{n-2}} \SkSeq_{2} \isoto \CChp_{1}\] in the
  second factor, so this follows immediately from \cref{obs:Omega
    bottom commute}.
\end{proof}

\begin{cor}\label{cor:skseq cch nateq}
  There is a commutative square
  \[
    \begin{tikzcd}
     \SkSeq_{n} \ar[r, "\sim"] \ar[d, "p"'] & \CChp_{n-1}  \ar[d, "b"] \\
     \SkSeq_{n-1} \ar[r, "\sim"'] & \CChp_{n-2}.
    \end{tikzcd}
  \]
\end{cor}
\begin{proof}
  We induct on $n$. The case $n = 2$ is obvious, and combining the
  case for $n-1$ with
  \cref{propn:bottom commute eqce} we get a commutative diagram
  \[
    \begin{tikzcd}
     \SkSeq_{n} \ar[r, "\sim"] \ar[d, "p"'] &
     \SkSeq_{n-1}\times_{\CChp_{0}} \CChp_{1}  \ar[d, "p \times_{\id}
     \id"] \ar[r, "\sim"] & \CChp_{n-1} \ar[d, "b"]\\
     \SkSeq_{n-1} \ar[r, "\sim"'] & \SkSeq_{n-2} \times_{\CChp_{0}}
     \CChp_{1} \ar[r, "\sim"] & \CChp_{n-2},
    \end{tikzcd}
  \]
  which completes the proof.
\end{proof}

From \cref{cor:skseq cch nateq} we get a natural equivalence
$\SkSeq_{\bullet} \simeq \CChp_{\bullet-1}$; in the limit we therefore get the
following:
\begin{cor}\label{cor:unfolding infinite}
  The equivalences of \cref{cor:skseq cch eq} combine to equivalences
  \[ \SkSeq_{\infty} \simeq
    \CChp_{\infty}, \quad \SMCat_{(\infty,\infty)}^{\lrig}
    \simeq \CChpc_{\infty}.\]
\end{cor}

\subsection{Extending down to $(\infty,\infty)$-spectra}
\label{sec:unbounded}

In this subsection we extend the unfolding equivalence to chain
complexes that are unbounded both above and below. These will be
equivalent to \emph{delooping sequences} of infinite skeletal
sequences (or of left rigid symmetric monoidal \iicats{} in the
complete case).

\begin{observation}
    It follows from \cref{obs:Omega bottom commute} that the looping
  functors give a natural transformation $\mOm \colon \SkSeq_{\bullet} \to
  \SkSeq_{\bullet-1}$, and so in the limit a functor
  \[ \mOm \colon \SkSeq_{\infty}  \to \SkSeq_{\infty}.\]
  These functors restrict on full subcategories to the looping
  functors
  \[ \mOm \colon \SMCatIn^{\lrig} \to \SMCatInm^{\lrig},\]
  so in the limit we get a restriction to 
  \[ \mOm \colon \SMCatII^{\lrig} \to
    \SMCatII^{\lrig}.\]
\end{observation}

\begin{defn}
  We define
  \[ \PSpIIl := \lim \left(\cdots \to \SkSeq_{\infty} \xto{\mOm}
      \SkSeq_{\infty} \xto{\mOm} \SkSeq_{\infty}
       \to \cdots \right) \]
  and refer to its objects as
  \emph{left-rigid $(\infty,\infty)$-prespectra}.
  Restricting to symmetric monoidal \iicats{} we get a full subcategory
  \[ \SpIIl := \lim \left(\cdots \to \SMCatII^{\lrig} \xto{\mOm}
      \SMCatII^{\lrig} \xto{\mOm}   \SMCatII^{\lrig} \to \cdots \right) \]
  of \emph{left-rigid $(\infty,\infty)$-spectra}.
  We denote the projection from each of these limits to the $n$th term in the
  diagram by
  \[ \mOm^{\infty-n} \colon \PSpIIl \to \SkSeq_{\infty}, \qquad \mOm^{\infty-n}
    \colon \SpIIl \to \SMCatII^{\lrig}.\]
\end{defn}

\begin{observation}
  A left-rigid $(\infty,\infty)$-prespectrum $\spX$ consists
  of a sequence \[\spX^{(-n)} = \mOm^{\infty-n}\spX = (\spX^{(-n)}_{1} \to \spX^{(-n)}_{2} \to
    \dots) \in \SkSeqI,\]
  $n = 0,1,\ldots$,
  such that $\mOm \spX^{(-n-1)} = \spX^{(-n)}$, \ie{} $\spX^{(-n)}_{k}
  = \mOm\spX^{(-n-1)}_{k+1}$.
\end{observation}

\begin{remark}
  Although this definition is convenient for our purposes, note that
  the symmetric monoidal structures in the definition of left-rigid
  $(\infty,\infty)$-(pre)spectra are actually redundant: If we keep
  the notation $\mOm$ for the functor
  $\Cat_{(\infty,\infty),*} \to \Cat_{(\infty,\infty),*}$ that takes a
  pointed \iicat{} $(\nC,x)$ to the endomorphisms $(\nC(x,x), \id_{x})$ then
  $\mOm^{n}(\nC,x)$ is automatically an $E_{n}$-monoidal \iicat{}. Thus
  all objects in an infinite delooping sequence of \iicats{} have
  canonical symmetric monoidal structures.  
\end{remark}

\begin{remark}
  Similar notions of categorified spectra, obtained by delooping \iicats{}, have been
  studied in the theses of Stefanich~\cite[Ch.~13]{Stefanich} and
  Masuda~\cite{Masuda}; both consider the more general variant of our
  definition that does not assume any dualizability (and so includes
  objects that do not correspond to categorical chain
  complexes). Heine~\cite{HeineSp} has recently set up an analogue of
  the theory of stable \icats{} in this context and proved
  $(\infty,\infty)$-versions of a number of fundamental results in stable
  homotopy theory.
\end{remark}

\begin{defn}
  For chain complexes, we have functors $t \colon \CChp_{n} \to
  \CChp_{n-1}$ that forget the bottom layer in the chain
  complex. These are compatible via commutative squares
  \[
    \begin{tikzcd}
     \CChp_{n} \ar[r, "t"] \ar[d, "b"'] & \CChp_{n-1}  \ar[d, "b"] \\
     \CChp_{n-1} \ar[r, "t"'] & \CChp_{n-2}
    \end{tikzcd}
  \]
  and so in the limit induce a functor
  \[ t \colon \CChp_{\infty} \to \CChp_{\infty}\]
  that again forgets the bottom layer.
  The functor $t$ restricts to the full subcategories of complete
  objects in the finite case, hence in the limit we also have
  \[ t \colon \CChpc_{\infty} \to \CChpc_{\infty}.\]  
  We define the \icat{} of \emph{(unbounded) categorical chain
    complexes} by
    \[ \CCh := \lim \left(\cdots \to \CChp_{\infty} \xto{t}
      \CChp_{\infty} \xto{t}
      \CChp_{\infty} \to \cdots \right) \]
  and its full subcategory of complete objects by
    \[ \CChc := \lim \left(\cdots \to \CChpc_{\infty} \xto{t}
      \CChpc_{\infty} \xto{t}
      \CChpc_{\infty} \to \cdots \right); \]
  we denote by $t^{\infty-n} \colon \CCh \to \CChp_{\infty}$ the
  projection to the $n$th term in the diagram.
\end{defn}

\begin{propn}\label{propn:omega top cube}
  The equivalences of \cref{cor:skseq cch eq} fit in commutative
  cubes
  \[
    \begin{tikzcd}[row sep=small,column sep=small]
      \SkSeq_{n} \ar[rr, "\sim"] \ar[dr, "p"] \ar[dd, "\mOm"] & &
      \CChp_{n-1} \ar[dr, "b"] \ar[dd, "t"{near start}] \\
      & \SkSeq_{n-1} \ar[rr,crossing over, "\sim"{near start}]  & & \CChp_{n-2}
      \ar[dd, "t"] \\
      \SkSeq_{n-1} \ar[rr, "\sim"] \ar[dr, "p"] & & \CChp_{n-2} \ar[dr, "b"] \\
      &  \SkSeq_{n-2} \ar[rr, "\sim"] \ar[uu,leftarrow,crossing over,
      "\mOm"{near end}]  & & \CChp_{n-3}.
    \end{tikzcd}
  \]
\end{propn}
\begin{proof}
  Unpacking the construction of these equivalences, we see that we
  have, firstly, a cube
  \[
    \begin{tikzcd}[row sep=small,column sep=small]
      \SkSeq_{n} \ar[rr, "\sim"] \ar[dr, "p"] \ar[dd, "\mOm"] & &
      \SkSeq_{n-1} \times_{\CChp_{0}} \CChp_{1} \ar[dr, "p
      \times_{\id} \id"] \ar[dd, "\text{proj}"{near start}] \\
      & \SkSeq_{n-1} \ar[rr,crossing over, "\sim"{near start}]  & &
      \SkSeq_{n-2} \times_{\CChp_{0}} \CChp_{1}
      \ar[dd, "\text{proj}"] \\
      \SkSeq_{n-1} \ar[rr, "="] \ar[dr, "p"] & & \SkSeq_{n-1} \ar[dr, "p"] \\
      &  \SkSeq_{n-2} \ar[rr, "="] \ar[uu,leftarrow,crossing over,
      "\mOm"{near end}]  & & \SkSeq_{n-2},
    \end{tikzcd}
  \]
  since the equivalence of \cref{propn:1-step m-rig pinned} is defined using
  $\mOm$. Next, the inductive construction of the equivalence with
  chain complexes gives a commutative cube
    \[
    \begin{tikzcd}[row sep=small,column sep=small]
      \SkSeq_{n-1} \times_{\CChp_{0}} \CChp_{1} \ar[rr, "\sim"]
      \ar[dr, "p \times_{\id}\id"] \ar[dd, "\text{proj}"] & &
      \CChp_{n-1} \ar[dr, "b"] \ar[dd, "t"{near start}] \\
      & \SkSeq_{n-2} \times_{\CChp_{0}} \CChp_{1} \ar[rr,crossing over, "\sim"{near start}]  & & \CChp_{n-2}
      \ar[dd, "t"] \\
      \SkSeq_{n-1} \ar[rr, "\sim"] \ar[dr, "p"] & & \CChp_{n-2} \ar[dr, "b"] \\
      &  \SkSeq_{n-2} \ar[rr, "\sim"] \ar[uu,leftarrow,crossing over,
      "\text{proj}"{near end}]  & & \CChp_{n-3},
    \end{tikzcd}
  \]
  where the top square is a pullback of the bottom. Putting these
  together gives the cube we want.
\end{proof}

In the limit, the cubes of \cref{propn:omega top cube} give the following:
\begin{cornoproof}
  The equivalence of \cref{cor:unfolding infinite} fits in a
  commutative square
    \[
    \begin{tikzcd}
      \SkSeq_{\infty} \ar[r] \ar[d,
   "\mOm"] & \CChp_{\infty}  \ar[d, "t"] \\
     \SkSeq_{\infty} \ar[r] & \CChp_{\infty}. 
    \end{tikzcd} \qedhere
  \]
\end{cornoproof}

Taking the limit over the vertical functors in this square, we then
get:
\begin{cornoproof}\label{cor:unfolding unbdd}
  The equivalences of \cref{cor:unfolding infinite} induce
  equivalences
  \[ \PSpIIl \simeq \CCh, \qquad \SpIIl \simeq \CChc.\]
  Moreover, these fit in commutative squares
  \[
    \begin{tikzcd}
     \PSpIIl \ar[r, "\sim"] \ar[d, "\mOm^{\infty-n}"'] & \CCh  \ar[d, "t^{\infty-n}"] \\
     \SkSeq_{\infty} \ar[r, "\sim"'] & \CChp_{\infty},
   \end{tikzcd}
   \quad
    \begin{tikzcd}
     \SpIIl \ar[r, "\sim"] \ar[d, "\mOm^{\infty-n}"'] & \CChc  \ar[d, "t^{\infty-n}"] \\
     \SMCatII^{\lrig} \ar[r, "\sim"'] & \CChpc_{\infty}
   \end{tikzcd}   
 \]
 with the previous equivalences.
\end{cornoproof}

We now want to show that we can embed $\CChp_{\infty}$ into $\CCh$ and
$\SkSeq_{\infty}$ into $\PSpIIl$, compatibly with the equivalences
between them. This is easy to set up for chain complexes, where we
just need to add zeros at the bottom:
\begin{observation}
  The functor $t \colon \CChp_{n} \to
  \CChp_{n-1}$ has a fully faithful left adjoint $\sigma$, given at an
  object $C$ by adding
  the pullback square
  \[
    \begin{tikzcd}
     \oC_{0} \ar[r, "="] \ar[d] & \oC_{0}  \ar[d] \\
     0 \ar[r] & 0
    \end{tikzcd}
  \]
  at the bottom; this further satisfies $t \sigma \simeq \id$. In the
  limit, the functors $\sigma$ induce a fully faithful functor $\sigma \colon
  \CChp_{\infty} \to \CChp_{\infty}$, left adjoint to $t$, and still
  satisfying $t\sigma \simeq \id$. We thus have a sequence of
  commutative squares
  \[
    \begin{tikzcd}
     \CChp_{\infty} \ar[r, hookrightarrow, "\sigma^{n}"] \ar[d, "="'] & \CChp_{\infty}  \ar[d, "t"] \\
     \CChp_{\infty} \ar[r, hookrightarrow, "\sigma^{n-1}"'] & \CChp_{\infty}.
    \end{tikzcd}
  \]
  In the limit we therefore get a fully faithful functor
  \[ \sigma^{\infty} \colon \CChp_{\infty} \hookrightarrow \CCh,\]
  which just extends a non-negatively graded chain complex down by
  adding $0$s.
\end{observation}

For skeletal sequences, the analogous embedding uses the canonical
delooping of a monoidal $(\infty,n)$-category to an
$(\infty,n+1)$-category with a single object:
\begin{observation}
  From \cite{enr}*{Theorem 6.3.2} we know that there is an adjunction
  \[ \mB : \Alg(\CatIn) \rightleftarrows \Cat_{(\infty,n+1),*} : \mOm\]
  between monoidal \incats{} and pointed $(\infty,n+1)$-categories,
  where $\mB$ forms the 1-object $(\infty,n+1)$-category that deloops
  a monoidal $(\infty,n)$-category and $\mOm$ takes the endomorphism
  $(\infty,n)$-category of the specified object, with monoidal
  structure given by composition. Here $\mB$ is fully
  faithful. Moreover, it is strong symmetric
  monoidal, so we get an induced adjunction
  \[ \mB : \SMCatIn \rightleftarrows \SMCat_{(\infty,n+1)} :
    \mOm.\]
  This restricts on full subcategories to an adjunction
  \[ \mB : \SMCatIn^{\lrig} \rightleftarrows \SMCat^{\lrig}_{(\infty,n+1)} :
    \mOm,\]
  since $\mB\nC$ always has duals for objects (as the only object is
  the symmetric monoidal unit). 
\end{observation}

\begin{propn}
  The functor $\mOm \colon \SkSeq_{n+1} \to \SkSeq_{n}$ has a fully faithful left adjoint
  \[\mB \colon \SkSeq_{n} \to \SkSeq_{n+1},\]
  which takes $\nA_{1} \to \cdots \to \nA_{n}$ to
  \[ 0 \to \mB \nA_{1} \to \cdots \to \mB \nA_{n}. \]
  Moreover, the mate of the commutative square \cref{obs:Omega bottom
    commute} is a commutative square
    \[
    \begin{tikzcd}
     \SkSeq_{n-1} \ar[r, "\mB"] \ar[d, "p"'] & \SkSeq_{n}  \ar[d, "p"] \\
     \SkSeq_{n-2} \ar[r, "\mB"'] & \SkSeq_{n-1}.
    \end{tikzcd}
  \]
\end{propn}
\begin{proof}
  We can regard $\SkSeq_{n+1}$ as a full subcategory of $\Fun([n],
  \SMCat_{(\infty,n+1)})$, and the functor $\mOm$ as the restriction
  to full subcategories of
  \[ \Fun([n], \SMCat_{(\infty,n+1)}) \xto{d_{0}^{*}} \Fun([n-1],
    \SMCat_{(\infty,n+1)}) \xto{\mOm_{*}} \Fun([n-1],\SMCatIn).\]
  This then has a left adjoint given by
  \[ \Fun([n-1],\SMCatIn) \xto{\mB_{*}}
    \Fun([n-1],\SMCat_{(\infty,n+1)})
    \xto{d_{0,!}} \Fun([n], \SMCat_{(\infty,n+1)})
  \]
  where the left Kan extension $d_{0,!}$ is given by inserting the
  initial object $0 \in \SMCat_{(\infty,n+1)}$ in the first
  position. Both of these functors are fully faithful, and since the
  composite clearly takes the full subcategory $\SkSeq_{n}$ into
  $\SkSeq_{n+1}$ this restricts to an adjoint to $\mOm$ on skeletal
  sequences, as required. To see that the mate of the square from
  \cref{obs:Omega bottom commute} commutes, it suffices to observe
  that the Beck--Chevalley transformation is obtained by restricting
  that for the composite square in \cref{eq:Omega bottom commute diagram},
  and the mates of both squares in the diagram clearly commute.
\end{proof}

This shows that the functor $\mB$ gives a natural transformation
$\SkSeq_{\bullet} \to \SkSeq_{\bullet+1}$; in the limit, we therefore get:
\begin{cornoproof}
  The functor $\mOm \colon \SkSeq_{\infty} \to \SkSeq_{\infty}$ has a
  fully faithful left adjoint
  \[\mB \colon \SkSeq_{\infty} \to \SkSeq_{\infty}. \qedhere\]
\end{cornoproof}

Now we observe that iterating the functor $\mB$ on $\SkSeq_{\infty}$
gives an embedding into $\PSpIIl$:
\begin{cor}
  The functor $\mOm^{\infty} \colon \PSpIIl \to \SkSeq_{\infty}$ has a fully faithful left adjoint
  $\mB^{\infty}$, so that $\mOm^{\infty-n}\mB^{\infty} \simeq \mB^{n}$.
\end{cor}
\begin{proof}
  The functor $\mOm^{\infty}$ is the functor on limits induced by the
  natural transformation
  \[
    \begin{tikzcd}
      \SkSeq_{\infty} \ar[r, "\mOm^{n+1}"] \ar[d,
      "\mOm"'] & \SkSeq_{\infty}   \ar[d, "="] \\
      \SkSeq_{\infty}\ar[r, "\mOm^{n}"'] & \SkSeq_{\infty}.
    \end{tikzcd}
  \]
  Taking horizontal mates gives a natural transformation of left
  adjoints, since $\mOm\mB^{n+1} \simeq \mB^{n}$. We therefore get an
  adjunction in the limit, as required.
\end{proof}

\begin{propn}
  The equivalences of \cref{cor:unfolding infinite} and
  \cref{cor:unfolding unbdd} fit in a commutative square
  \[
    \begin{tikzcd}
     \CChp_{\infty} \ar[r, "\sim"] \ar[d, hookrightarrow, "\sigma^{\infty}"] &
     \SkSeq_{\infty}  \ar[d, hookrightarrow, "\mB^{\infty}"] \\
     \CCh \ar[r, "\sim"'] & \PSpIIl.
    \end{tikzcd}
  \]
\end{propn}
\begin{proof}
  We already have an equivalence between the corresponding right
  adjoints from \cref{cor:unfolding unbdd}, so the left adjoints must
  also correspond.
\end{proof}

\begin{warning}
  The functors \[\mB \colon \SkSeq_{\infty} \to \SkSeq_{\infty}, \quad
  \sigma \colon \CChp_{\infty} \to \CChp_{\infty}\] do \emph{not}
  restrict to the full subcategories $\SMCatII^{\lrig}$ or $\CChpc_{\infty}$.
  The left adjoints of the functors \[\mOm \colon \SMCatII^{\lrig}
  \to \SMCatII^{\lrig}, \quad t \colon \CChpc_{\infty} \to
  \CChpc_{\infty}\] are therefore given by applying $\mB$ and $\sigma$, respectively,
  and then localizing to these full subcategories.
\end{warning}

\subsection{Completion}
\label{sec:completion}

In this subsection we will give an inductive description of the
localization functor
\[ \CChp_{n} \to \CChpc_{n}\]
for finite $n$, and then extend this to the infinite-length and unbounded
cases.

\begin{notation}
  For $k = 0,1,\dots,n$, let $\CChpcgr{k}_{n}$ denote the full subcategory of $\CChp_{n}$
  consisting of chain complexes \cref{eq:ch cx pb}
  such that $\oC_{i} \to \oZ_{i-1}$ is complete for $i > k$; then
  we have
  \[ \CChpc_{n} = \CChpcgr{0}_{n} \subseteq 
    \cdots \subseteq \CChpcgr{n}_{n} = \CChp_{n}.\]
  Similarly, we define $\SkSeq_{n+1}^{k\dcomp} \subseteq
  \SkSeq_{n+1}$ to be the full
  subcategory on sequences
  \[ \nA_{1} \to \cdots \to \nA_{n} \to \nA_{n+1}\]
  such that $\nA_{i} \to \core[i] \nA_{i+1}$ is an equivalence for $i
  > k$, so that we have
  \[ \SMCat_{(\infty,n+1)}^{\lrig} = \SkSeq_{n+1}^{0\dcomp} \subseteq \cdots
    \subseteq \SkSeq_{n+1}^{n\dcomp} = \SkSeq_{n+1}.\]
\end{notation}

\begin{observation}
  The equivalence of \cref{cor:unfold pinned} restricts on full
  subcategories to equivalences
  \[ \SkSeq_{n+1}^{k\dcomp} \simeq
    \CChpcgr{k}_{n} \]
  for $0 \leq k \leq n$.
\end{observation}

The inclusion
$\SkSeq_{n+1}^{(k-1)\dcomp} \hookrightarrow  \SkSeq_{n+1}^{k\dcomp}$ has an obvious left
adjoint, given by replacing  $\nA_{k}$ in a sequence
\[ \nA_{1} \to \cdots \to \nA_{n} \to \nA_{n+1}\]
by $\core[k]\nA_{k+1}$. We now want to describe the corresponding
adjoint on the chain complex side, starting with the case $n = 2$:
\begin{observation}
  Consider an object of $\CChpcgr{2}_{2} = \CChp_{2}$, that is a chain complex
  \[
    \begin{tikzcd}
      \oZ_{2} \ar[r] \ar[d] \drpullback & \oC_{2} \ar[d] \\
      0 \ar[r] & \oZ_{1} \ar[d] \ar[r]\drpullback  & \oC_{1} \ar[d] \\
      &  0 \ar[r] & \oZ_{0}
    \end{tikzcd}
  \]
  of length $2$. This
  corresponds to a 2-pinned $(\infty,3)$-category
  \[ \oZ_{0} \to \nB \to \nC\]
  where
  \[ \oC_{1} = \mP\nB|_{\oZ_{0}},
    \quad \oZ_{1} = \mOm\nB, \quad \oC_{2} = \mP\mOm\nC|_{\mOm\nB}, \quad \oZ_{2} = \mOm^{2}\nC.\]
  The corresponding $1$-complete $2$-pinned $(\infty,3)$-category is
  \[ \oZ_{0} \to \core[2]\nC \to \nC,\]
  with corresponding chain complex
  \[
    \begin{tikzcd}
      \oZ_{2} \ar[r] \ar[d] \drpullback & \oC'_{2} \ar[d] \\
      0 \ar[r] & \oZ'_{1} \ar[d] \ar[r]\drpullback  & \oC'_{1} \ar[d] \\
      &  0 \ar[r] & \oZ_{0}
    \end{tikzcd}
  \]
  where
  \[ \oC'_{1} = \mP\core[2]\nC|_{\oZ_{0}} = (\core \mP\nC)|_{\oZ_{0}},  
    \quad \oZ'_{1} = \mOm \core[2] \nC = \core\mOm\nC, \quad \oC'_{2} = \mP \mOm \nC.\]
  In particular, this localization does \emph{not} change the end terms
  $\oZ_{2}$ and $\oZ_{0}$ of the chain complex. Since, for $2 \leq k
  \leq n$, we have decompositions
  \[
    \begin{split}
      \CChpcgr{k}_{n} & \simeq \CChpc_{n-k}
                        \times_{\CChp_{0}} \CChp_{k} \\
       & \simeq \CChpc_{n-k}
      \times_{\CChp_{0}} \CChp_{2} \times_{\CChp_{0}} \CChp_{k-2},
    \end{split}
  \]
  the localization $\CChpcgr{k}_{n} \to \CChpcgr{k-1}_{n}$ is
  obtained by applying the localization $\CChp_{2} \to
  \CChpcgr{1}_{2}$ to the middle factor. The final step is induced by
  $\CChp_{1} \to \CChpc_{1}$. We thus obtain the localization
  $\CChp_{n} \to \CChpc_{n}$ by inductively completing each map in
  the chain complex from the top down.
\end{observation}

Now we consider completion in $\SkSeq_{\infty}$ and $\CChp_{\infty}$,
using the following lemma. This is a special case of
\cite[Theorem~A]{GeissCofiltered}; we include the short argument for
convenience.
\begin{lemma}\label{lem:seq of locns}
  Suppose $\oC$ is an \icat{} with sequential colimits and a nested sequence of
  subcategories
  \[ \dots \subseteq \oC_{n+1} \subseteq \oC_{n} \subseteq \cdots
    \subseteq \oC_{1} \subseteq \oC\]
  such that each $\oC_{n}$ is closed under sequential colimits and the
  inclusion $\oC_{n} \hookrightarrow \oC$ has a left adjoint
  $L_{n}$. If $\oC_{\infty} = \lim_{n \rightarrow \infty} \oC_{n}$
  denotes the intersection of these subcategories, then the inclusion
  $\oC_{\infty} \hookrightarrow \oC$ also has a left adjoint
  $L_{\infty}$, given by the colimit
  \[ L_{\infty}x = \colim \left(L_{1}x \to L_{2}x \to \cdots\right).\]
\end{lemma}
\begin{proof}
  Since $\oC_{n}$ is closed under sequential colimits, we see that
  $L_{\infty}x$ is contained in $\oC_{n}$ (as the terms $L_{k}x$ for
  $k \geq n$ lie in $\oC_{n}$). As this holds for all $n$, the object
  $L_{\infty}x$ lies in $\oC_{\infty}$. Moreover, for $y \in
  \oC_{\infty}$ we have
  \[ \oC_{\infty}(L_{\infty}x,y) \simeq \limP_{n \rightarrow
      \infty}\oC(L_{n}x, y) \simeq \limP_{n \rightarrow
      \infty}\oC(x, y) \simeq \oC(x,y) \]
  as the limit is constant, which gives the required adjunction.
\end{proof}

\begin{notation}
  Let $\SkSeqIcl{k}$ denote the full subcategory of $\SkSeqI$
  containing the skeletal sequences $\nA_{1} \to \nA_{2}
  \to \cdots$ where $\nA_{i} \to \core[i] \nA_{i+1}$ is an equivalence for $i
  \leq k$, so that we have $\SkSeqIcl{k+1} \subseteq \SkSeqIcl{k}$ with intersection
  \[ \SMCatII^{\lrig} = \lim_{k \rightarrow \infty}\SkSeqIcl{k} \]
  and
   $\SkSeqIcl{0} = \SkSeqI$.
\end{notation}

\begin{observation}
  Each full subcategory $\SkSeqIcl{k}$ is closed under sequential
  colimits in $\SkSeqI$ (as $\core[i]$ preserves filtered colimits for
  all $i$), and the
  inclusion $\SkSeqIcl{k} \hookrightarrow \SkSeqI$ has a left adjoint
  \[ \gamma_{[0,k]} \colon \SkSeqI \to \SkSeqIcl{k},\]
  which puts the cores of $\nA_{k+1}$ in the $k$ lowest degrees.
\end{observation}

\begin{propn}\label{propn:complete iicat}
  The inclusion $\SMCatII^{\lrig} \hookrightarrow \SkSeq_{\infty}$ has
  a left adjoint $\gamma_{\geq 0}$, which takes a skeletal sequence
  $\nA_{1} \to \nA_{2} \to \cdots$ to the symmetric monoidal \iicat{}
  $\gamma_{\geq 0}(\nA_{\bullet})$ with $(\infty,n)$-component
  \[ \gamma_{\geq 0}(\nA_{\bullet})_{n} = \colimP_{k \geq n} \core[n]\nA_{k}. \]
\end{propn}
\begin{proof}
 This follows from \cref{lem:seq of locns}.
\end{proof}

\begin{observation}
  Applying the same lemma similarly to chain complexes, we see that
  the localization $\CChp_{\infty} \to \CChpc_{\infty}$ is computed by
  a sequential colimit over the chain complexes obtained by completing
  in degrees $< k$ for all $k$.
\end{observation}

\begin{propn}
  There is a commutative square
  \[
    \begin{tikzcd}
     \SkSeqI \ar[r, "\gamma_{\geq 0}"] \ar[d, "\mOm"'] & \SMCatII^{\lrig}  \ar[d, "\mOm"] \\
     \SkSeqI \ar[r, "\gamma_{\geq 0}"'] & \SMCatII^{\lrig},
    \end{tikzcd}
  \]
  obtained as the horizontal mate of the commutative square
  restricting $\mOm$ from $\SkSeqI$ to the full subcategory $\SMCatII^{\lrig}$.
\end{propn}
\begin{proof}
  We must show that the Beck--Chevalley transformation $\gamma_{\geq 0} \mOm
  \to \mOm \gamma_{\geq 0}$ is invertible. At an object $\nA_{\bullet} \in \SkSeqI$,
  this is given in its $(\infty,n)$-component by
  \[
    \begin{split}
      (\gamma_{\geq 0} \mOm \nA_{\bullet})_{n}
      & = \colimP_{k \geq n} \core[n] \mOm(\nA_{k+1}) \\
      & \simeq \colimP_{k \geq n} \mOm \core[n+1] \nA_{k+1} \\
      & \longrightarrow \mOm \left(\colimP_{k \geq n}
        \core[n+1] \nA_{k+1}\right) \\
      & = \mOm(\gamma_{\geq 0}\nA_{\bullet})_{n}.
    \end{split}
  \]
  This is an equivalence since $\mOm$ preserves sequential colimits.
\end{proof}

We can thus take the limit along $\mOm$ and obtain:
\begin{cornoproof}\label{cor:gamma psp geq 0}
  The fully faithful inclusion $\SpIIl \hookrightarrow \PSpIIl$ has a
  left adjoint $\gamma$, given for all $n$ by
  \[ \mOm^{\infty-n}\gamma\spX = \gamma_{\geq
      0}(\mOm^{\infty-n}\spX). \qedhere\]
\end{cornoproof}

We can also describe this localization a bit more concretely using
\cref{lem:seq of locns} again:
\begin{defn}
  Let $\PSpIIlcl{k}$ denote the full subcategory of $\PSpIIl$ spanned
  by the objects $\spX$ such that $\mOm^{\infty-i}\spX$ lies in
  $\SkSeqIcl{k+i}$ for all $i$ (where we interpret $\SkSeqIcl{j}$ as
  all of $\SkSeqI$ for $j \leq 0$).
\end{defn}

\begin{observation}
  We have a commutative square
  \[
    \begin{tikzcd}
     \SkSeqI \ar[r, "\gamma_{[0,k+1]}"] \ar[d, "\mOm"'] & \SkSeqIcl{k+1}  \ar[d, "\mOm"] \\
     \SkSeqI \ar[r, "\gamma_{[0,k]}"'] & \SkSeqIcl{k}
    \end{tikzcd}
  \]
  and so in the limit the full subcategory inclusion $\PSpIIlcl{k}
  \hookrightarrow \PSpIIl$ has a left adjoint $\gamma_{\leq k}$, given
  by
  \[ \mOm^{\infty-n}\gamma_{\leq k}\spX = \gamma_{[0,k+n]} \mOm^{\infty-n}\spX.\]
\end{observation}

\begin{propn}\label{propn:locn to iisp colim}
  The localization $\gamma \colon \PSpIIl \to \SpIIl$ can be computed
  as
  \[ \gamma \spX \simeq \colimP_{k \rightarrow \infty} \gamma_{\leq k}\spX.\]
\end{propn}
\begin{proof}
  The only thing needed to apply \cref{lem:seq of locns} is to observe
  that the full subcategory $\PSpIIlcl{k}$ is closed under sequential
  colimits, which again holds since $\core[n]$ always preserves
  filtered colimits.
\end{proof}

\begin{observation}
  Transporting \cref{propn:locn to iisp colim} to
  categorical chain complexes, we see that the localization
  \[ \CCh \to \CChc\]
  is given by taking the sequential colimit over $n$ of the chain complex obtained
  by completing downwards starting in level $n$.
\end{observation}

\subsection{Chain complexes of $\infty$-groupoids and localization}
\label{sec:gpdchcx}

In this section we first discuss how the unfolding equivalence
restricts to chain complexes of \igpds{} and describe how the general
case localizes to this one. In particular, we obtain a description of
the spectrum obtained by localizing an $(\infty,\infty)$-(pre)spectrum.

\begin{defn}
  For $n \leq \infty$, we write $\GpdChp_{n}$ for the full subcategory
  of $\CChp_{n}$ spanned by those chain complexes of symmetric
  monoidal \icats{} whose terms are all \igpds{}, and $\GpdChpc_{n}$
  for the further full subcategory of such chain complexes that are
  also complete. Similarly, we write $\GpdCh \subseteq \CCh$, $\GpdChc
  \subseteq \CChc$ for the full subcategories of (complete) unbounded
  chain complexes consisting of \igpds{}.
\end{defn}

\begin{defn}
  A $0$-pinning of a grouplike $E_{\infty}$-space $C$ is, by
  definition, no additional data.  For $k\geq1$, a
  \emph{$k$-pinning} of $C$ is a
  sequence
  \[ A_{1} \to \cdots \to A_{k} \to C\] of morphisms of grouplike
  $E_{\infty}$-spaces, where $A_{i} \to A_{i+1}$ is
  $(i-1)$-surjective\footnote{Recall from \cref{rmk:k-surj gpd} that
    for maps of \igpds{} this means the map is an isomorphism on all
    homotopy groups in degree $< i-1$ and surjective in degree $i-1$.}
  and $A_{k} \to C$ is $(k-1)$-surjective.  We write $\CGrp(\Spc)$ for
  the \icat{} of grouplike $E_{\infty}$-spaces, and
  $\PCGrp[k] \subseteq \Fun([k], \CGrp(\Spc))$ for the full
  subcategory of \emph{$k$-pinned} grouplike $E_{\infty}$-spaces,
  meaning grouplike $E_{\infty}$-spaces equipped with a $k$-pinning.
  In particular, $\PCGrp[0]=\CGrp(\Spc)$.
\end{defn}

\begin{observation}
  A symmetric monoidal \igpd{} has duals for objects \IFF{} it is
  grouplike, and always has all adjoints. Since loop spaces are always
  grouplike, this means that the \icat{}
  $\PCGrp[k]$ is a full subcategory of $\PSMCatIn[k]^{m\dlrig}$ for
  all $n$ and $m$.
\end{observation}

\begin{propn}\label{propn:pinned gpd ch cx eq}
  The equivalence of \cref{cor:skseq cch eq} restricts to an
  equivalence
  \[ \PCGrp[n] \simeq \GpdChp_{n}\]
  for $0 \leq n \leq \infty$.
\end{propn}

The proof boils down to the following observation:
\begin{lemma}\label{propn:compgpd}
  Suppose $\pi \colon C \to Z$ is a morphism of symmetric monoidal
  \igpds{} with duals (\ie{} that are grouplike). If $p \colon Z \to
  \Spc$ is the lax symmetric monoidal functor obtained by
  straightening $\pi$, then the \icat{} $p_{*}\overline{Z}$ is an
  \igpd{}.
\end{lemma}
\begin{proof}
  Since $C$ has duals, we know from \cref{propn:PCdualadj} that the
  \icat{} $p_{*}\overline{Z}$ must be left 2-rigid when regarded as an
  \itcat{}. Then \cref{propn:left 2-rigid cat is gpd} shows that it is
  in fact an \igpd{}.
\end{proof}

\begin{remark}
  If we drop the assumption that $C$ is grouplike in
  \cref{propn:compgpd}, then $p_{*}\overline{Z}$ may no longer be an
  \igpd{}. For example, for $Z = *$ the map $p \colon C \to *$
  produces the \icat{} $BC$ with a single object that has $C$ as its
  endomorphisms; if $C$ is not grouplike then this is not an \igpd{}.
\end{remark}

\begin{proof}[Proof of \cref{propn:pinned gpd ch cx eq}]
  It is clear that an object of $\PCGrp[n]$ corresponds to a chain
  complex of \igpds{}. The converse follows by applying
  \cref{propn:compgpd} iteratively.
\end{proof}

\begin{notation}
  Let $\PSpIIgpd$ denote the limit
  \[ \lim \left( \cdots \to \PCGrp[\infty] \xto{\Omega} \PCGrp[\infty]
      \to \cdots
    \right)\]
  where $\Omega$ is the restriction of $\mOm \colon \SkSeqI \to
  \SkSeqI$, which is just given by forgetting the first term and
  taking loop spaces on the rest. We note that the intersection of
  $\PSpIIgpd$ with $\SpIIl$ is equivalent to the limit
    \[ \lim \left( \cdots \to \CGrp(\Spc) \xto{\Omega} \CGrp(\Spc)
      \to \cdots
    \right),\]
  which gives the \icat{} $\Sp$ of spectra.
\end{notation}

\begin{cornoproof}
  The equivalences of \cref{cor:unfolding unbdd} restrict to
  compatible equivalences
  \[
    \begin{tikzcd}
     \Sp \ar[r, "\sim"] \ar[d, hookrightarrow] & \GpdChc  \ar[d, hookrightarrow] \\
     \PSpIIgpd \ar[r, "\sim"'] & \GpdCh. \qedhere
    \end{tikzcd}
  \]
\end{cornoproof}

\begin{observation}
  By definition, a non-negatively graded chain complex of symmetric monoidal \igpds{}
  is complete \IFF{} for each degree $i$, the functor
  $(Z_{i-1})_{\bbone/} \to C_{i}$ exhibits the source as the
  underlying left fibration over $Z_{i-1}$ of the target. Since
  $C_{i}$ is a space, this amounts to requiring this functor to be an
  equivalence, and since $Z_{i-1}$ is also a space, the slice
  $(Z_{i-1})_{\bbone/}$ is contractible. A complete chain complex
  therefore consists of pullback squares of the form
  \[
    \begin{tikzcd}
     Z_{i} \ar[r] \ar[d] & 0  \ar[d] \\
     0 \ar[r] & Z_{i-1},
    \end{tikzcd}
  \]
  so that $Z_{i} \simeq \Omega Z_{i-1}$; the chain complex is
  therefore completely determined by $Z_{0} \simeq C_{0}$. We conclude that
  for every $n \leq \infty$ the \icat{} $\GpdChpc_{n}$ is
  equivalent to $\CGrp(\Spc)$.
\end{observation}

In summary, we have shown that the unfolding equivalence gives
identifications of commutative squares:
\[
  \begin{tikzcd}
   \CGrp(\Spc) \ar[r, hookrightarrow] \ar[d, hookrightarrow] &
   \SMCatIn^{\lrig}  \ar[d, hookrightarrow] \\
   \PCGrp[n-1] \ar[r, hookrightarrow] & \SkSeq_{n}
 \end{tikzcd}
 \quad\simeq\quad
   \begin{tikzcd}
   \GpdChpc_{n-1} \ar[r, hookrightarrow] \ar[d, hookrightarrow] &
   \CChpc_{n-1}  \ar[d, hookrightarrow] \\
   \GpdChp_{n-1} \ar[r, hookrightarrow] & \CChp_{n-1},
 \end{tikzcd}
\]
\[
  \begin{tikzcd}
   \Sp \ar[r, hookrightarrow] \ar[d, hookrightarrow] &
   \SpIIl  \ar[d, hookrightarrow] \\
   \PSpIIgpd \ar[r, hookrightarrow] & \PSpIIl
 \end{tikzcd}
 \quad\simeq\quad
   \begin{tikzcd}
   \GpdChc \ar[r, hookrightarrow] \ar[d, hookrightarrow] &
   \CChc  \ar[d, hookrightarrow] \\
   \GpdCh \ar[r, hookrightarrow] & \CCh.
 \end{tikzcd}
\]
All of these inclusions moreover have left adjoints. For the right
vertical maps these are the completion functors we studied in the
previous subsection. Note that for the left-hand squares it is clear
from the definitions of these functors that they restrict to the
groupoidal case; we can also see this for chain complexes by applying
\cref{propn:compgpd}. 

Let us now look at the left adjoints of the horizontal inclusions. We
will see that in the left-hand squares we can describe the top
horizontal left adjoints fairly explicitly,
but the bottom horizontal left adjoints are less explicit. On the
other hand, the latter admit concrete descriptions in the right-hand
squares. We start with
the chain complex case:
\begin{propn}
  All of the inclusions $\GpdChp_{n} \hookrightarrow \CChp_{n}$ for $n \leq \infty$
  and $\GpdCh \hookrightarrow \CCh$
  have left
  adjoints, which are simply given by
  applying the localization \[\|\blank\| \colon \CatI \to \Spc\]
  termwise. 
\end{propn}
\begin{proof}
  It suffices  to check that for each square
  \[
    \begin{tikzcd}
     \oZ_{n} \ar[r] \ar[d] & \oC_{n}  \ar[d, "\partial"] \\
     0 \ar[r] & \oZ_{n-1}
    \end{tikzcd}
  \]
  in a chain
  complex, we get a pullback square after applying $\|\blank\|$. Since
  $\partial$ is both a cartesian and cocartesian fibration, this
  follows from \cite[Theorem 2.3]{SteimleAdd}. 
\end{proof}

\begin{warning}
  The localization $\CChp_{n} \to \GpdChp_{n}$ does \emph{not}
  preserve complete chain complexes --- indeed, $\|\core\mP\nC\|$ is
  typically far from contractible! To describe the localization to
  $\GpdChpc_{n}$ we therefore have to complete the result; we will
  study this completion procedure in the groupoidal case in the next subsection.
\end{warning}

\begin{observation}\label{obs:smcatinlrig loc}
  Let
  \[ \|\blank\| \colon \CatIn \rightleftarrows \Spc \]
  be the left adjoint to the inclusion $\Spc \hookrightarrow
  \CatIn$. This preserves products, and so also gives a left adjoint
  \[ \SMCatIn \to \CMon(\Spc),\]
  which restricts on full subcategories to a left adjoint
  \[ \SMCatIn^{\lrig} \to \CGrp(\Spc)\]
  to the inclusion of $\CGrp(\Spc)$ into $\SMCatIn^{\lrig}$.
\end{observation}

\begin{propn}\label{lem:locn catii spc}
  The inclusion $\Spc \hookrightarrow \CatII$ has a left adjoint given
  by
  \[ \|\nC\| \simeq \colimP_{n\rightarrow \infty} \|\core[n] \nC\|,\]
  which also induces left adjoints of the inclusions $\CMon(\Spc)
  \hookrightarrow \SMCatII$ and $\CGrp(\Spc) \hookrightarrow \SMCatII^{\lrig}$.
\end{propn}
\begin{proof}
  The formula for $\|\nC\|$ clearly defines a functor. We have
  \[
    \begin{split}
      \Map_{\Spc}(\|\nC\|, X) & \simeq \limP_{n \rightarrow \infty} 
      \Map_{\Spc}(\|\core[n]\nC\|, X) \\ & \simeq \limP_{n \rightarrow \infty}
      \Map_{\CatIn}(\core[n]\nC, X) \\ & \simeq \Map_{\CatII}(\nC, X),
    \end{split}
  \]
  so this gives a left adjoint as required. Since the component
  functors on \incats{} preserve products, so does this functor, and
  we therefore get induced adjunctions on (left rigid) symmetric monoidal
  \iicats{} as in \cref{obs:smcatinlrig loc}.
\end{proof}

\begin{cor}\label{cor:locn skseq formula}
  The left adjoint of the inclusion $\Spc \hookrightarrow \SkSeqI$ is
  given by the formula
  \[ \|\gamma_{\geq 0}\nA_{\bullet}\| \simeq \colimP_{n \rightarrow \infty} \|\nA_{n}\|\]
\end{cor}
\begin{proof}
  Combining \cref{propn:complete iicat} and \cref{lem:locn catii spc},
  we have
  \[
    \begin{split}
      \|\gamma_{\geq 0}\nA_{\bullet}\| & \simeq \colimP_{n \rightarrow
                                   \infty} \| \core[n]\gamma_{\geq 0}
                                   \nA_{\bullet}\| \\
      & \simeq \colimP_{n \rightarrow \infty}\colimP_{m \geq n}
        \| \core[n]\nA_{m}\|. \\
    \end{split}
  \]
  Here the diagonal subsequence is cofinal, which gives the result.
\end{proof}

\begin{observation}
  Given $\nC \in \SMCatII^{\lrig}$, the canonical map $\nC \to
  \|\nC\|$ induces a map $\mOm \nC \to \mOm\|\nC\| \simeq \Omega
  \|\nC\|$. For $\spX \in \SpIIl$ we therefore have a prespectrum in
  $\CGrp(\Spc)$ with
  components $\|\mOm^{\infty-n}\spX\|$ and
  structure maps
   \begin{equation}
     \label{eq:prespm on realization}
     \|\mOm^{\infty-n}\spX\| \simeq \|\mOm \mOm^{\infty-n-1}\spX\| \to \Omega\|\mOm^{\infty-n-1}\spX\|.
   \end{equation}
\end{observation}

\begin{propn}\label{propn:L SpIIl formula}
  The inclusion $\Sp \hookrightarrow \SpIIl$ has a left adjoint
  $L_{\Sp}$, given by
  \[ L_{\Sp}\spX \simeq \colimP_{n \rightarrow \infty} \Sigma^{-n}
    \|\mOm^{\infty-n}\spX\| \simeq \colimP_{n,m \rightarrow \infty} \Sigma^{-n}
    \|\core[m] \mOm^{\infty-n}\spX\|, \]
  with the colimit computed in $\Sp$ using the structure maps \cref{eq:prespm on realization}.
\end{propn}
\begin{proof}
  The given formula for $L_{\Sp}$ clearly defines a functor $\SpIIl
  \to \Sp$. For a spectrum $E$ we have
  \[
    \begin{split}
      \Map_{\Sp}(L_{\Sp}\spX, E)
      & \simeq \limP_{n \rightarrow \infty}
        \Map_{\Sp}(\Sigma^{-n}\|\mOm^{\infty-n}\spX\|, E) \\
      & \simeq \limP_{n \rightarrow \infty}
        \Map_{\CGrp(\Spc)}(\|\mOm^{\infty-n}\spX\|,
        \Omega^{\infty-n}E)) \\
      & \simeq \limP_{n \rightarrow \infty}
        \Map_{\SMCatII^{\lrig}}(\mOm^{\infty-n}\spX,\Omega^{\infty-n}E))
      \\
      & \simeq \Map_{\SpIIl}(\spX, E),
    \end{split}
  \]
  as required.
\end{proof}

\begin{observation}
  We can also express the formula for $L_{\Sp}$ as
  \[ \Omega^{\infty-n}L_{\Sp}\spX \simeq \colimP_{k \rightarrow
      \infty} \Omega^{k}\| \mOm^{\infty-n-k} \spX\|. \]
\end{observation}

Combining \cref{propn:L SpIIl formula} with \cref{cor:gamma psp geq 0}
and \cref{cor:locn skseq formula}, we get:
\begin{cornoproof}
  The left adjoint of the inclusion $\Sp \hookrightarrow \PSpIIl$ is
  given by the formula
  \[ L_{\Sp} \gamma \spX \simeq \colimP_{n,m \rightarrow \infty}
    \Sigma^{-n}\|(\mOm^{\infty-n}\spX)_{m}\|. \qedhere\]
\end{cornoproof}

\subsection{Completion for $\infty$-groupoids and localization to
  spectra}
\label{sec:comp igpd}

In this subsection we give an alternative description of completion
for chain complexes of \igpds{} in terms of quotients in spectra. This
allows us to give a formula for the homotopy groups of the completion,
in a sense by taking the ``homology'' of the chain complex. We then
use this to describe the homotopy groups of the spectrum obtained by
localizing a left rigid $(\infty,\infty)$-(pre)spectrum in terms of
the associated categorical chain complex.

The completion functor for chain complexes of \igpds{} is obtained by
restricting the completion functor on categorical chain complexes (as
follows from \cref{propn:compgpd}, for instance). We conclude that it is
determined by iterating the completion $\GpdChp_{2} \to
\GpdChpcgr{1}_{2}$, which we can describe more concretely as follows:
\begin{propn}\label{propn:gpdchcomp}
  Given an object
  \begin{equation}
    \label{eq:len 2 ch cx gpd}
    \begin{tikzcd}
      Z_{2} \ar[r] \ar[d] \drpullback & C_{2} \ar[d] \\
      0 \ar[r] & Z_{1} \ar[d] \ar[r]\drpullback  & C_{1} \ar[d] \\
      &  0 \ar[r] & Z_{0}
    \end{tikzcd}
  \end{equation}
  in $\GpdChp_{2}$, its localization to $\GpdChpcgr{1}_{2}$ is
  \begin{equation}
    \label{eq:gpdcompch}
    \begin{tikzcd}
      Z_{2} \ar[r] \ar[d] \drpullback & 0 \ar[d] \\
      0 \ar[r] & Z_{1}/C_{2} \ar[d] \ar[r]\drpullback  & C_{1}/C_{2} \ar[d] \\
      &  0 \ar[r] & Z_{0}.
    \end{tikzcd}
  \end{equation}
\end{propn}

For the proof we need some simple observations about connective spectra:
\begin{lemma}\label{lem:connquot}
  Given a map of spectra $C \to X$ with $C$ connective, we have
  $\tau_{\geq 0}(X)/C \isoto \tau_{\geq 0}(X/C)$.
\end{lemma}
\begin{proof}
  Applying the long exact sequence of homotopy groups to the map of
  cofibre sequences of spectra
  \[
    \begin{tikzcd}
      C \ar[r] \ar[d, equals] & \tau_{\geq 0}X \ar[d] \ar[r]
        & \tau_{\geq 0}(X)/C \ar[d] \\
      C \ar[r] & X \ar[r] & X/C
    \end{tikzcd}
  \]
  we see that $\tau_{\geq 0}(X)/C \to X/C$ is an isomorphism on
  homotopy groups in dimension $\geq 0$.
\end{proof}

\begin{lemma}\label{lem:conn sp pullback}
  Given a pullback square in $\Sp^{\geq 0}$
  \[
    \begin{tikzcd}
     F \ar[r] \ar[d] \drpullback & Y  \ar[d, "f"] \\
     0 \ar[r] & X
    \end{tikzcd}
  \]
  and a map $C \to F$ of connective spectra, we get another pullback
  square
  \[
    \begin{tikzcd}
     F/C \ar[r] \ar[d] \drpullback & Y/C  \ar[d] \\
     0 \ar[r] & X
    \end{tikzcd}
  \]
  in $\Sp^{\geq 0}$.
\end{lemma}
\begin{proof}
  Let $Z$ be the fibre of $f$ in $\Sp$ (so $F \simeq \tau_{\geq
    0}Z$). In $\Sp$ we then have (using the 3-for-2 property of pushouts) a
  pullback square
  \[
    \begin{tikzcd}
     Z/C \ar[r] \ar[d] & Y/C  \ar[d] \\
     0 \ar[r] & X,
    \end{tikzcd}
  \]
  and truncating this then gives a pullback
  \[
    \begin{tikzcd}
     \tau_{\geq 0}(Z/C) \ar[r] \ar[d] & Y/C  \ar[d] \\
     0 \ar[r] & X,
    \end{tikzcd}
  \]
  in $\Sp^{\geq 0}$. Now \cref{lem:connquot} identifies this with the
  square we wanted to show was a pullback.
\end{proof}

\begin{proof}[Proof of \cref{propn:gpdchcomp}]
  We first observe that if there is a chain complex of the form
  \cref{eq:gpdcompch} then this must be the completion of the one we
  start with --- the universal property of the quotients implies that
  this diagram gives the localization from the larger \icat{} of
  all diagrams
  \[
    \begin{tikzcd}
      Y_{2} \ar[r] \ar[d] & B_{2} \ar[d] \\
      0 \ar[r] & Y_{1} \ar[d] \ar[r]  & B_{1} \ar[d] \\
      &  0 \ar[r] & Y_{0}
    \end{tikzcd}
  \]
  (with the squares not required to be pullbacks) to its subcategory
  with $B_{2} \simeq 0$. It thus suffices to show that the two squares
  in \cref{eq:gpdcompch} are indeed pullbacks.

  For the top square this follows because in the diagram
  \[
    \begin{tikzcd}
      Z_{2} \ar[r] \ar[d] \drpullback & C_{2} \ar[d] \ar[r] \arrow[phantom]{dr}[very near end,description]{\ulcorner} & 0 \ar[d]
      \\
      0 \ar[r]  & Z_{1} \ar[r] & Z_{1}/C_{2}
    \end{tikzcd}
  \]
  the right square is a pushout of connective spectra, and so also a
  pullback. The bottom square is a pullback by \cref{lem:conn sp
    pullback}.
\end{proof}

\begin{cor}\label{cor:complete htpy finite}
  Suppose $P_{1} \to P_{2} \to \cdots \to P_{n} \to Q$ in $\PCGrp[n]$ is an
  $n$-pinned grouplike $E_{\infty}$-space, with corresponding chain complex
  $C=(C_{\bullet}, Z_{\bullet}, \partial_{\bullet})$ in
  $\GpdChp_{n}$. Then the
  homotopy groups of $Q$ are given by
  \[ \pi_{i}Q \cong
    \begin{cases}
      \pi_{0}Z_{i}/\pi_{0}C_{i+1}, & 0 \leq i < n \\
      \pi_{i-n} Z_{n}, & i \geq n.
    \end{cases}
  \]
\end{cor}
\begin{proof}
  The description is evidently correct if the chain complex $C$ is
  complete, since then $Q = Z_{0}$, $C_{i} = 0$ for $i > 0$, and
  $Z_{i} \simeq \Omega^{i}Z_{0}$. It thus suffices to check these
  groups are unchanged at each step in our iterative completion
  construction. This reduces to showing that for a length-2 chain
  complex \cref{eq:len 2 ch cx gpd}, we have
  \[ \pi_{0}(Z_{0})/\pi_{0}(C_{1}) \cong \pi_{0}(Z_{0}/C_{1}) \cong
    \pi_{0}(Z_{0}/(C_{1}/C_{2})) \cong
    \pi_{0}(Z_{0})/\pi_{0}(C_{1}/C_{2}).
  \]
  Here the first and third isomorphisms are immediate since
  $\pi_{0} \colon \Sp^{\geq 0} \to \Ab$ is a left adjoint, while the
  second holds since the map $C_{1} \to Z_{0}$ factors through
  $C_{1}/C_{2}$ and
  $\pi_{0}C_{1} \to \pi_{0}(C_{1}/C_{2}) \cong
  \pi_{0}C_{1}/\pi_{0}C_{2}$ is surjective.
\end{proof}

\begin{cor}\label{cor:skseq gpd htpy}
  Suppose $P_{\bullet} \in \PCGrp[\infty]$ corresponds to the chain complex
  $C=(C_{\bullet}, Z_{\bullet}, \partial_{\bullet})$ in
  $\GpdChp_{\infty}$.
  Then the
  homotopy groups of $\gamma_{\geq 0} P_{\bullet}$ are given by
  \[ \pi_{i}\gamma_{\geq 0}P_{\bullet} \cong \pi_{0}Z_{i}/\pi_{0}C_{i+1}
  \]
  for all $i \geq 0$.
\end{cor}
\begin{proof}
  By \cref{propn:complete iicat} we have $\gamma_{\geq 0} P_{\bullet} \simeq
  \colimP_{n \rightarrow \infty} P_{n}$, and so
  \[ \pi_{i} \gamma_{\geq 0} P_{\bullet} \simeq \colimP_{n \rightarrow \infty}
    \pi_{i} P_{n}.\] Here we see from \cref{cor:complete htpy finite}
  that $\pi_{i} P_{n}$ stabilizes at $\pi_{0}Z_{i}/\pi_{0}C_{i+1}$ for
  $n > i$.
\end{proof}

\begin{cor}\label{cor:psp gpd htpy}
  Suppose $\spP \in \PSpIIgpd$
  corresponds to the chain complex
  $C=(C_{\bullet}, Z_{\bullet}, \partial_{\bullet})$ in
  $\GpdCh$.
  Then the
  homotopy groups of the spectrum $\gamma \spP$ are given by
  \[ \pi_{i}\gamma \spP \cong \pi_{0}Z_{i}/\pi_{0}C_{i+1}
  \]
  for all $i \in \mathbb{Z}$.
\end{cor}
\begin{proof}
  For a fixed $i \in \mathbb{Z}$, choose $n$ such that $i+n \geq 0$.
  We have
  \[
    \pi_{i}\gamma \spP = \pi_{i+n} \Omega^{\infty-n}\gamma\spP \simeq
    \pi_{i+n} \gamma_{\geq 0}(\Omega^{\infty-n}\spP).
  \]
  Now we just apply \cref{cor:skseq gpd htpy} to $\Omega^{\infty-n}\spP$,
  which corresponds to the chain complex obtained by shifting $C$ up
  by $n$ and truncating.
\end{proof}

\begin{defn}
  We say a chain complex \cref{eq:ch cx pb}
  $(\oC_{\bullet},\oZ_{\bullet},\partial_{\bullet})$ (of finite or
  infinite length, or unbounded) is
  \emph{object-complete} if for every $i$ the map $(\oZ_{i})_{\bbone/}
  \to \oC_{i+1}$ is surjective on objects.
\end{defn}

\begin{lemma}
  Suppose
  \[ \nA_{1} \to \nA_{2} \to \cdots \to \nA_{n} \]
  is a skeletal sequence $(n \leq \infty)$. If the functor $\nA_{i}
  \to \nA_{i+1}$ is $i$-surjective for all $i$, then the corresponding
  chain complex is object-complete.
\end{lemma}
\begin{proof}
  We note that if the given property holds for $\nA_{\bullet}$, then it also
  holds for $\mOm \nA_{\bullet}$. It therefore suffices to observe that if
  $\nA_{1} \to \nA_{2}$ is surjective on morphisms, then
  \[(\nA_{1})_{\bbone/} \to \mP \nA_{2}|_{\nA_{1}}\] is surjective on
  objects.
\end{proof}

Putting together our results on homotopy groups of completions with
the description of localizations in the previous subsection, we obtain
the following:
\begin{cor}\label{cor:pi of completion of locn}\ 
  \begin{enumerate}[(i)]
  \item\label{it:pi formula skseq finite} Suppose $\nA_{\bullet} \in \SkSeq_{n}$ ($n < \infty$) corresponds
    to the chain complex $\oC = (\oC_{\bullet}, \oZ_{\bullet},
    \partial_{\bullet}) \in \CChp_{n-1}$. Then
    \[ \pi_{i} \Lgpd{ \nA_{n}}
      \cong
      \begin{cases}
        \pi_{0}\Lgpd{\oZ_{i}}/\pi_{0}\Lgpd{\oC_{i+1}}, & 0 \leq i < n-1 \\
        \pi_{i-n+1} \Lgpd{\oZ_{n-1}}, & i \geq n-1.
      \end{cases}
    \]
  \item In the situation of \ref{it:pi formula skseq finite}, if the
    chain complex $\oC$ is object-complete, then the formula simplifies to
        \[ \pi_{i} \Lgpd{ \nA_{n}}
      \cong
      \begin{cases}
        \pi_{0}\Lgpd{\oZ_{i}}, & 0 \leq i < n-1 \\
        \pi_{i-n+1} \Lgpd{\oZ_{n-1}}, & i \geq n-1.
      \end{cases}
    \]
  \item\label{it:pi formula skseq infinite} Suppose $\nA_{\bullet} \in \SkSeqI$ corresponds
    to the chain complex $\oC = (\oC_{\bullet}, \oZ_{\bullet},
    \partial_{\bullet}) \in \CChp_{\infty}$. Then
    \[ \pi_{i} \Lgpd{\gamma_{\geq 0}\nA_{\bullet}} 
      \cong \pi_{0}\Lgpd{\oZ_{i}}/\pi_{0}\Lgpd{\oC_{i+1}}.
    \]
  \item In the situation of \ref{it:pi formula skseq infinite}, if the
    chain complex $\oC$ is object-complete, then the formula simplifies to
    \[ \pi_{i} \Lgpd{\gamma_{\geq 0}\nA_{\bullet}} 
      \cong \pi_{0}\Lgpd{\oZ_{i}}.
    \]
  \item\label{it:pi formula presp} Suppose $\spX \in \PSpIIl$ corresponds
    to the chain complex $\oC = (\oC_{\bullet}, \oZ_{\bullet},
    \partial_{\bullet}) \in \CCh$. Then
    \[ \pi_{i} L_{\Sp}(\gamma \spX) \cong 
      \pi_{0} \Lgpd{\oZ_{i}}/\pi_{0} \Lgpd{\oC_{i+1}}.
    \]
  \item In the situation of \ref{it:pi formula presp}, if the
    chain complex $\oC$ is object-complete, then the formula simplifies to
    \[ \pi_{i} L_{\Sp}(\gamma \spX) \cong 
      \pi_{0} \Lgpd{\oZ_{i}}.
    \]
  \end{enumerate}
\end{cor}
\begin{proof}
  Parts \ref{it:pi formula skseq finite}, \ref{it:pi formula skseq
    infinite}, and \ref{it:pi formula presp} hold by
  \cref{cor:complete htpy finite}, \cref{cor:skseq gpd htpy}, and
  \cref{cor:psp gpd htpy}, respectively.
  The object-complete cases follow because under this assumption
  $\pi_{0}\Lgpd{\oC_{i+1}} \to \pi_{0}\Lgpd{\oZ_{i}}$ takes every
  component of $\Lgpd{\oC_{i+1}}$ into the component of $\bbone$.
\end{proof}

Let us also spell out the complete case explicitly:
\begin{cornoproof}\label{cor:complete loc isos}\
  \begin{enumerate}[(i)]
  \item For $\nC \in \SMCatIn^{\lrig}$, we have
    \[ \pi_{i} \Lgpd{\nC} \cong
      \begin{cases}
        \pi_{0} \Lgpd{\core \mOm^{i}\nC}, & i \leq n-1, \\
        \pi_{i-n+1} \Lgpd{\mOm^{n-1}\nC}, & i \geq n-1.
      \end{cases}
      \]
    \item For $\nC \in \SMCatII^{\lrig}$, we have
      \[ \pi_{i} \Lgpd{\nC} \cong \pi_{0} \Lgpd{\core \mOm^{i}\nC}.\]
  \item For $\spX \in \SpIIl$, we have
    \[ \pi_{i} L_{\Sp}\spX \cong \pi_{0} \Lgpd{\core
        \mOm^{\infty+i}\spX}. \quad (i \in \mathbb{Z}) \qedhere\]
  \end{enumerate}
\end{cornoproof}

\begin{remark}
  Note that for a left rigid symmetric monoidal \iicat{} $\nC$, the
  isomorphisms of \cref{cor:complete loc isos} imply that we can
  fairly concretely describe $\pi_{i} \Lgpd{\nC}$: it is the quotient
  of the set $\pi_{0}\core[0] \mOm^{i}\nC$ (of equivalence classes of
  $i$-fold endomorphisms of $\bbone$) where we identify the source and
  target of any $(i+1)$-morphism between these.
\end{remark}

\begin{remark}\label{rmk:homotopy posets}
  The preceding corollary can be reformulated in terms of the
  homotopy posets of Gepner--Heine
  \cite[Definitions~4.1.1 and~4.2.4]{GepnerHeineHomotopy}. To
  distinguish their poset-valued invariant from ordinary homotopy
  groups, write
  \[
    \pi_i^{\mathrm{pos}}(\nD,z)
      := \tau_{\leq 0}\operatorname{Mor}_{\nD}^{i}(z)
  \]
  for the $i$th homotopy poset of an \iicat{} $\nD$ at an oriented
  base point $z$. For $i \geq 1$, let
  \[
    z_{\bbone,i} \colon \partial\mathbb{D}^{i} \longrightarrow \nC
  \]
  be the oriented base point determined by the tensor unit and its
  iterated identity morphisms. Then
  \[
    \operatorname{Mor}_{\nC}^{i}(z_{\bbone,i}) \simeq \mOm^{i}\nC.
  \]
  This homotopy poset is discrete. Indeed, by transposition, left
  $2$-rigidity of $\mOm^{i}\nC$ implies that for every morphism
  $x \to y$ in $\mOm^{i}\nC$ there also exists a morphism $y \to x$.

  Consequently, the preceding corollary says equivalently that the
  canonical functor $\nC \to \Lgpd{\nC}$ induces an isomorphism
  \[
    \pi_i^{\mathrm{pos}}(\nC,z_{\bbone,i})
      \xrightarrow{\;\cong\;}
    \pi_i^{\mathrm{pos}}(\Lgpd{\nC},\bar z_{\bbone,i})
      \cong \pi_i\Lgpd{\nC},
  \]
  where $\bar z_{\bbone,i}$ denotes the image of $z_{\bbone,i}$ in
  $\Lgpd{\nC}$ and the last term is regarded as a discrete poset. The
  same argument applies in degree zero.
\end{remark}

\section{Examples and applications}
\label{sec:examples}

In this section we discuss  three examples of categories constructed using categorical chain complexes. In the first two, the terms $\oZ_i$ consist informally of `closed'
objects, the terms $\oC_i$ consist of objects that are allowed to have
boundary, and the differential is literally the boundary operation.
Iterated spans give a different, but particularly explicit, complete
example.  In each case we also work out what the homotopy-group formula
of the previous section says.

\subsection{Geometric bordisms}
\label{sec:geometric-bordisms}

In this section we work without tangential structure for simplicity, but the same discussion applies, for instance, to oriented manifolds.
Let $\Bord_{n,n+1}$ denote the symmetric monoidal \icat{} of closed
$n$-manifolds and $(n+1)$-dimensional cobordisms \cite{GMTW}.
Similarly, let
$\Bord^{\partial}_{n,n+1}$ have $n$-manifolds with boundary as
objects.  Its morphisms are cobordisms with corners of codimension at
most $2$ \cite{Genauer}.  Both constructions are obtained from
topological categories.  All symmetric monoidal structures in this
subsection are given by disjoint union.

\begin{propn}\label{propn:geometric bordism chain complex}
  The boundary functors fit into pullback squares
  \[
    \begin{tikzcd}
      \Bord_{n,n+1} \drpullback \ar[r] \ar[d]
        & \Bord^{\partial}_{n,n+1} \ar[d,"\partial"] \\
      0 \ar[r] & \Bord_{n-1,n}
    \end{tikzcd}
  \]
  for $n\geq1$.  Together with
  $\oC_0=\oZ_0=\Bord_{0,1}$, these squares form an
  object-complete categorical chain complex
  \[
    \Bord_{\mathrm{ch}}\in\CChp_{\infty}.
  \]
\end{propn}

\begin{proof}
  The right vertical functor takes the boundary of a manifold and of
  a cobordism with corners.  Its fibre over the empty manifold is the
  category of closed manifolds, giving the displayed pullback square.
  Steimle proves that the boundary functor is a cocartesian fibration
  \cite[Lemma~4.4]{SteimleAdd}.  The construction is compatible with
  disjoint union.  Every object is dualizable: in the unoriented case
  its dual is the same manifold, while for an oriented manifold the
  orientation is reversed.  Evaluation and coevaluation are
  represented by the cylinder, regarded in the two directions.  Thus
  all terms in the square are rigid symmetric monoidal \icats{}.

  Finally, an $n$-manifold $M$ with boundary, regarded as an object of
  $\Bord^{\partial}_{n,n+1}$, is itself a bordism from the empty
  manifold to $\partial M$.  It therefore lies in the image of
  \[
    (\Bord_{n-1,n})_{\varnothing/}
      \longrightarrow \Bord^{\partial}_{n,n+1}
  \]
  on objects.  This proves object-completeness.
\end{proof}

\begin{remark}
  The chain complex $\Bord_{\mathrm{ch}}$ is not complete in general: there are
  invertible cobordisms that are not represented by diffeomorphisms.
  The skeletal sequence associated to $\Bord_{\mathrm{ch}}$ is expected to agree with
  the sequence
  \[
    \Bord_{0,1}\longrightarrow\Bord_{0,2}
      \longrightarrow\Bord_{0,3}\longrightarrow\cdots
  \]
  of bordism \incats{} constructed in \cite{CalaqueScheimbauer}.  We
  do not prove this comparison here.
\end{remark}

\begin{defn}\label{defn:geometric bordism unfolding}
  Let $\Bord_{0,\infty}$ denote the left rigid symmetric
  monoidal \iicat{} corresponding, under the unfolding equivalence of
  \cref{cor:unfolding infinite}, to the completion of the categorical
  chain complex $\Bord_{\mathrm{ch}}$.
\end{defn}

\begin{cor}\label{cor:geometric bordism localization}
  For every $n\geq0$ there is a natural isomorphism
  \[
    \pi_n\Lgpd{\Bord_{0,\infty}}
      \cong \pi_0\Lgpd{\Bord_{n,n+1}}.
  \]
  The group on the right consists of equivalence classes of closed
  $n$-manifolds modulo bordism.
\end{cor}

\begin{proof}
  This is the object-complete case of
  \cref{cor:pi of completion of locn}, applied to
  \cref{propn:geometric bordism chain complex}.
\end{proof}

In the unoriented case, these groups are the homotopy groups of
$\mathrm{M}O$.  Thus the corollary is compatible with Schommer-Pries's
identification of the spectrum associated to the geometric bordism
$(\infty,\infty)$-category with the unoriented cobordism spectrum
\cite{SPinv}.

\subsection{Algebraic bordisms}
\label{sec:algebraic-bordisms}

We next give an algebraic analogue of the preceding example following the
notation and constructions of the series of papers
\cites{CalmesEtAlI, CalmesEtAlII, CalmesEtAl}. We freely use the notation from the series of papers here.

Let
$(\mathcal C,\Qoppa)$ be a  Poincar\'e \icat{}, and recall that $\Qoppa^{[i]} := \Qoppa[i]$. 
Let
$
  \operatorname{Cob}(\mathcal C,\Qoppa)
$
denote the algebraic cobordism \icat{} as in  \cite[Definition~2.3.2]{CalmesEtAlII}. This  is the \icat{} associated
to the hermitian $Q$-construction.  Its objects are the Poincar\'e
objects of $(\mathcal C,\Qoppa\qshift{1})$, and its morphisms are algebraic
Poincar\'e cobordisms.
The Poincar\'e-category $\operatorname{Met}(\mathcal C,\Qoppa)$ has underlying stable \icat{}
$\Ar(\mathcal C)$ and Poincar\'e structure
\[
  \Qoppa^{\mathrm{met}}(W\to X)
    :=\operatorname{fib}\bigl(\Qoppa(X)\longrightarrow \Qoppa(W)\bigr).
\]
The target projection is part of a
Poincar\'e--Verdier sequence
\begin{equation}\label{eq:metabolic pv sequence}
  (\mathcal C,\Qoppa\qshift{-1})\longrightarrow
  \operatorname{Met}(\mathcal C,\Qoppa)\longrightarrow
  (\mathcal C,\Qoppa);
\end{equation}
see \cite[Example~1.2.5]{CalmesEtAlII}.  The algebraic cobordism
category with boundary is defined by
$ 
  \operatorname{Cob}^{\partial}(\mathcal C,\Qoppa)
    :=\operatorname{Cob}\bigl(
      \operatorname{Met}(\mathcal C,\Qoppa\qshift{1})\bigr).
$
For $k\in\mathbb Z$, we set
\[
  \operatorname{Cob}_{k,k+1}(\mathcal C,\Qoppa)
    :=\operatorname{Cob}
      (\mathcal C,\Qoppa\qshift{-k-1})
\]
and
\[
  \operatorname{Cob}^{\partial}_{k,k+1}(\mathcal C,\Qoppa)
    :=\operatorname{Cob}^{\partial}
      (\mathcal C,\Qoppa\qshift{-k-1}).
\]
Thus the objects of
$\operatorname{Cob}_{k,k+1}(\mathcal C,\Qoppa)$ are the
$\Qoppa\qshift{-k}$-Poincar\'e objects.

\begin{propn}\label{propn:algebraic bordism chain complex}
  The sequence \ref{eq:metabolic pv sequence} gives a pullback square
  \[
    \begin{tikzcd}[column sep=small]
      \operatorname{Cob}_{k,k+1}(\mathcal C,\Qoppa)
        \drpullback \ar[r] \ar[d]
        & \operatorname{Cob}^{\partial}_{k,k+1}
          (\mathcal C,\Qoppa)
          \ar[d,"\partial_k"] \\
      0 \ar[r]
        & \operatorname{Cob}_{k-1,k}(\mathcal C,\Qoppa)
    \end{tikzcd}
  \]
  for all $k\in\mathbb Z$.  The functors $\partial_k$ are
  bicartesian fibrations, and the squares form an unbounded
  categorical chain complex
  \[
    \operatorname{Cob}_{\mathrm{ch}}(\mathcal C,\Qoppa)\in\CCh.
  \]
\end{propn}

\begin{proof}
  Apply the algebraic cobordism construction to the shift of
  \cref{eq:metabolic pv sequence}
  \[
    (\mathcal C,\Qoppa\qshift{-k-1})\longrightarrow
    \operatorname{Met}(\mathcal C,\Qoppa\qshift{-k})\longrightarrow
    (\mathcal C,\Qoppa\qshift{-k}).
  \]
  By \cite[Theorem~2.5.3]{CalmesEtAlII}, the functor induced by the
  second arrow is a bicartesian fibration.  The cartesian square in
  the proof of \cite[Theorem~2.5.1]{CalmesEtAlII} identifies its fibre
  over the zero object with the cobordism category of the first term,
  namely $\operatorname{Cob}_{k,k+1}(\mathcal C,\Qoppa)$.
  This proves the pullback assertion.
  The cobordism construction preserves finite products, so the terms
  carry their natural symmetric monoidal structures induced by
  orthogonal direct sum.  The description of the cocartesian edges in
  \cite[Theorem~2.5.3]{CalmesEtAlII} is stable under direct sum; hence
  $\partial_k$ is a cocartesian fibration of symmetric monoidal
  \icats{}.

  It remains to note that these symmetric monoidal \icats{} are
  rigid.  More generally, let $(\oD,\Phi)$ be a Poincar\'e
  \icat{}.  An object of $\operatorname{Cob}(\oD,\Phi)$ is a
  Poincar\'e object
  \[
    (X,q)\in\operatorname{Pn}(\oD,\Phi\qshift{1}).
  \]
  For such an object, $(X,-q)$ is a dual.
  Evaluation is represented by the null cobordism of
  \cite[Construction~2.3.8]{CalmesEtAlII}, applied to
  $(\oD,\Phi\qshift{1})$, and coevaluation by its image under the reversal
  equivalence of \cite[Example~2.3.3(v)]{CalmesEtAlII}.  Unwinding
  the pullback used to compose $Q$-diagrams we get the triangle identities.  Applying this argument with
  \[
    (\mathcal C,\Qoppa\qshift{-k-1})
    \qquad\text{and}\qquad
    \operatorname{Met}(\mathcal C,\Qoppa\qshift{-k})
  \]
  proves rigidity of
  $\operatorname{Cob}_{k,k+1}(\mathcal C,\Qoppa)$ and
  $\operatorname{Cob}^{\partial}_{k,k+1}
    (\mathcal C,\Qoppa)$.
\end{proof}

The algebraic cobordism categories with boundary admit an alternative
description.  The
algebraic Thom construction gives natural equivalences
\begin{equation}\label{eq:algebraic boundary span description}
  \operatorname{Cob}^{\partial}_{k,k+1}
    (\mathcal C,\Qoppa)\simeq
  \Span\bigl(
    \operatorname{He}(\mathcal C,\Qoppa\qshift{-k-1})\bigr)
\end{equation}
by \cite[Proposition~2.3.12]{CalmesEtAlII}.  Thus an object of
$\operatorname{Cob}^{\partial}_{k,k+1}
  (\mathcal C,\Qoppa)$ can be viewed as a possibly degenerate
$\Qoppa\qshift{-k-1}$-hermitian
object.  Under this description, $\partial_k$ is the functor induced
by the metabolic target projection.  For
example, taking $\mathcal C=\operatorname{Perf}(R)$ with a standard
Poincar\'e structure associated to a ring with involution recovers
the familiar algebraic bordism categories of non-degenerate perfect
complexes over $R$. 

\begin{defn}\label{defn:algebraic bordism unfolding}
  Let
  \[
    \operatorname{Cob}_{-\infty,\infty}
      (\mathcal C,\Qoppa)\in\SpIIl
  \]
  denote the completion of the $(\infty,\infty)$-prespectrum
  corresponding to
  $\operatorname{Cob}_{\mathrm{ch}}(\mathcal C,\Qoppa)$ under the
  unfolding equivalence of \cref{cor:unfolding unbdd}.
\end{defn}

\begin{cor}\label{cor:algebraic bordism homotopy}
  For every $i\in\mathbb Z$ there are natural isomorphisms
  \[
    \pi_i L_{\Sp}\bigl(
      \operatorname{Cob}_{-\infty,\infty}
        (\mathcal C,\Qoppa)\bigr)
      \cong \pi_0\Lgpd{\operatorname{Cob}_{i,i+1}
        (\mathcal C,\Qoppa)}
      \cong \mathrm{L}_i(\mathcal C,\Qoppa).
  \]
\end{cor}

\begin{proof}
  By \cite[Corollary~2.3.10]{CalmesEtAlII}, the cobordism category of
  a metabolic Poincar\'e \icat{} is connected.  Hence
  \[
    \pi_0\Lgpd{\operatorname{Cob}^{\partial}_{i+1,i+2}
      (\mathcal C,\Qoppa)}
      =\pi_0\Lgpd{\operatorname{Cob}\bigl(
        \operatorname{Met}(\mathcal C,
          \Qoppa\qshift{-i-1})\bigr)}=0.
  \]
  The first isomorphism now follows from
  \cref{cor:pi of completion of locn}.  The second is the
  cobordism description of the zeroth $L$-group
  \cite[Proposition~2.3.7 and Definition~4.2.6]{CalmesEtAlII}, applied
  after shifting $\Qoppa$.
\end{proof}

\begin{propn}\label{propn:algebraic bordism realization}
  There is a natural equivalence of spectra
  \[
    L_{\Sp}\operatorname{Cob}_{-\infty,\infty}
      (\mathcal C,\Qoppa)
      \simeq
    \mathrm{L}(\mathcal C,\Qoppa),
  \]
  where $\mathrm{L}(\mathcal C,\Qoppa)$ denotes the $L$-theory
  spectrum of the Poincar\'e \icat{} $(\mathcal C,\Qoppa)$.
\end{propn}

\begin{proof}
  Regarding a $\Qoppa$-Poincar\'e object as an object of
  $\operatorname{Cob}_{0,1}(\mathcal C,\Qoppa)$ defines a
  natural transformation
  \[
    \operatorname{Pn}(\mathcal C,\Qoppa)
      \longrightarrow
    \Omega^\infty L_{\Sp}
      \operatorname{Cob}_{-\infty,\infty}
        (\mathcal C,\Qoppa).
  \]
  The functor
  $
    (\mathcal C,\Qoppa)\longmapsto
    L_{\Sp}\operatorname{Cob}_{-\infty,\infty}
      (\mathcal C,\Qoppa)
  $
  is additive and bordism invariant.  Indeed, by
  \cref{cor:algebraic bordism homotopy}, the maps expressing these
  properties induce on every homotopy group the corresponding maps on 
  $
    \mathrm{L}_i(\mathcal C,\Qoppa),
  $
  which are isomorphisms.
  By the universal property of $L$-theory
  \cite[Theorem~4.4.11]{CalmesEtAlII}, there is therefore a natural
  transformation
  \[
    \vartheta\colon
    \mathrm{L}(\mathcal C,\Qoppa)
      \longrightarrow
    L_{\Sp}\operatorname{Cob}_{-\infty,\infty}
      (\mathcal C,\Qoppa).
  \]
  Under the canonical identifications
  $
    \pi_i\mathrm{L}(\mathcal C,\Qoppa)
      \cong
    \mathrm{L}_0(\mathcal C,\Qoppa\qshift{-i})
  $
  of \cite[Corollary~4.4.6]{CalmesEtAlII} and
  $
    \pi_iL_{\Sp}\operatorname{Cob}_{-\infty,\infty}
      (\mathcal C,\Qoppa)
      \cong
    \mathrm{L}_0(\mathcal C,\Qoppa\qshift{-i})
  $
  of \cref{cor:algebraic bordism homotopy}, the map induced by
  $\vartheta$ is the identity: both transformations send a
  Poincar\'e object to its algebraic cobordism class.  Hence
  $\vartheta$ induces an isomorphism on every homotopy group and is
  therefore an equivalence.
\end{proof}

\begin{remark}
  There are also natural levelwise equivalences
  \[
    \mathrm{GW}(\mathcal C,\Qoppa\qshift{-k-1})
      \simeq \Omega\Lgpd{\operatorname{Cob}_{k,k+1}
        (\mathcal C,\Qoppa)}
  \]
  by \cite[Corollary~4.1.5]{CalmesEtAlII}, where $\mathrm{GW}$ denotes
  the Grothendieck--Witt space of
  \cite[Definition~4.1.1]{CalmesEtAlII}. 
\end{remark}

\subsection{Equivariant bordisms and automorphisms}
\label{sec:equivariant-bordisms}

Let $\Glo$ be the full $(2,1)$-subcategory of groupoids spanned by the
one-object groupoids $BG$ for discrete groups $G$.  We use $\Glo$ for
this discrete-group version throughout this subsection.  Thus
the mapping groupoid from $BH$ to $BG$ has as objects the homomorphisms
$H\to G$ and as morphisms conjugations in $G$.

As in \cref{sec:geometric-bordisms}, we work with unoriented
manifolds. The same construction applies verbatim after imposing any
compatible family of tangential structures. For a discrete group $G$, let
\[
  \Bord^G_{n,n+1}
  \quad\text{and}\quad
  \Bord^{\partial,G}_{n,n+1}
\]
denote the equivariant analogues of the categories in
\cref{sec:geometric-bordisms}. Here a $G$-action is an actual smooth
action on the manifold. It is not an action introduced only after passing
to the bordism \icat{} or an action on the underlying homotopy types of
the manifolds. In the presence of boundary and corners we include
invariant product collars in the data.

\begin{propn}\label{propn:equivariant bordism functoriality}
  For every discrete group $G$ and every $n\geq1$, taking boundaries
  gives a pullback square
  \[
    \begin{tikzcd}[column sep=small]
      \Bord^G_{n,n+1}
        \drpullback \ar[r] \ar[d]
        & \Bord^{\partial,G}_{n,n+1}
          \ar[d,"\partial_n"] \\
      0 \ar[r]
        & \Bord^G_{n-1,n}.
    \end{tikzcd}
  \]
  Together with $\Bord^G_{0,1}$ in degree zero, these squares
  form an object-complete categorical chain complex
  $\Bord_{\mathrm{ch}}^G\in\CChp_{\infty}$.  As $G$ varies,
  restriction of actions assembles these chain complexes into a
  contravariant functor
  \[
    \Glo^{\op}\longrightarrow\CChp_{\infty}.
  \]
\end{propn}

\begin{proof}
  The proof of \cref{propn:geometric bordism chain complex} is
  equivariant: boundary gives the pullback squares, the equivariant
  cylinder gives duals, and a $G$-manifold with boundary is a
  $G$-bordism from the empty manifold to its boundary.  The chosen
  collars ensure that all gluings are equivariant.

  Restriction along a homomorphism $\varphi\colon H\to G$ gives a
  functor $\varphi^*$ from $G$-manifolds to $H$-manifolds.  If
  $\psi=c_g\varphi$ for $g\in G$, then the action of $g$ on every
  $G$-manifold defines a monoidal natural equivalence
  $\varphi^*\simeq\psi^*$.  These are precisely the $2$-morphisms in
  $\Glo$, which proves the claimed $(2,1)$-functoriality.
\end{proof}

We let
$
  \Bord_{\mathrm{unf}}^G\in\PSpIIl
$
be the prespectrum obtained by extending $\Bord_{\mathrm{ch}}^G$ by zero in
negative degrees and unfolding it.  Its completed degree-zero unfolding
is denoted $\Bord^{G}_{0,\infty}$, in accordance
with the notation $\Bord_{0,\infty}$ of
\cref{sec:geometric-bordisms}.  Define
$
  E_{\Bord}^G
    :=L_{\Sp}\bigl(\gamma
      \Bord_{\mathrm{unf}}^G\bigr)\in\Sp.
$
Equivalently, $E_{\Bord}$ is the composite
\[
  \Glo^{\op}\longrightarrow\CChp_{\infty}
    \longrightarrow\PSpIIl\xto{\gamma}\SpIIl
    \xto{L_{\Sp}}\Sp.
\]
This functoriality gives a natural coassembly map
\begin{equation}\label{eq:bordism coassembly}
  \operatorname{coass}_G\colon E_{\Bord}^G\longrightarrow
    \bigl(E_{\Bord}^1\bigr)^{BG}.
\end{equation}
By \cref{cor:pi of completion of locn}, object-completeness gives
\begin{equation}\label{eq:equivariant bordism homotopy groups}
  \pi_nE_{\Bord}^G
    \cong\pi_0\Lgpd{\Bord^G_{n,n+1}}
    \quad(n\geq0),
  \qquad
  \pi_nE_{\Bord}^G=0\quad(n<0).
\end{equation}

For a finite group $G$, Galatius--Sz\H{u}cs \cite[Theorem~1.1]{GalatiusSzucs}
construct a genuine orthogonal $G$-spectrum $\mathrm{MTO}_d$ and prove
for every $n$ a weak equivalence
$
  \Lgpd{\Bord^G_{n,n+1}}
    \simeq\bigl(\Omega^{\infty-1}\mathrm{MTO}_{n+1}\bigr)^G . 
$
The right-hand side uses actual
$G$-fixed points.  Their \cite[Theorem~6.4]{GalatiusSzucs} gives the
corresponding statement with tangential structures.  In particular,
\cref{eq:equivariant bordism homotopy groups} recovers the usual
equivariant unoriented bordism group
$\pi_n^G(\mathrm{mO})$.  

The construction above has an oriented variant, obtained by requiring all
actions and bordisms to preserve orientations; we denote the resulting spectra
by $E_{\Bord}^{G,\mathrm{SO}}$. For $G=\mathbb Z$, such an action on a
compact manifold is exactly an orientation-preserving diffeomorphism
$f\colon M\to M$.
After the standard collar-straightening, equivariant bordism agrees
with the relation defined by pairs $(W,F)$ in which $F$ is an
orientation-preserving diffeomorphism extending the boundary
diffeomorphisms. Applying the oriented analogue of
\cref{eq:equivariant bordism homotopy groups} gives the exact identification
\begin{equation}\label{eq:Z bordism Kreck}
  \pi_nE_{\Bord}^{\mathbb Z,\mathrm{SO}}
    \cong\Delta_n\qquad(n\geq0),
\end{equation}
where $\Delta_n$ is Kreck's oriented bordism group of
orientation-preserving diffeomorphisms of closed oriented
$n$-manifolds \cite{KreckDiffeomorphisms}.  Restriction along
$\mathbb Z\to\mathbb Z$, $1\mapsto k$, induces
$(M,f)\mapsto(M,f^k)$ on these groups.

The oriented analogue of \cref{eq:bordism coassembly} gives a map
\[
  \operatorname{coass}^{\mathrm{SO}}_{\mathbb Z}\colon
  E_{\Bord}^{\mathbb Z,\mathrm{SO}}\longrightarrow
    \bigl(E_{\Bord}^{1,\mathrm{SO}}\bigr)^{B\mathbb Z}
\]
with a concrete description. Choose $B\mathbb Z\simeq S^1$ using the
positive generator of $\mathbb Z$.
The stable splitting
$\Sigma^\infty_+S^1\simeq\mathbb S\oplus \Sigma\mathbb S$ identifies
\[
  \pi_n\bigl(E_{\Bord}^{1,\mathrm{SO}}\bigr)^{B\mathbb Z}
    \cong\Omega_n^{\mathrm{SO}}\oplus
      \Omega_{n+1}^{\mathrm{SO}}
    \cong\Omega_{n+1}^{\mathrm{SO}}(S^1).
\]
Under \cref{eq:Z bordism Kreck}, the induced map on homotopy groups is
\begin{equation}\label{eq:Z coassembly mapping torus}
  \pi_n(\operatorname{coass}^{\mathrm{SO}}_{\mathbb Z})\colon
  \Delta_n\longrightarrow
  \Omega_n^{\mathrm{SO}}\oplus\Omega_{n+1}^{\mathrm{SO}},
  \qquad
  [M,f]\longmapsto\bigl([M],[T_f]\bigr),
\end{equation}
where
$
  T_f:=M\times[0,1]/(x,1)\sim(f(x),0)
$
is the mapping torus.  Before applying the stable splitting, the
corresponding loop in the bordism category is represented by the
mapping cylinder of $f$, regarded as a bordism from $M$ to itself;
closing its incoming and outgoing boundaries gives $T_f$.  More
precisely, the reduced coordinate is
$[T_f]-[M\times S^1]=[T_f]$, since $M\times S^1$ bounds
$M\times D^2$.  Equivalently,
the second description above sends $[M,f]$ to the mapping-torus class
$[T_f\to S^1]\in\Omega_{n+1}^{\mathrm{SO}}(S^1)$
\cite[Definition~28.8, Remark~28.9, and
Proposition~28.10(ii)]{RanickiKnot}. \\

There is a parallel algebraic construction.  For a Poincar\'e
\icat{} $(\mathcal C,\Qoppa)$ and a discrete group $G$, write
\[
  \mathcal C^{BG}:=\Fun(BG,\mathcal C),
  \qquad
  \Qoppa_G(X):=\lim_{BG^{\op}}(\Qoppa\circ X^{\op})
    \simeq\Qoppa(X)^{hG}.
\]
The inversion equivalence $BG^{\op}\simeq BG$ and the duality on
$\mathcal C$ give $\mathcal C^{BG}$ a pointwise duality $D_G$.  The
polarization of $\Qoppa_G$ satisfies
\[
  B_{\Qoppa_G}(X,Y)\simeq B_{\Qoppa}(X,Y)^{hG}
    \simeq\Map_{\mathcal C^{BG}}(X,D_GY),
\]
so $\Qoppa_G$ is a Poincar\'e structure.  Its Poincar\'e objects are
Poincar\'e objects of $(\mathcal C,\Qoppa)$ equipped with a coherent
$G$-action by isometries.  Restriction along $H\to G$ defines a
Poincar\'e functor
$
  (\mathcal C^{BG},\Qoppa_G)\longrightarrow
  (\mathcal C^{BH},\Qoppa_H),
$
so the algebraic bordism chain complexes are contravariantly
functorial in $BG$ as well.
For $G=\mathbb Z$, define
\begin{equation*}
  L^{\mathrm{Aut}}(\mathcal C,\Qoppa)
    :=L_{\Sp}\bigl(
      \operatorname{Cob}_{-\infty,\infty}
      (\mathcal C^{B\mathbb Z},\Qoppa_{\mathbb Z})\bigr).
\end{equation*}
By \cref{propn:algebraic bordism realization}, this spectrum is
canonically equivalent to
$\mathrm L(\mathcal C^{B\mathbb Z},\Qoppa_{\mathbb Z})$.
Moreover, \cref{cor:algebraic bordism homotopy} gives
$
  \pi_nL^{\mathrm{Aut}}(\mathcal C,\Qoppa)
    \cong\mathrm L_n\bigl(
      \mathcal C^{B\mathbb Z},\Qoppa_{\mathbb Z}\bigr).
$
For $\mathcal C=\operatorname{Perf}(A)$ with its symmetric Poincar\'e
structure, the standard comparison with classical symmetric
Poincar\'e complexes identifies the right-hand side, by unwinding
objects and cobordisms, with Ranicki's projective symmetric
automorphism group $\mathrm{LAut}^{n}_{p}(A)$
\cite[Definition~28.16 and Terminology~28.18]{RanickiKnot}.  Over
$A=\mathbb Z$ the projective, free, and simple decorations agree since
$\widetilde K_0(\mathbb Z)=\mathrm{Wh}_1(\mathbb Z)=0$.  We then write
simply $\mathrm{LAut}^{n}(\mathbb Z)$.

At the level of bordism groups, Ranicki's simple automorphism symmetric
signature is
\begin{equation}\label{eq:Ranicki automorphism signature}
  \sigma^{*}_{\mathrm{aut}}\colon\Delta_n\longrightarrow
  \mathrm{LAut}^{n}(\mathbb Z),\qquad
  [M,f]\longmapsto
  [C_*(M;\mathbb Z),\varphi_M,(f_*,\chi_f)],
\end{equation}
where $\varphi_M$ is the symmetric Poincar\'e structure and $\chi_f$
is the  homotopy expressing that $f_*$ preserves it
\cite[Definition~28.20 and Proposition~28.22]{RanickiKnot}.  In
dimension $n=2k$ this is the invariant
\[
  [H_k(M;\mathbb Z)/\mathrm{tors},\lambda_M,f_*]
    \in W_{(-1)^k}(\mathbb Z;\mathbb Z)
\]
of L\'opez de Medrano and Kreck
\cite[Example~28.23]{RanickiKnot}.  Thus the relevant algebraic target
is symmetric automorphism $L$-theory. This is different from ordinary
$L$-theory of $\mathbb Z[z,z^{-1}]$:
Ranicki identifies automorphism $L$-groups with relative terms
for a  localization and places
\cref{eq:Ranicki automorphism signature}, together with the mapping
torus class, in a localization exact sequence
\cite[Propositions~28.17 and~28.22]{RanickiKnot}.  

\begin{remark}
  The functor
  $
    E_{\Bord}\colon\Glo^{\op}\longrightarrow\Sp
  $
  constructed above only records restriction of actual group actions.
  It is not a genuine global equivariant spectrum.  Such a refinement
  would have to add induction or transfer maps.
  For arbitrary discrete groups induction also fails to preserve
  compactness unless the relevant subgroup has finite index.  \end{remark}

\subsection{Iterated spans}
\label{sec:iterated-spans}

Let $\oE$ be an \icat{} with finite limits.  For $n\geq1$, let
$\Span_n(\oE)$ denote Haugseng's $(\infty,n)$-category of iterated
spans in $\oE$ \cite{HaugsengSpans}.  Its objects are those of
$\oE$, its $1$-morphisms are spans, and its higher morphisms are
iterated spans.  There are natural equivalences
\begin{equation}\label{eq:span mapping category}
  \Span_n(\oE)(X,Y)
    \simeq \Span_{n-1}(\oE_{/X\times Y})
\end{equation}
by \cite[Proposition~8.3]{HaugsengSpans}; here we use the convention
$\Span_0(\oD):=\core[0]\oD$.  Cartesian product in $\oE$ induces a
symmetric monoidal structure on $\Span_n(\oE)$.

\begin{propn}\label{propn:iterated spans rigid}
  For every $n\geq1$, $\Span_n(\oE)$ is fully dualizable, and in
  particular left rigid.  More precisely, every $k$-morphism for
  $1\leq k<n$ has both a left and a right adjoint, and every object is
  dualizable.
\end{propn}

\begin{proof}
  Haugseng proves the assertion about adjoints in
  \cite[Theorem~12.4]{HaugsengSpans} and the assertion about duals in
  \cite[Corollary~12.5]{HaugsengSpans}.  Concretely, the adjoint of a
  span is obtained by reversing it.  Every object $X$ is self-dual,
  with evaluation and coevaluation represented by
 $
    X\times X\xleftarrow{\Delta}X\longrightarrow *$ and $
    *\longleftarrow X\xrightarrow{\Delta}X\times X$
  respectively.
\end{proof}

Inserting a constant final simplicial direction in Haugseng's
iterated Segal-space model gives strong symmetric monoidal functors
\[
  j_n\colon\Span_n(\oE)\longrightarrow\Span_{n+1}(\oE).
\]

\begin{lemma}\label{lem:span core tower}
  For every $n\geq1$, the functor $j_n$ induces an equivalence
  \[
    \Span_n(\oE)\xrightarrow{\;\sim\;}
      \core[n]\Span_{n+1}(\oE).
  \]
\end{lemma}

\begin{proof}
  We argue by induction on $n$.  The two higher categories have the
  same objects.  For $n=1$, \cref{eq:span mapping category} identifies
  the mapping space in $\Span_1(\oE)$ with
  $
    \core[0](\oE_{/X\times Y}).
  $
  The corresponding mapping space in
  $\core[1]\Span_2(\oE)$ is
  $
    \core[0]\Span_1(\oE_{/X\times Y}).
  $
  A morphism in a span category is an equivalence precisely when both
  legs of its representing span are equivalences
  \cite[Lemma~8.2]{HaugsengSpans}.  Consequently
  $\core[0]\Span_1(\oD)\simeq\core[0]\oD$, which proves the base
  case.

  Under \cref{eq:span mapping category}, the functor induced by $j_n$
  on the mapping categories from $X$ to $Y$ identifies with
  $j_{n-1}$ for the slice $\oE_{/X\times Y}$.  For the inductive step,
  the target of this mapping functor is therefore
  \[
    \begin{split}
      \bigl(\core[n]\Span_{n+1}(\oE)\bigr)(X,Y)
        &\simeq \core[n-1]\Span_n(\oE_{/X\times Y}) \\
        &\simeq \Span_{n-1}(\oE_{/X\times Y}),
    \end{split}
  \]
  and the second equivalence is the inductive hypothesis.  Its source
  is the mapping category in $\Span_n(\oE)$ by
  \cref{eq:span mapping category}.  Hence $j_n$ is fully faithful,
  and it is essentially surjective since it is the identity on
  objects.  The construction is compatible with Cartesian products.
\end{proof}

\begin{cor}\label{cor:span skeletal sequence}
  The sequence
  \[
    \Span_1(\oE)\longrightarrow\Span_2(\oE)
      \longrightarrow\cdots
  \]
  is a complete infinite skeletal sequence.  It therefore determines
  a left rigid symmetric monoidal \iicat{}, which we denote by
  $\Span_{\infty}(\oE)$.
\end{cor}

\begin{proof}
  Left rigidity follows from
  \cref{propn:iterated spans rigid}.  By
  \cref{lem:span core tower}, each transition identifies its source
  with the appropriate core of its target, which is precisely the
  completeness condition.
\end{proof}

Write $(\oC_{\bullet},\oZ_{\bullet},\partial_{\bullet})$ for the
categorical chain complex obtained by unfolding this skeletal
sequence. We can identify this quite explicitly using the following identification of $\mP\Span_{2}(\oE)$, which follows from recent work of G\"odicke, Ho, and Stern in \cite{SpanII}:
\begin{propn}\label{propn:oplax slice spans}
  We have a natural equivalence of cocartesian fibrations of symmetric monoidal \icats{}
  \[
    \begin{tikzcd}
     \mP\Span_{2}(\oE) \ar[rr, "\sim"] \ar[dr, "\partial"'] & & \Span_{1}(\oE^{[1]}) \ar[dl, "\ev_{1}"] \\
      & \Span_{1}(\oE).
    \end{tikzcd}
  \]
\end{propn}
\begin{proof}
  By definition, the left-hand side is given by a pullback
  \[
    \begin{tikzcd}
      \mP\Span_{2}(\oE) \ar[r] \ar[d, "\partial"] & \core \FUN([1], \Span_{2}(\oE))^{\oplax}  \ar[d, "{(\ev_{0}, \ev_{1})}"] \\
      \{*\} \times \Span_{1}(\oE) \ar[r] & \Span_{1}(\oE) \times \Span_{1}(\oE),
    \end{tikzcd}
  \]
  where $* \in \oE$ is the terminal object (the monoidal unit).
  As a special case of \cite{SpanII}*{Theorem A}, we have a natural
  equivalence of \icats{}
  \[ \core \FUN([1], \Span_{2}(\oE))^{\oplax} \simeq
    \Span_{1}(\Fun(\Lambda, \oE)),\] where $\Lambda$ is the category
  $\ell \from t \to r$, \ie{} the barycentric subdivision of $[1]$,
  with evaluation at $0$ and $1$ induced by evaluation at $\ell$ and
  $r$; since both sides preserve products as functors of $\oE$, this
  also gives an equivalence of symmetric monoidal \icats{}. Thus the
  pullback square defining $\mP\Span_{2}(\oE)$ is obtained by applying
  $\Span_{1}(\blank)$ to the pullback square
  \[
    \begin{tikzcd}
     \oE^{[1]} \ar[r] \ar[d, "\ev_{1}"'] & \Fun(\Lambda, \oE)  \ar[d, "{(\ev_{\ell},\ev_{r})}"] \\
     \{*\} \times \oE \ar[r] & \oE \times \oE,
    \end{tikzcd}
  \]
  where we use that $\Lambda \simeq [1] \amalg_{[0]} [1]$ and that the pullback of $\ev_{1} \colon \oE^{[1]} \to \oE$ to $\{*\}$ gives $\oE_{/*} \simeq \oE$ as $*$ is the terminal object.
\end{proof}

\begin{cor}\label{propn:span chain complex}
  There are natural symmetric monoidal equivalences
  \[
    \oZ_i\simeq\Span_1(\oE)\quad(i\geq0),
    \qquad
    \oC_i \simeq \mP\Span_{2}(\oE) \simeq \Span_1(\oE^{[1]})\quad(i\geq1),
  \]
  and $\oC_0\simeq\oZ_0$.  Under the second equivalence the boundary
  functor is induced by target evaluation
  $\ev_1\colon\oE^{[1]}\to\oE$.  Thus every non-degenerate
  square in the unfolded complex is naturally equivalent to
  \[
    \begin{tikzcd}
      \Span_1(\oE) \drpullback \ar[r] \ar[d]
        & \Span_1(\oE^{[1]})
          \ar[d,"\Span_1(\ev_1)"] \\
      0 \ar[r] & \Span_1(\oE).
    \end{tikzcd}
  \]
  This categorical chain complex is complete.
\end{cor}

\begin{proof}
  Let $*\in\oE$ be terminal, hence the tensor unit in the span
  categories.  Applying \cref{eq:span mapping category} at $(*,*)$
  gives
  \[
    \mOm\Span_{r+1}(\oE)
      \simeq\Span_r(\oE_{/*\times *})
      \simeq\Span_r(\oE).
  \]
  Iterating this equivalence and applying
  \cref{cor:skseq cch eq} gives the first description of the terms, and the identification with spans in the arrow \icat{} is then \cref{propn:oplax slice spans}.
  Finally, completeness follows from
  \cref{cor:span skeletal sequence}, since the unfolding of a complete
  skeletal sequence is complete.
\end{proof}

\begin{remark}\label{rmk:span boundary fibration}
  The description by the arrow category makes the boundary fibration
  quite concrete.  The functor $\partial=\Span_1(\ev_1)$ sends an
  arrow $S\to X$ to its target $X$, and its fibre over $X$ is
  \[
    \Span_2(\oE)(*,X)\simeq\Span_1(\oE_{/X}).
  \]
  Cocartesian transport along a span
  $X\xleftarrow{}M\xrightarrow{}Y$ is composition of spans. It sends
  $S\to X$ to
  $
    S\times_XM\longrightarrow M\longrightarrow Y.
  $
  Thus every degree of the complex is governed by the single family
  $
    X\longmapsto\Span_1(\oE_{/X}).
  $
\end{remark}

\begin{cor}\label{cor:span localization}
  For $n\geq1$ there are natural isomorphisms
  \[
    \pi_i\Lgpd{\Span_n(\oE)}\cong
    \begin{cases}
      0, & 0\leq i<n,\\
      \pi_{i-n+1}\Lgpd{\Span_1(\oE)}, & i\geq n.
    \end{cases}
  \]
  Moreover $
    \Lgpd{\Span_{\infty}(\oE)}\simeq *.$
\end{cor}

\begin{proof}
  For any $X,Y\in\oE$, the product gives a span
  \[
    X\longleftarrow X\times Y\longrightarrow Y.
  \]
  Consequently, the localization of $\Span_1(\oE)$ is connected, so
  $\pi_0\Lgpd{\Span_1(\oE)}=0$.  The finite-dimensional formula now
  follows from \cref{cor:pi of completion of locn} and
  \cref{propn:span chain complex}.  Applying the infinite complete
  case of the same corollary shows that all homotopy groups of
  $\Lgpd{\Span_{\infty}(\oE)}$ vanish.  This space is nonempty and
  connected, hence contractible.
\end{proof}

\begin{remark}
  The last assertion does not say that
  $\Lgpd{\Span_1(\oE)}$ is contractible.  For example, if $\oE$ is a
  small stable \icat{}, the usual $Q$-construction gives
  \[
    \Omega^\infty K(\oE)\simeq\Omega\Lgpd{\Span_1(\oE)};
  \]
  see \cite[Example~2.3.3(iii) and Observation~2.3.11]{CalmesEtAlII}.
  In this case \cref{cor:span localization} reads
  \[
    \pi_i\Lgpd{\Span_n(\oE)}\cong
    \begin{cases}
      0, & 0\leq i<n,\\
      K_{i-n}(\oE), & i\geq n.
    \end{cases}
  \]
  Thus each finite stage retains a shifted copy of the algebraic
  $K$-groups, whereas each fixed homotopy group disappears in the
  infinitely iterated case.
\end{remark}

\end{document}